\magnification=\magstep1
\input amstex
\documentstyle{amsppt}

\define\defeq{\overset{\text{def}}\to=}

\define\Gal{\operatorname{Gal}}

\def \isom {\overset \sim \to \rightarrow}

\define\Edg{\operatorname{Edg}}
    
\define\Spec{\operatorname{Spec}}
\define\id{\operatorname{id}}
\define\Ver{\operatorname{Ver}}

\define\nor{\operatorname{nor}}
\define\Ker{\operatorname{Ker}}
\define\op{\operatorname{op}}

\def \cd{\operatorname {cd}}
\def \Ker{\operatorname {Ker}}
\def \id{\operatorname {id}}
\def \Fr{\operatorname {Fr}}
\def \char{\operatorname {char}}
\def \Spec{\operatorname {Spec}}
\def \Spf{\operatorname {Spf}}
\def \Sp{\operatorname {Sp}}
\def \geo{\operatorname {geo}}
\def \Exp{\operatorname {Exp}}
\def \Spf{\operatorname {Spf}}

\NoRunningHeads
\NoBlackBoxes
\topmatter

\title 
Fake Liftings of Galois Covers between Smooth Curves
\endtitle

\author
Mohamed Sa\"\i di
\endauthor

\abstract
In this paper we investigate the problem of lifting of Galois covers between algebraic curves from characteristic $p>0$ to characteristic
$0$.  We prove a refined version of the main result of Garuti concerning this problem in [Ga]. We formulate a refined version
of the Oort conjecture on liftings of cyclic Galois covers between curves. We introduce the notion of fake liftings of cyclic Galois covers between curves,
their existence would contradict the Oort conjecture, and we study the geometry of their semi-stable models. 
Finally, we introduce and investigate on some examples
the smoothening process, which ultimately aims to show that fake liftings do not exist. This in turn would imply the Oort conjecture.
\endabstract

\toc

\subhead
\S 0. Introduction
\endsubhead

\subhead
\S 1. Background 
\endsubhead

\subhead
\S 2. Pro-p Quotients of the Geometric Galois Group of a $p$-adic Open Disc
\endsubhead

\subhead
\S 3. Fake Liftings of Cyclic Covers between Smooth Curves
\endsubhead

\subhead
\S 4. The Smoothening Process
\endsubhead

\endtoc

\endtopmatter

\document 

\subhead
\S 0. Introduction
\endsubhead
In this paper we investigate the following problem, known as the problem of lifting of Galois covers between algebraic curves from
characteristic $p>0$ to characteristic $0$.

In what follows $R$ will denote a complete discrete valuation ring 
of unequal characteristics, $K\defeq \Fr (R)$ the quotient field of $R$, $\char (K)=0$,
and $k$ the residue field of $R$, which we assume to be algebraically 
closed of characteristic $p>0$.

\definition {Problem I} Let 
$$f_k:Y_k\to X_k$$
be a finite Galois cover between smooth $k$-curves, with Galois group $G$. Is it possible to lift the Galois 
cover $f_k$ to a Galois cover
$$f:Y'\to X'$$
between smooth $R'$-curves, where $R'/R$ is a finite extension?
\enddefinition

In the original version of this problem, one doesn't fix $R$, but fixes $k$, $f_k:Y_k\to X_k$, and 
asks for the existence of a local domain $R$ dominating the ring of Witt vectors $W(k)$ 
over which a lifting of $f_k$ exists as part of the problem (cf. [Oo]).

One can formulate Problem I in terms of lifting of curves and their automorphism groups from positive to zero characteristics (cf. [Oo], [Oo1]).
One can also formulate the following variant of the above problem, where one fixes a lifting of the curve $X_k$.

\definition{Problem II} Let $X$ be a proper, smooth, geometrically connected $R$-curve, and 
$$f_k:Y_k\to X_k\defeq X\times _Rk$$ 
a finite Galois cover between smooth 
$k$-curves, with Galois group $G$. Is it possible to lift the Galois cover $f_k$ to a Galois cover 
$$f:Y\to X'\defeq X\times _RR'$$ 
where $R'/R$ is a finite extension, and $Y$ is a smooth $R'$-curve?
\enddefinition

We shall refer to a lifting $f$ as above, if it exits, as a smooth lifting of the Galois cover $f_k$.

The above two problems are in fact equivalent (cf. discussion in 3.1), and one may well consider Problem II, instead.

This problem has been considered successfully by Grothendieck in the case where $f_k$ is a tamely ramified cover.
In this case a smooth lifting $f$ as above exists over $R$ (cf. [Gr]). 

The answer to this problem is however No in general.  Indeed, in the case where $G$ is the full automorphism group of $Y_k$, there are examples where
the size of $G$ exceeds the Hurwitz bound for the size of automorphism groups of curves in characteristic zero (cf. [Ro]), and the cover $f_k$ can not be lifted 
in this case. Also it is in general necessary to perform a finite extension of $R$ in order to solve this problem (cf. [Oo], 1). 

In the case where $f_k$ is wildly ramified there are non liftable examples with Galois groups as simple as $G\isom \Bbb Z/p\Bbb Z\times
\Bbb Z/p\Bbb Z$ (cf. [Gr-Ma], 5).  See also [Oo], 1, for an example of a genus $2$ curve in characteristic $5$, 
and an automorphism group of cardinality $20$, which cannot lift to characteristic $0$. 

The most general result one can hope for, in the case where $f_k$ is wildly ramified,  is the following which was conjectured by F. Oort.

\definition{Oort conjecture [Conj-O]} Problem I, or equivalently Problem II,  has a positive answer if $G\isom \Bbb Z/m\Bbb Z$ is a cyclic group. 
Moreover, in this case one can choose $R'$ in Problem I, and Problem II, to be the minimal
extension of $R$ which contains the $m$-th roots of $1$.
\enddefinition

In order to solve this conjecture one may reduce to the case where $G\isom \Bbb Z/p^n\Bbb Z$ is a cyclic $p$-group (cf. Lemma 3.1.1).
In this case Oort conjecture has been verified when $n\le 2$ (cf. [Se-Oo-Su] for the case $n=1$, and [Gr-Ma] for the case $n=2$). 
In the approach of Oort, Sekiguchi, Suwa,
Green, and Matignon one uses the Oort-Sekiguchi-Suwa theory, which provides [explicit] equations describing the degeneration of the 
Kummer equations in characteristic $0$ to the Artin-Schreier-Witt
equations in characteristic $p>0$. 
In general, this conjecture is still widely open. To the best knowledge of the author, no concrete liftable 
examples are known in the case where $n\ge 3$. 

However, if in Problem II one relaxes the requirement that $Y$ in the lifting $f:Y\to X'\defeq X\times _RR'$ is smooth over $R'$, one has the following rather 
general result where one allows introducing singularities in $Y_k$, and which is du to Garuti (cf. [Ga], and Theorem 2.5.1). 

\proclaim {Theorem A (Garuti)}
There exists a finite extension $R'/R$, and a finite Galois cover
$f':Y'\to X'\defeq X\times _RR'$, with Galois group $G$, where $Y'$ is a normal $R'$-curve (which in general need not be smooth over $R'$), 
and the natural morphism $f'_k:Y'_k\defeq Y'\times _{R'}k\to X_k$ between special fibres is generically Galois with Galois group $G$. Moreover,  there exists a factorisation
$f_k:Y_k@>\nu>> Y'_k@>{f'_k}>> X_k$, where the morphism $\nu:Y_k\to Y'_k$ is a morphism of normalisation, which is an isomorphism outside the ramified points, and
$Y'_k$ is unibranch.
\endproclaim

We call $f'$ as in Theorem A a Garuti lifting of the Galois cover $f_k$. 

In the first part of this paper, in $\S2$, we revisit Garuti's theory. We prove the following refined version of Theorem A 
(cf. Theorem 2.5.3).

\proclaim {Theorem B} We use the same notations as in Theorem A. Let $H$ be a quotient of $G$, 
and $g_k:Z_k\to X_k$ the corresponding Galois sub-cover of $f_k$ with Galois group $H$.
Let $h':Z'\to X'\defeq X\times _RR'$ be a  Garuti lifting of the Galois cover $h_k$, defined over the finite extension $R'/R$. 
Then there exists a finite extension $R''/R'$, and a Garuti lifting 
$f'':Y''\to X''\defeq X\times _RR''$ 
of the Galois cover $f_k$ over $R''$, which dominates $h'$, i.e. we have a factorisation 
$f'':Y''@>{g''}>> Z''\defeq Z'\times _{R'}R'' @>{h''\defeq h'\times _{R'}R''}>> X''$,
where $g'':Y''\to Z''$ is a finite morphism between normal $R''$-curves.
\endproclaim

In the course of proving this result we prove a structure theorem concerning 
a certain quotient of the "geometric Galois group" of a $p$-adic open disc, 
which is the most relevant to the lifting problem. 
This result might be of interest independently from the lifting problem. 

Let $\Tilde {X}\defeq \Spf R[[T]]$, $\Tilde {X}_K\defeq \Spec (R[[T]]\times _RK)$ [a $p$-adic
open disc over $K$], and $\Cal X\defeq \Spf R[[T]]\{T^{-1}\}$ the formal boundary of $\Tilde X$ (cf. $\S2$). Let $\Delta$ (resp. $\Delta'$) be the maximal pro-$p$ 
group which classifies geometric Galois covers of $\Tilde X$ (resp. of $\Cal X$) which are pro-$p$, and which are generically \'etale at the level of special fibres 
(cf. 2.3, and 2.4, for more precise definitions). 

\proclaim {Theorem C} (cf. Theorem 2.3.1, and Theorem 2.4.1) The profinite group $\Delta$ is a free pro-$p$ group. Moreover, there exists a natural morphism
$\Delta'\to \Delta$ which makes $\Delta'$ into a direct factor $\Delta$ (cf. 1.1, for the definition, and characterisation,  of a direct factor of a free pro-$p$ group). 
 \endproclaim

In light of Theorem B, we revisit in $\S3$, 3.1, the Oort conjecture. We formulate the following refined version of this conjecture (cf. 3.1, for more details). 

\definition {Oort Conjecture Revisited [Conj-O-Rev]} We use the same notations as in Problem II.
Assume that $G\isom \Bbb Z/m\Bbb Z$ is a cyclic group. Let $H$ be a quotient of $G$, and $g_k:Z_k\to X_k$ the Galois sub-cover of $f_k$ with 
Galois group $H$. 
Then there exists a smooth Galois lifting 
$g:Z'\to X'\defeq X\times _RR'$ 
of $g_k$  over some finite extension $R'/R$.

Furthermore, for every smooth lifting $g$ of the Galois sub-cover $g_k$ of $f_k$ as above,
there exists a smooth lifting
$f:Y''\to X''\defeq X\times _RR''$ 
of $f_k$, over some finite extension $R''/R'$, 
such that $f$ dominates $g$, i.e. we have a factorisation 
$f:Y''\to Z''\defeq Z'\times _{R'}R''@>{g \times _{R'}R''}>> X''$. 
Moreover, $R''$ can be chosen to be the minimal extension 
of $R'$ which contains a primitive $m$-th root of $1$. 
\enddefinition

As for the original Oort conjecture, to prove this revisited version one may reduce to the case where $G\isom \Bbb Z/p^n\Bbb Z$.
In the case $n=1$, both [Conj-O] and [Conj-O-Rev] are clearly equivalent. In 3.2, and in the case where $n=2$, we verify
[Conj-O-Rev] in some cases (cf. Lemma 3.2.1, and Lemma 3.2.2).

The second main part of this paper is motivated by the idea of the search for a path, or a bridge, between Garuti's theory
and the (revisited) Oort conjecture, which may lead to the solution of this conjecture.  We introduce in $\S3$ the notion of fake liftings of cyclic Galois covers
between curves, with the purpose of establishing such a bridge.

Next, we explain the definition of fake liftings, and the simple idea which leads to their existence.

Assume that $G\isom \Bbb Z/p^n\Bbb Z$, $n\ge 1$. Let $H$ be the unique quotient of $G$ with cardinality $p^{n-1}$.  
We use the notations in Problem II, and assume that 
$X=\Bbb P^1_R$. In fact one can reduce the solution of Problem II to this case,
(cf. 3.1, and Lemma 3.1.1)

Let $f_k:Y_k\to \Bbb P^1_k$ be a finite ramified Galois cover with Galois group $G$, and $g_k:X_k\to\Bbb P^1_K$ the Galois sub-cover of $f_k$ 
with Galois group $H$. 
In order to solve [Conj-O-Rev] for the Galois cover $f_k$, and the sub-cover $g_k$, one may proceed by induction on the cardinality of the group $G$.
The case where $G$ has cardinality $p$ is solved in [Se-Oo-Su]. 

So we may 
assume, by an induction hypothesis,  that $g_k$ admits a smooth lifting $g:\Cal X\to \Bbb P^1_R$ defined over $R$, where $\Cal X$ is a smooth $R$-curve, 
i.e. we assume that [Conj-O] holds for the Galois sub-cover $g_k$ of $f_k$. We would like to show that [Conj-O-Rev] is true for $f_k$, and 
the smooth lifting $g$ of the sub-cover $g_k$, i.e. show that $g$ can be dominated by a smooth lifting of $f_k$, after eventually a finite extension of $R$.

Consider all possible Garuti liftings $f:\Cal Y \to \Cal X@>g>> \Bbb P^1_R$ of $f_k$, which dominate the smooth lifting
$g$ of $g_k$.
These Garuti liftings exist by the above refined version of Garuti's theory, in Theorem B, and are a priory defined over a finite extension of $R$. 
For a Garuti lifting $f$ as above, which we can suppose defined over $R$ without loss of generality, the degree of the different in the morphism
$f_K:\Cal Y_K\defeq \Cal Y\times _RK \to \Bbb P^1_K$ between generic fibres is greater than the degree of the different in the morphism $f_k:Y_k\to \Bbb P^1_k$.
Moreover, $\Cal Y$ is smooth over $R$, which implies that [Conj-O-Rev] holds in this case,
if and only if these degrees of different are equal.

Next, we argue by contradiction.
Assume that [Conj-O-Rev] doesn't hold for the Galois cover $f_k$, and the smooth lifting $g$ of the sub-cover $g_k$.
In particular, for all possible Garuti liftings $f$ as above, $\Cal Y$ is not smooth over $R$.
A Garuti lifting $f:\Cal Y\to \Bbb P^1_R$ as above such that the degree of the different in the morphism
$f_K:\Cal Y_K\defeq \Cal Y\times _RK \to \Bbb P^1_K$ between generic fibres
is minimal among all possible $f$'s, is called a fake lifting of the Galois cover $f_k$, relative to the smooth lifting $g$ of $g_k$ (cf. Definition 3.3.2). 

Fake liftings won't exist if [Conj-O-Rev] is true. In fact in order to prove [Conj-O-Rev] for the Galois cover $f_k$, 
relative to the smooth lifting $g$ of $g_k$, it suffices to show that fake liftings $f$ as above do not exist (cf. Remark 3.3.3).

One expects fake liftings to have very special properties, 
which eventually may lead to their non existence. 
Special properties of fake liftings should be encoded in their semi-stable models.

Let $f:\Cal Y\to \Cal X@>g>> \Bbb P^1_R$ be a fake lifting as above, assuming it exists.
In $\S3$ we study the geometry of a minimal semi-stable model ${\Cal Y}'\to \Cal Y$ of $\Cal Y$, which we suppose defined over $R$,
and in which the ramified points in the morphism $f_K:\Cal Y_K\to \Bbb P^1_K$ specialise in smooth distinct points of ${\Cal Y}'_k\defeq \Cal Y'\times _Rk$. 
It turns out that these semi-stable models have indeed very specific properties, which are in some sense reminiscent to the properties of the minimal semi-stable models of smooth liftings of cyclic Galois covers between curves. 

We prove, among other facts, that the configuration of the special fibre $\Cal Y'_k\defeq \Cal Y'\times _Rk$ of the semi-stable model $\Cal Y'$, of the fake lifting $f$, is a tree-like
(cf. Theorem 3.5.4). Moreover, all the irreducible components of positive genus in $\Cal Y'_k$, and which contribute to the difference between the generic and special different
in the morphism $f:\Cal Y\to \Bbb P^1_R$,  are end vertices of the tree associated to ${\Cal Y}'_k$ with special properties (cf. loc. cit).  
In the course of proving this result we establish some of the properties of the minimal semi-stable model of an order $p^n$
automorphism of a $p$-adic open disc, with no inertia at the level of special fibres, that were established in the case $n=1$ in [Gr-Ma1] (cf. Proposition 3.5.3).

Finally, in $\S4$, we introduce the smoothening process for a fake lifting $f:\Cal Y\to \Cal X@>g>> \Bbb P^1_R$ as above. 
The ultimate aim of this process is to show that fake liftings do not exist. This in turn would prove [Conj-O-Rev].

The basic idea of smoothening of the fake lifting $f$ is to construct, starting from $f$, 
a new Garuti lifting $f_1:\Cal Y_1\to \Cal X@>g>> \Bbb P^1_R$, which dominates the smooth lifting $g$ of $g_k$, and such that 
the degree of the different in the morphism $f_{1,K}:\Cal Y_{1,K}\defeq \Cal Y_{1}\times _RK \to \Bbb P^1_K$ between generic fibres is smaller than 
the degree of the different in the morphism $f_K:\Cal Y_K\defeq \Cal Y'\times _RK \to \Bbb P^1_K$ .  We call such $f_1$ a smoothening of $f$.
If this construction is possible, it would imply that
the fake lifting $f$ doesn't exist. Indeed, this would contradict the minimality of the generic different in $f$. Hence this will prove [Conj-O-Rev] for
the Galois cover $f_k$, and the smooth lifting $g$ of the sub-cover $g_k$. 

We describe a formal way, using formal patching techniques, 
to construct a smoothening $f_1$ of the fake lifting $f:\Cal Y\to \Bbb P^1_R$, starting from the  minimal semi-stable model 
$\Cal Y'\to \Cal Y$ of $\Cal Y$ (cf. 4.1). This construction is related to the existence of (internal) irreducible components in 
the special fibre $\Cal P_k$ of the quotient
$\Cal P\defeq \Cal Y'/G$ of the semi-stable model $\Cal Y'$ by $G$, which satisfy certain technical 
conditions arising from the geometry of the semi-stable model $\Cal Y'$,
and the Galois cover $\Cal Y'\to \Cal P$.  
We call such a component a removable vertex of the tree associated
to $\Cal P_k\defeq \Cal P\times _Rk$ (cf. Definition 4.1.2). The existence of a removable vertex in $\Cal P_k$ leads immediately to the existence
of a smoothening $f_1$ of the fake lifting $f$ as above (cf. Definition 4.1.3). 

We show that the smoothening process is possible in the case where $G\isom \Bbb Z/p\Bbb Z$, i.e. 
$n=1$ (cf Proposition 4.2.2).
This gives an alternative proof of the Oort conjecture in this case.
This proof, though simple, is striking in the view of the author in many respects.

First, this proof is not explicit, in the sense that it doesn't produce an explicit lifting of the Galois cover $f_k$. It doesn't even produce any 
automorphism of order $p$ of a $p$-adic open disc.  
Second, the proof doesn't rely (in any form) on the degeneration of the Kummer equation
to the Artin-Schreier equation as in [Se-Oo-Su] (cf. also [Gr-Ma]), 
but rather on the degeneration of the Kummer equation to a radicial equation (cf. proof of Proposition 4.2.2). 
This suggests the possibility of proving [Conj-O] without using the Oort-Sekiguchi-Suwa theory.

In the case where $n=2$, i.e. $G\isom \Bbb Z/p^2\Bbb Z$, we give, in 4.3,  some sufficient conditions for the existence of removable vertices 
which lead to the execution of the smoothening process (cf. Proposition
4.3.1).

Next, we briefly review the content of each section of this paper. In $\S1$ we collect some background material 
which is used in this paper. 
This includes background material on pro-$p$ groups, and formal patching techniques.
In $\S2$ we revisit the theory of Garuti. We prove Theorem B (Theorem 2.5.3), and our main 
result Theorem C concerning the structure of a certain quotient of the geometric Galois group of a $p$-adic open disc (Theorem 2.3.1, and Theorem 2.4.1). 
Both $\S1$ and $\S2$ can be read independently of the rest of the paper.
In $\S3$ we revisit Oort conjecture, and introduce the notion of fake liftings of cyclic Galois covers between curves. 
We then establish the main properties of their minimal semi-stable models in Theorem 3.5.4. 
In $\S4$ we introduce the notion of the smoothening process for fake liftings, and
we investigate on some examples, in degree $p$ and $p^2$, this process. 
This section is dependent on $\S3$.

\subhead
\S 1 Background
\endsubhead
In this section, and for the convenience of the reader,  we collect some background material which is used in this paper. 
We recall some well-known facts on pro-$p$ groups, that will be used in $\S2$ (cf. 1.1). We 
state the main result of formal patching that we use in this paper, and which plays a crucial role in the proof of the main results in 
$\S2$, $\S3$, and $\S4$ (cf. 1.2).
Finally, we recall the degeneration of $\mu_p$-torsors above boundaries of formal fibres at closed points of formal curves (cf. 1.3).

\subhead
1.1  Complements on pro-$p$ Groups
\endsubhead
In this sub-section we fix a prime integer $p>1$. We recall some well-known facts 
on profinite pro-$p$ groups that will be used  in $\S2$.

First, we recall the following characterisations of free pro-$p$ groups.

\proclaim {Proposition 1.1.1} Let $G$ be a profinite pro-$p$ group. Consider
the following properties:

\noindent
(i)\ $G$ is a free pro-$p$ group.

\noindent
(ii)\ The $p$-cohomological dimension of $G$ satisfies $\cd_p(G)\le 1$.

\noindent
(iii)\ Given a surjective homomorphism $\sigma: Q\twoheadrightarrow P$ between 
finite $p$-groups, and a continuous surjective homomorphism $\phi:
G\twoheadrightarrow  P$, there exists a continuous homomorphism 
$\psi:G\to Q$ such that the following diagram is commutative:
$$
\CD
G    @>\id>>  G \\
@V{\psi}VV     @V{\phi}VV \\
Q   @>{\sigma}>>   P  \\
\endCD
$$

Then the following equivalences hold:
$$(i)\Longleftrightarrow (ii) \Longleftrightarrow (iii).$$
\endproclaim

\demo{Proof} Well-known (cf. [Se], and [Ri-Za], Theorem 7.7.4).
\qed
\enddemo

Next, we recall the notion of a direct factor of a free pro-$p$ group 
(cf. [Ga], 1, the discussion preceding  Proposition 1.8).

\definition {Definition 1.1.2 (Direct Factors of Free  pro-$p$ Groups)} 
Let $F$ be a free pro-$p$ group, and $H$ a closed 
subgroup of $F$. Let $\iota: H\to F$ be the natural homomorphism. We say that 
$H$ is a direct factor of $F$ if there exists a continuous homomorphism 
$s:F\to H$ such that $s\circ \iota=\id_H$ [$s$ is necessarily surjective]. 
We then have a natural exact sequence 
$$1\to N\to F @>{s}>> H\to 1,$$ 
where $N\defeq \Ker s$, and $F$ is isomorphic to the free direct 
product 
$$H*N$$
In particular, $N$ is also a direct factor of $F$ (cf. loc. cit).
\enddefinition

One has the following cohomological characterization of direct factors of free 
pro-$p$ groups.

\proclaim {Proposition 1.1.3}
Let $H$ be a pro-$p$ group, and $F$ a free pro-$p$ group. Let 
$\sigma :H\to F$ 
be a continuous homomorphism. Assume that the map induced by $\sigma$ on 
cohomology: 
$$h^1(\sigma):H^1(F,\Bbb Z/p\Bbb Z)\to H^1(H,\Bbb Z/p\Bbb Z)$$
is surjective. [Here $\Bbb Z/p\Bbb Z$ is considered as a trivial discrete module]. 
Then $H$ is a direct factor of $F$.
\endproclaim

\demo{Proof} cf. [Ga], Proposition 1.8.
\qed
\enddemo

\subhead
1.2 Formal Patching
\endsubhead
In this sub-section we explain the procedure which allows to construct [Galois] covers
of curves in the setting of formal geometry, by patching together covers of formal affine
curves with covers of formal fibres at closed points of the special fibre (cf. [Sa], 1, for more details).
We also recall the [well-known] local-global principle for liftings of Galois covers of curves.

Let $R$ be a complete discrete valuation ring, with fraction field $K$, residue field $k$, and 
uniformiser $\pi$. Let $X$ be an admissible formal $R$-scheme which is an $R$-curve, by which
we mean that the special fibre $X_k\defeq X\times _Rk$ is a reduced one-dimensional scheme
of finite type over $k$. Let $Z$ be a finite set of closed points of $X$. For a point $x\in Z$,
let $X_x\defeq \Spf \hat {\Cal O}_{X,x}$ be the formal completion of $X$ at $x$, which is the formal fibre at the point $x$.
Let $X'$ be a formal open sub-scheme of $X$ whose special fibre is $X_k\setminus Z$. 

For each closed point $x\in Z$,
let $\{\Cal P_i\}_{i=1}^n$ be the set of minimal prime ideals of $\hat {\Cal O}_{X,x}$ which contain $\pi$; they correspond 
to the branches $\{\eta_i\}_{i=1}^n$ of the completion of $X_k$ at $x$, and let $X_{x,i}\defeq \Spf \hat {\Cal O}_{x,\Cal P_i}$
be the formal completion of the localisation of $X_x$ at $\Cal P_i$. The local ring $\hat {\Cal O}_{x,\Cal P_i}$ is a complete 
discrete valuation ring. The set $\{X_{x,i}\}_{i=1}^n$ is the set of boundaries of the formal fibre $X_x$. For each $i\in \{1,...,n\}$,
we have a canonical morphism $X_{x,i}\to X_x$.

\definition
{Definition 1.2.1}
With the same notations as above, a ($G$-)cover patching data for the pair $(X,Z)$
consists of the following.

\noindent
(i)\  A finite (Galois) cover $Y'\to X'$ (with Galois group $G$).

\noindent
(ii)\ For each point $x\in Z$, a finite (Galois) cover $Y_x\to X_x$ (with Galois group $G$).

The above data (i) and (ii) must satisfy the following compatibility condition.

\noindent
(iii)\ If $\{X_{x,i}\}_{i=1}^n$ are the boundaries of the formal fibre at the point $x$, then for each $i\in \{1,...,n\}$
is given a ($G$-equivariant) $X_x$-isomorphism 
$$\sigma _i:Y_x\times _{X_x}X_{x,i}\isom Y'\times _{X'} X_{x,i}.$$
Property (iii) should hold for each $x\in Z$.
\enddefinition

The following is the main patching result that we will use in this paper.

\proclaim {Proposition 1.2.2} With the same notations as above. Given a ($G$-)cover patching data as in Definition 1.2.1, 
there exists a unique, up to isomorphism, (Galois) cover $Y\to X$ (with Galois group $G$) which induces the above ($G$-)cover
in Definition 1.2.1, (i), when restricted to $X'$, and induces the above ($G$-)cover in Definition 1.2.1, (ii), 
when pulled-back to $X_x$,  for each point $x\in Z$.
\endproclaim

\definition
{1.2.3}  With the same notations as above, let $x\in Z$, and $\Tilde X_k$ the normalisation of $X_k$. There is a one-to-one correspondence
between the set of points of $\Tilde X_k$ above $x$, and the set of boundaries of the formal fibre at the point $x$. Let $x_i$ be the point
of $\Tilde X_k$ above $x$ which corresponds to the boundary $X_{x,i}$, for $i\in \{1,...,n\}$. Assume that the point $x\in X_k(k)$ is rational.
Then the completion of $\Tilde X_k$ at $x_i$ is isomorphic to the spectrum of a ring of formal power series $k[[t_i]]$ in one variable over $k$,
where $t_i$ is a local parameter at $x_i$.  

The complete local ring $\hat {\Cal O}_{x,\Cal P_i}$ is a discrete valuation ring with residue field isomorphic
to $k((t_i))$. Let $T_i\in \hat {\Cal O}_{x,\Cal P_i}$ be an element which lifts $t_i$. Such an element is called a parameter of $\hat {\Cal O}_{x,\Cal P_i}$.
Then there exists an isomorphism $\hat {\Cal O}_{x,\Cal P_i}\isom R[[T_i]]\{T_i^{-1}\}$, where  
$$R[[T]]\{T^{-1}\} \defeq \{\sum _{i=-\infty}^{\infty}a_iT^i,\ \underset 
{i\to -\infty} \to  \lim \vert a_i \vert=0\},$$ 
and $\vert \ \vert$ is a normalised absolute value of $R$.
\enddefinition

As a direct consequence of the above patching result, and the theorems of liftings of \'etale covers (cf. [Gr]),
one obtains the following  [well-known] local-global principle for liftings of (Galois) covers of curves.

\proclaim
{Proposition 1.2.4}
Let $X$ be a proper, flat, algebraic (or formal) $R$-curve, and let $Z\defeq \{x_i\}_{i=1}^n$ be a finite set of closed points of $X$.
Let $f_k:Y_k\to X_k$ be a finite generically separable (Galois) cover (with Galois group $G$), whose branch locus is contained in
$Z$. Assume that for each $i\in \{1,...,n\}$, there exists a (Galois) cover $f_i:Y_i\to \Spf \hat {\Cal O} _{X,x_i}$ (with Galois group $G$)
which lifts the cover $\hat Y_{k,x_i}\to \Spec \hat {\Cal O} _{X_k,x_i}$ induced by $f_k$, where $\hat {\Cal O}_{X_k,x_i}$ (resp. $\hat Y_{k,x_i}$)
 denotes the completion of $X_k$ at $x_i$ (resp. the completion of $Y_k$ above $x_i$). Then 
 there exists a unique, up to isomorphism, (Galois) cover $f:Y\to X$ (with Galois group $G$)
which lifts the cover $f_k$, and which is isomorphic to the cover $f_i$ when pulled back to $\Spf \hat {\Cal O}_{X,x_i}$, for each $i\in \{1,...,n\}$.
\endproclaim

\subhead
1.3 Degeneration of $\mu_p$-torsors
\endsubhead
In this sub-section, we recall the [well-known] degeneration of $\mu_p$-torsors from zero to positive characteristics, 
above the boundaries of formal fibres of formal $R$-curves at closed points. 

Here $R$ denotes a complete discrete valuation ring 
of unequal characteristics, with fraction field $K$, residue field $k$ of characteristic $p>0$, uniformiser $\pi$, and which contains $\zeta$: a primitive 
$p$-th root of $1$. We write $\lambda \defeq \zeta-1$. We denote by $v_K$ the valuation of $K$, which is normalised by $v_K(\pi)=1$.

First, we recall the definition of a certain class of $R$-group schemes (cf. [Se-Oo-Su], for more details).

\subhead
{1.3.1 Torsors under finite and flat $R$-group schemes of rank $p$: the group schemes $\Cal G_n$ and $\Cal H_n$} 
\endsubhead
Let $n\ge 1$ be an integer. Define the affine $R$-group scheme 
$$\Cal G_{n,R}\defeq \Spec (A_n)$$ 
as follows.

\par
\noindent
(i) $A_n\defeq R[X,\frac {1}{1+\pi ^nX}].$

\noindent
(ii)\ The comultiplication $c_n:A_n\to A_n\otimes_R A_n$ is defined by
$c_n(X)\defeq X\otimes 1+1\otimes X+ \pi^n X\otimes X.$

\noindent
(iii)\ The coinverse $i_n:A_n\to A_n$ is defined by 
$i_n(X)\defeq -\frac {X}{1+\pi ^n X}.$

\noindent
(iv)\ The counit $\epsilon_n :A_n\to R$ is defined by 
$\epsilon_n (X)\defeq 0.$

One verifies [easily] that $\Cal G_n\defeq \Cal G_{n,R}$ is an affine, commutative, and smooth $R$-group scheme, 
with generic fibre $(\Cal G_n)_K\isom \Bbb G_{m,K}$, and special fibre 
$(\Cal G_n)_k\isom \Bbb G_{a,k}$. 

Next, we introduce some finite and flat group schemes of rank $p$.
Assume that $n$ satisfies the following condition
$$0<n(p-1)\le v_K(p).\tag *$$
Assuming that the above condition $(*)$ holds, consider the map
$$\phi_{n}:\Cal G_n\to \Cal G_{pn},$$
given by
$$X\mapsto \frac {(1+\pi^n X)^{p}-1}{\pi ^{pn}}.$$
Then $\phi_{n}$ is a surjective homomorphism of $R$-group schemes.
Denote by
$$\Cal H_{n}\defeq \Cal H_{n,R}\defeq \Ker (\phi_{n}).$$
It is a finite and flat commutative group scheme of rank $p$. Under the assumption
(*) one verifies [easily] that the generic fibre 
$\Cal H_{n,K}\defeq \Cal H_{n}\otimes _R K \isom \mu _{p,K}$
is \'etale, and the special fibre
$\Cal H_{n,k}\defeq \Cal H_{n}\otimes _R k \isom \alpha _{p,k}$
is radicial of type $\alpha_{p}$, if $n(p-1)<v_K(\lambda)$, and
$\Cal H_{n,k}\defeq \Cal H_{n}\otimes _R k \isom (\Bbb Z/p\Bbb Z)_k$
is \'etale,  if $n(p-1)=v_K(\lambda)$.

Let $\Cal U\defeq \Spf A$ be a formal affine $R$-scheme, and 
$f:\Cal V\to \Cal U$
a torsor under the finite group scheme $\Cal H_{n}$, for some $n$ as above satisfying (*).
Then there exists a regular function $u\in A$, such that the image $\bar u$ of $u$ in $\overline A\defeq {A}/{\pi A}$
is not a $p$ power if $n(p-1)<v_K(\lambda)$,
$1+\pi^{pn}u$ is defined up to
multiplication by a $p$-th power of the form $(1+\pi^n v)^{p}$, and the torsor
$f$ is given by an equation 
${(X')}^{p}=(1+\pi^n X)^{p}=1+\pi ^{np}u,$
where $X'$ and $X$ are indeterminates. Moreover, the natural morphism 
$f_k:\Cal V_k\to \Cal U_k$
between the special fibres is either the $\alpha_{p}$-torsor given by the equation
${x}^{p}=\bar u,$
where $x= X\mod \pi$, and $\bar u= u \mod \pi$, if $n(p-1)<v_K(\lambda)$. Or,
is the $\Bbb Z/p\Bbb Z$-torsor given by the equation
$x^p-x=\bar u$
where $x=X\mod \pi$, and $\bar u= u \mod \pi$, if $n(p-1)=v_K(\lambda)$

Next, we recall the degeneration of $\mu_p$-torsors on the boundary 
$\Cal X\defeq \Spf R[[T]]\{T^{-1}\}$ of formal fibres of germs of formal $R$-curves. Here  
$R[[T]]\{T^{-1}\}$ is as in 1.2.3.
Note that $R[[T]]\{T^{-1}\}$ is a complete discrete valuation ring, with uniformising parameter $\pi$, and residue field 
$k((t))$, where $t= T\mod \pi$. The following result will be used in $\S3$ and $\S4$.

 \proclaim
 {Proposition 1.3.2} Let $A\defeq R[[T]]\{T^{-1}\}$ (cf. above definition), and $f:\Spf B\to \Spf A$
a non-trivial Galois cover of degree $p$. Assume that the ramification index of the corresponding 
 extension of discrete valuation rings equals $1$.
 Then $f$ is a torsor under a finite and flat $R$-group scheme $G$ of rank $p$. Let $\delta$ be the degree of the different in the above extension.
 The following cases occur.
 
 (a)\ $\delta=v_K(p)$. Then $f$ is a torsor under the group scheme 
 $G=\mu_{p,R}$, and two cases occur.
 
 (a1)\ For a suitable choice of the parameter $T$ of $A$ the torsor $f$ is given, after eventually a finite extension of $R$, by
 an equation $Z^p=T^h$. In this case we say that the torsor $f$ has a degeneration of type $(\mu_p,0,h)$.

 (a2)   \ For a suitable choice of the parameter $T$ of $A$ the torsor $f$ is given, after eventually a finite extension of $R$, by
 an equation $Z^p=1+T^m$ where $m$ is a positive integer prime to $p$. In this case we say that the torsor $f$ has a degeneration of type $(\mu_p,-m,0)$.

 (b)\ $0<\delta <v_K(p)$. Then $f$ is a torsor under the group scheme $\Cal H_{n,R}$, where $n$ is such that
 $\delta=v_K(p)-n(p-1)$. Moreover, for a suitable choice of the parameter $T$ the torsor $f$ is given, after eventually a finite extension of $R$,
 by an equation $Z^p=1+\pi ^{pn}T^m$, with $m\in \Bbb Z$ prime to $p$. In this case we say that the torsor $f$ has a degeneration of type $(\alpha_p,-m,0)$.

 (c)\ $\delta=0$. Then $f$ is an \'etale torsor under the $R$-group scheme $G=\Cal H_{v_K(\lambda),R}$, and is given, 
 after eventually a finite extension of $R$, by an equation $Z^p=\lambda^pT^m+1$, where $m$ is a negative integer prime to $p$, for 
 a suitable choice of the parameter $T$ of $A$.
 In this case we say that the torsor $f$ has a degeneration of type $(\Bbb Z/p\Bbb Z,-m,0)$.
 \endproclaim

\demo {Proof}
See [Sa], Proposition 2.3.
\qed
\enddemo

\subhead
\S 2. Pro-p Quotients of the Geometric Galois Group of a $p$-adic Open Disc
\endsubhead

\definition {2.1 Notations}
The following notations will be used in this section and the subsequent ones, unless we specify otherwise.
\enddefinition

\noindent
$p>1$ is a fixed prime integer.

\noindent
$R$ will denote a complete discrete valuation ring 
of unequal characteristics, with uniformising parameter $\pi$.

\noindent
$K\defeq \Fr (R)$ is the quotient field of $R$, $\char (K)=0$.

\noindent
$k\defeq {R}/{\pi R}$ is the residue field of $R$, which we assume to be algebraically 
closed of characteristic $p>0$.

\noindent
$v_K$ will denote the valuation of $K$, which is normalised by $v_K(\pi)=1$.

\noindent
For an $R$-(formal) scheme $X$ we will denote by
$X_K\defeq X\times _RK$
(resp. $X_k\defeq X\times _Rk)$
the generic (resp. special) fibre of $X$.

\definition {2.2}
Next, we would like to state a result of Garuti, in [Ga], which concerns the structure
of the pro-$p$ geometric fundamental group of a $p$-adic annulus of thickness zero (cf. Proposition 2.2.3). 
First, we recall how one defines the fundamental group of a rigid analytic affinoid space.
\enddefinition

Let  $\Cal X=\Spf \Cal A$ 
be an affine $R$-formal scheme which is topologically of finite type. Thus, 
$\Cal A$ is a $\pi$-adically complete noetherian $R$-algebra. Let 
$A\defeq \Cal A\otimes _RK$ be the corresponding Tate algebra, and 
$X\defeq \Sp A$ 
the associated  rigid analytic affinoid space, which is the generic fiber of 
$\Cal X$ in the sense of Raynaud (cf. [Ab]). 

Assume that $X$ is integral and 
geometrically connected. Let $\eta$ be a geometric point of the affine 
scheme $\Spec A$ above the generic point of $\Spec A$. Then $\eta$ 
determines naturally an algebraic closure $\overline K$ of $K$, and a geometric 
point of $\Spec (A\times _K\overline K)$, which we will also denote by $\eta$.

\definition {Definition 2.2.1 (\'Etale Fundamental Groups of Affinoid Spaces)} 
(See also [Ga], D\'efinition 2.2, and D\'efinition 2.3).
We define the \'etale fundamental group of $X$ with base point $\eta$ by
$$\pi_1(X,\eta)\defeq \pi_1(\Spec A,\eta),$$
where $\pi_1(\Spec A,\eta)$ is the \'etale fundamental group of the connected 
scheme $\Spec A$ with base point $\eta$ in the sense of Grothendieck (cf. [Gr]).
Thus, $\pi_1(X,\eta)$ naturally classifies rigid analytic coverings 
$Y\to X$, where $Y=\Sp B$, and $B$ is a finite $A$-algebra which is \'etale over $A$.

There exists a natural continuous surjective homomorphism
$$\pi_1(X,\eta)\twoheadrightarrow \Gal (\overline K,K).$$
We define the geometric fundamental group $\pi_1(X,\eta)^{\geo}$ of $X$ so 
that the following sequence is exact:
$$1\to \pi_1(X,\eta)^{\geo}\to \pi_1(X,\eta)\to \Gal (\overline K,K)\to 1.$$
\enddefinition

\definition {Remark 2.2.2} If $L/K$ is a finite field extension contained in 
$\overline K/K$, and $X_L\defeq X\times _KL$ 
is the affinoid rigid analytic 
space obtained from $X$ by extending scalars, then we have a natural 
commutative diagram:
$$
\CD 
1 @>>>  \pi _1(X_L,\eta)^{\geo} @>>> \pi_1(X_L,\eta) @>>> \Gal (\overline K/L) 
@>>> 1 \\
@.  @VVV   @VVV     @VVV \\
1 @>>> \pi _1 (X,\eta)^{\geo} @>>> \pi_1(X,\eta) @>>> \Gal (\overline K/K) 
@>>> 1\\
\endCD
$$ 
where the two right vertical maps are injective homomorphisms, and the left 
vertical map is an isomorphism. 

The geometric fundamental group $\pi_1(X,\eta)
^{\geo}$ is strictly speaking not the fundamental group of a rigid analytic 
space [since $\overline K$ is not complete]. It is, however, the projective limit 
of fundamental groups of rigid affinoid spaces. More precisely, there exists 
a natural isomorphism
$$\pi_1(X,\eta)^{\geo} \isom \underset{L/K} \to {\varprojlim}\pi_1(X\times _KL
,\eta),$$
where the limit is taken over all finite extensions $L/K$ contained in 
$\overline K$.
\enddefinition

Next, we introduce some notations involved in the statement of Garuti's result.

For a finite field extension $L/K$ contained in $\overline K/K$ we will denote 
by 
$$D_0\defeq D_{0,L}\defeq \Sp L<X>$$ 
the unit closed disc [centred at $X=0$], and 
$$C_L\defeq \Sp \frac {L<X,Y>}{(XY-1)}$$ 
the annulus of thickness $0$ which is the ``boundary'' of $D_0$. Here $L<X>$ (resp. 
$L<X,Y>$) denotes the Tate algebra in the variable $X$ (resp. the variables 
$X$ and $Y$).

We denote by $\Bbb P^1_L$ the rigid analytic projective line 
over $L$ which is obtained by patching the closed discs  $D_0=D_{0,L}\defeq \Sp L<X>$, 
and $D_{\infty}=D_{\infty,L}\defeq \Sp L<Y>$, along the annulus $C_L$ (see above), 
via the identification $X\mapsto \frac {1}{Y}.$

Let $S=\{a_1,a_2,...,a_n\}$ 
be a finite set of closed points of $\Bbb P^1_K$ which contains $\{0,\infty\}$, and such that
$S\cap C_K=\emptyset $. 
We view $S\subset \Bbb P^1_K$ as a closed subscheme of $\Bbb P^1_K$, and write 
$S_L\defeq S\times _KL$.
Let $\eta$ be a geometric point of $\Bbb P^1_{\overline K}\defeq \Bbb P^1_K\times _K {\overline K}$ 
above the generic point of $\Bbb P^1_{\overline K}$. We denote by 
$\pi_1(\Bbb P^1_L\setminus S_L, \eta)$ 
the algebraic \'etale fundamental group 
of  $\Bbb P^1_L\setminus S_L$ with base point $\eta$. Write 
$$C\defeq C_K\defeq \Sp \frac {K<X,Y>}{(XY-1)}.$$

The natural embedding 
$C\times _KL=C_L\to \Bbb P^1_L$ 
induces a natural continuous homomorphism 
$$\pi_1(C,\eta)^{\geo}\to \pi_1(\Bbb P^1_L\setminus S_L, \eta),$$ 
and by passing to the projective limit a continuous homomorphism 
$$\pi_1(C,\eta)^{\geo}\to \pi_1(\Bbb P^1_{\overline K}\setminus S_{\overline K}, \eta) 
\defeq \underset{L/K} \to {\varprojlim} \pi_1(\Bbb P^1_L\setminus S_L, \eta),$$
where $L/K$ runs over all finite extensions contained in $\overline K$.

Let 
$\pi_1(C,\eta)^{\geo,p}$ 
be the maximal pro-$p$ quotient of $\pi_1(C,\eta)^{\geo}$,
and 
$\pi_1(\Bbb P^1_{\overline K}\setminus S_{\overline K}, \eta)^p$ 
the maximal pro-$p$ 
quotient of $\pi_1 (\Bbb P^1_{\overline K}\setminus S_{\overline K}, \eta)$. The above 
homomorphism $\pi_1(C,\eta)^{\geo}\to \pi_1(\Bbb P^1_{\overline K}\setminus S_{\overline K}, 
\eta)$ induces a natural continuous homomorphism
$$\phi_S: \pi_1(C,\eta)^{\geo,p}\to \pi_1(\Bbb P^1_{\overline K}\setminus 
S_{\overline K}, \eta)^p,$$ 
which induces, by passing to the projective limit, a continuous homomorphism 
$$\phi\defeq \underset{S} \to {\varprojlim}\ \phi_S: \pi_1(C,\eta)^{\geo,p}\to 
\underset{S}\to {\varprojlim}\ \pi_1(\Bbb P^1_{\overline K}\setminus S_{\overline K}
,\eta)^p,$$ 
where the limit is taken over all finite set of closed points of 
$\Bbb P^1_{K}\setminus C$ which contain $\{0,\infty\}$.

The profinite pro-$p$ group 
$\underset{S} \to {\varprojlim}\  \pi_1(\Bbb P^1_{\overline K}\setminus S_{\overline K}, \eta)^p$ 
is a free pro-$p$ group, 
as follows from the well-known structure of algebraic fundamental groups in 
characteristic $0$ (cf. [Gr]).

The following result is one of the main technical results in [Ga], which we will use in this section 
(cf. proof of Theorem 2.3.1, and proof of Theorem 2.4.1).

\proclaim {Proposition 2.2.3 (Garuti)} The natural continuous homomorphism
$$\phi\defeq \underset{S} \to {\varprojlim}\ \phi_S: \pi_1(C,\eta)^{\geo,p}\to 
\underset{S} \to {\varprojlim}\ \pi_1(\Bbb P^1_{\overline K}\setminus S_{\overline K}, 
\eta)^p,$$
where the limit is taken over all finite set of closed points of 
$\Bbb P^1_{K}\setminus C$ which contain $\{0,\infty\}$, makes $\pi_1(C,\eta)^{\geo,p}$ into a direct factor of $\underset{S}\to
{\varprojlim}\ \pi_1(\Bbb P^1_{\overline K}\setminus S, \eta)^p$. 
In particular, $\pi_1(C,\eta)^{\geo,p}$ is a free pro-$p$ group.
\endproclaim

\demo{Proof}
See [Ga], Lemma 2.11.
\qed
\enddemo

\definition {2.3}
Next, we will investigate the structure of a certain quotient of the 
``geometric absolute Galois group'' of a $p$-adic open disc. First, we will define this quotient (see the definition of the profinite group
$\Delta$ below).
\enddefinition

Write 
$$\Tilde X\defeq \Spec R[[T]],$$
and
$$\Tilde X_K\defeq \Tilde X\times _RK=\Spec (R[[T]]\otimes_RK).$$ 
[$\Tilde X_K$ is what we shall refer to as a $p$-adic open disc (over $K$)]. 

Let
$\Tilde S\defeq \{ x_1,x_2,...,x_n\}\subset \Tilde X_K$ be a finite set of closed points of 
$\Tilde X_K$. We view $\Tilde S\subset \Tilde X_K$ as a closed subscheme of $\Tilde X_K$. Write 
$$U_{K,\Tilde S}\defeq \Tilde X_K\setminus \Tilde S,$$ 
and let $\eta$ be a geometric point of $\Tilde X_K$ above the generic point of 
$\Tilde X_K$. We have a natural exact sequence of profinite groups:
$$1\to \pi_1(U_{K,\Tilde S}\times _K\overline K,\eta)\to \pi_1(U_{K,\Tilde S},\eta)\to 
\Gal (\overline K,K)\to 1.$$
By passing to the projective limit over all finite set of closed points 
$\Tilde S\subset \Tilde X_K$, we obtain a natural exact sequence:
$$1\to \underset{\Tilde S} \to {\varprojlim}\ \pi_1(U_{K,\tilde S}\times _K\overline K,\eta)
\to \underset {\tilde S} \to {\varprojlim}\ \pi_1(U_{K,\Tilde S},\eta)\to \Gal (\overline 
K/K)\to 1.$$

Let 
$L\defeq \Fr (R[[T]])$ 
be the quotient field of the formal power series ring $R[[T]]$. 
The generic point $\eta$ determines an algebraic closure $\overline L$ of $L$.
We have a natural exact sequence of Galois groups:
$$1\to \Gal (\overline L/L.\overline K)\to \Gal (\overline L/L)\to 
\Gal (\overline K/K)\to 1.$$
Moreover, there exist natural identifications:
$$\Gal (\overline L/\overline K.L)\isom \underset{\Tilde S} \to {\varprojlim}\ 
\pi_1(U_{K,\tilde S}\times _K\overline K,\eta),$$
and
$$\Gal (\overline L/L)\isom \underset{\Tilde S} \to {\varprojlim}\ \pi_1(U_{K,\tilde S},
\eta),$$
where $\Tilde S$ is as above [$\Gal (\overline L/\overline K.L)$ is what we shall refer to as the geometric
Galois group of a $p$-adic open disc].

Let 
$$I\defeq I_{(\pi)}\subset \Gal (\overline L.\overline K.L)$$ 
be the normal closed subgroup which is generated by the inertia subgroups above the ideal 
$(\pi)$ of $R[[T]]$, which is generated by $\pi$. Write
$$\overline \Delta \defeq \Gal (\overline L/\overline K.L)/I.$$

Note that, by definition, the profinite group $\overline \Delta$ classifies
finite Galois covers 
$\Tilde Y_{K'}\to \Tilde X_{K'}\defeq \Tilde X\times _RK',$ 
where $K'$ is a finite extension of $K$ with valuation ring $R'$, $K'$ is algebraically 
closed in $\Tilde Y_{K'}$, $\pi'$ is a uniformising parameter of $K'$, and the 
natural morphism $\Tilde Y'\to \Tilde X'\defeq \Tilde X\times _RR'$, where $\Tilde Y'$ 
is the normalisation of $\Tilde X'$ in $\Tilde Y_{K'}$, is \'etale above the 
generic point of the special fiber $\Tilde X'_k\defeq \Tilde X'\times _{R'}k$
of $\Tilde X'$. In particular,
the special fiber $\Tilde Y'_k\defeq \Tilde Y'\times _{R'}k$ is reduced and the 
natural morphism $\Tilde Y'_k\to \Tilde X'_k$ is generically \'etale. 

Let 
$$\Delta \defeq \overline \Delta ^p$$ 
be the maximal pro-$p$ quotient of $\overline \Delta$. 
Our main technical result in this section is the following:

\proclaim {Theorem 2.3.1} The profinite group $\Delta$ is a free pro-$p$ group.
\endproclaim 

\demo {Proof} The two main ingredients of the proof are the technical result 
of Garuti in Proposition 2.2.3, and a result of Harbater, Katz, and Gabber (cf. [Ha]
, and [Ka]). 

We will show that the profinite pro-$p$ group  $\Delta$ satisfies 
property (iii) in Proposition 1.1.1. 
Let $Q \twoheadrightarrow P$ be a surjective 
homomorphism between finite $p$-groups. 
Let $\phi: \Delta \twoheadrightarrow P$ be a surjective homomorphism. We will 
show that $\phi$ lifts to a homomorphism $\psi: \Delta \to Q$. The homomorphism
$\phi$ corresponds to a finite Galois extension $\overline L'/\overline K.L$
with Galois group $P$. We can [without loss of generality] assume that this 
extension is defined over $K$, thus descends to a finite Galois extension
$L'/L$ with Galois group $P$ where $K$ is algebraically closed in $L'$.

Let  $A$ be the integral closure of $R[[T]]$ in $L'$. We have a finite 
morphism $f:\Spec A\to \Spec R[[T]]$ which is [by assumptions] Galois with Galois group $P$, 
is \'etale above the point $(\pi)\in \Spec R[[T]]$, and with $\Spec A$
geometrically connected. In particular, $f$ induces at the level of
special fibres a finite generically 
Galois cover $\bar f:\Spec (A/\pi A)\to \Spec k[[t]]$
[where $t=T \mod \pi$] with Galois group $P$.
We will assume [in order to simplify the arguments below] that $\Spec (A/\pi A)$
is connected, the general case is treated in a similar fashion.

By a result of Harbater, Katz and Gabber (cf. [Ha], and [Ka]) there exists a finite Galois cover 
$\bar g:\overline Y\to \Bbb P^1_k$ with Galois group $P$, which is \'etale 
outside a unique closed point $\infty$ of $\Bbb P^1_k$ with local parameter 
$t$, $\bar g$ is totally ramified above $\infty$, $\overline Y$ is connected,
and such that the Galois cover 
above the formal completion $\Spec \hat \Cal O_{\Bbb P^1_k,\infty}$ 
of $\Bbb P^1_k$ at $\infty$, which is naturally 
induced by $\bar g$, is isomorphic to $\bar f$. 
Let $\Bbb A^1_k\defeq\Bbb P^1_k\setminus \{\infty\}$. The restriction 
$\bar f':\overline Y'\to \Bbb A^1_k$ of $\bar g$ to $\Bbb A^1_k$ is an 
\'etale Galois cover with Galois group $P$, and $\overline Y'$ is connected. 

Consider the rigid analytic projective line
$\Bbb P^1_K$ which is obtained by patching the closed unit disc 
$D_{\infty}\defeq \Sp K<T>$ [centered at $T=\infty$] with the closed disc
$D_{0}\defeq \Sp K<S>$ [centered at $S=0$] along the annulus 
$C\defeq \Sp {K<T,S>}/{(ST-1)}$ [of thickness $0$], via the 
identification $S\mapsto \frac {1}{T}$. The \'etale Galois cover
$\bar f'$ lifts [uniquely up to isomorphism] to an \'etale Galois cover 
$f':Y'\to D_{0}$ by the theorems of liftings of \'etale covers (cf. [Gr]), 
whose restriction $\tilde f':\tilde Y'\to C$ to the annulus $C$ is an \'etale Galois 
cover with Galois group $P$.

By using formal patching techniques \`a la Harbater one can construct a [connected] 
rigid analytic Galois cover $g:Y\to \Bbb P^1_K$ with Galois group $P$ whose 
restriction to the annulus $C$ is isomorphic to $\tilde f'$, 
and which above the formal completion at $T=\infty$ induces the above Galois cover $f$.
(see the arguments used in [Ga], and Proposition 1.2.2).

The Galois cover $g$ is ramified above a finite 
set of closed points $\Tilde S\subset \Bbb P^1_K$ [which are contained
in the interior $D_{\infty}^{\op}$ of the closed disc $D_{\infty}$],
hence gives rise naturally to a surjective 
homomorphism $\phi_1:\pi_1(\Bbb P^1_{\overline K}\setminus \Tilde S_{\overline K},\eta)
\twoheadrightarrow P$, and also  a surjective homomorphism
$\underset{S} \to {\varprojlim}\ \pi_1(\Bbb P^1_{\overline K}\setminus S_{\overline K}, 
\eta)^p\twoheadrightarrow P$ [where $S$ and the projective limit are as in Proposition 2.2.3]. 
We also denote by $\phi_1:\pi_1(C,\eta)^{\geo,p}
\twoheadrightarrow P$ the corresponding homomorphism induced on the 
direct factor $\pi_1(C,\eta)^{\geo,p}$ of $\underset{S} \to {\varprojlim}\ \pi_1(\Bbb P^1_{\overline K}\setminus  S_{\overline K}, 
\eta)^p$.

Next, let $\bar f_1:\Spec \overline B\to \Spec k[[t]]$ be a finite connected 
Galois cover with Galois group $Q$ which dominates the above Galois cover $\bar f:\Spec
(A/\pi A)\to \Spec k[[t]]$ with Galois group $P$. Note that such $\bar f_1$ exists, since the maximal
pro-$p$ quotient of the absolute Galois group of $k((t))$ is a 
free pro-$p$ group. Let $\bar g_1:\overline Y_1\to \Bbb P^1_k$
be the finite Galois cover with Galois group $Q$ which is \'etale outside
$\infty$, which induces above $\Spec \hat \Cal O_{\Bbb P^1_k,\infty}$ a finite 
Galois cover which is isomorphic to $\bar f_1$, and let 
$\bar g_1':\overline Y_1'\to \Bbb A^1_k$ be its restriction to $\Bbb A^1_k$ 
[the Galois cover $\bar g_1$ exists by the above result of Harbater, Katz and Gabber
(cf. loc. cit.)]. By construction the Galois cover  $\bar g_1:\overline Y_1\to \Bbb P^1_k$
dominates the Galois cover $\bar g:\overline Y\to \Bbb P^1_k$.
The \'etale Galois cover $\bar g_1'$ lifts to a finite \'etale Galois cover
$f_1':Y_1'\to D_{0}$ with Galois group $Q$, which by construction
dominates the lifting $f':Y'\to D_{0}$ of $\bar f'$. The restriction
of $f_1'$ to the annulus $C$ is a finite \'etale Galois cover $\tilde f_1':
\tilde Y_1'\to C$ with Galois group $Q$ which dominates the Galois cover $\tilde f':\tilde Y'
\to C$.

Let $N$ be a complement of $\pi_1(C,\eta)^{\geo,p}$ in $\underset{S} 
\to {\varprojlim}\ \pi_1(\Bbb P^1_{\overline K}\setminus  S_{\overline K}, 
\eta)^p$ (cf. Proposition 2.2.3). The Galois cover 
$\tilde f'_1:\tilde Y'_1\to C$ (resp. $\tilde f':\tilde Y'\to C$) 
corresponds to the continuous homomorphism $\phi_2:\pi_1(C,\eta)^{\geo}
\to Q$ (resp. $\phi_1:\pi_1(C,\eta)^{\geo,p}\to P$), and $\phi_2$ dominates
$\phi_1$ [by construction]. Also the above homomorphism $\phi _1:
\varprojlim\ \pi_1(\Bbb P^1_{\overline K}\setminus  S_{\overline K}, 
\eta)^p \twoheadrightarrow  P$ 
induces naturally a continuous homomorphism
$\psi_1:N\to P$. 

The pro-$p$ group $N$ being free one can lift the homomorphism
$\psi_1$ to a homomorphism $\psi_2:N\to Q$ which dominates $\psi_1$.
The profinite pro-$p$ group 
$\underset{S} 
\to {\varprojlim}\ \pi_1(\Bbb P^1_{\overline K}\setminus  S_{\overline K}, 
\eta)^p$
being the free direct product of $N$ 
and $\pi_1(C,\eta)^{\geo,p}$ one can construct a continuous homomorphism
$\psi: \underset{S} 
\to {\varprojlim}\ \pi_1(\Bbb P^1_{\overline K}\setminus  S_{\overline K}, 
\eta)^p\to Q$ which restricts to $\phi_2$ on the factor $\pi _1(C,\eta)^{\geo,p}$, 
and to $\psi_2$ on the factor $N$. Moreover, $\psi:\pi_1(\Bbb P^1_{\overline K}\setminus  (S_{0})_{\overline K}, 
\eta)^p\to Q$ 
factors through
$\psi:\pi_1(\Bbb P^1_{\overline K}\setminus  (S_{0})_{\overline K}, 
\eta)^p$ for some set of closed points $S_0\subset \Bbb P^1_K$, and $S_0\cap C=\emptyset$.
The homomorphism $\psi$ corresponds [after eventually a finite extension 
of $K$] to a Galois cover $Y\to \Bbb P^1_K$
with Galois group $Q$, which induces naturally a finite Galois cover
$g:\Spec B\to \Spec R[[T]]$ with Galois group $Q$ (above the formal completion at $T=\infty$), 
which is by construction \'etale above the ideal 
$\pi$, and which dominates the Galois cover $f:\Spec A\to \Spec R[[T]]$ we 
started with. This in turn corresponds to a homomorphism $\psi:\Delta \to Q$
with the required properties.
\qed 
\enddemo

The author doesn't know, and is interested to know, the answer to the 
following question.

\definition {Questions 2.3.2}   Is the maximal pro-$p$ quotient $\Gal (\overline L/\overline K.L)^p$ of the 
[geometric] Galois group $\Gal (\overline L/\overline K.L)$ a free pro-$p$ group?
\enddefinition

\definition {2.4}
Next, we investigate a certain quotient of the "geometric absolute Galois group"
of the boundary of a $p$-adic open disc (see the definition of the profinite groups $\Delta '$, and $\Pi'$ below). 
\enddefinition

Let
$R[[T]]\{T^{-1}\} \defeq \{\sum _{i=-\infty}^{\infty}a_iT^i,\ \underset 
{i\to -\infty} \to  \lim \vert a_i \vert=0\}$ be as in 1.2.3.
Note that $R[[T]]\{T^{-1}\}$ is a complete discrete valuation ring, with uniformising 
parameter $\pi$, and residue field the formal power series field $k((t))$, where $t=T \mod \pi$. 
Write 
$$\Cal X\defeq \Spec R[[T]]\{T^{-1}\}.$$
[$\Cal X$ is what we shall refer to as the boundary of a $p$-adic open disc (over $K$)].
Let 
$$M\defeq \Fr (R[[T]]\{T^{-1}\})$$ 
be the quotient field of the discrete valuation ring $R[[T]]\{T^{-1}\}$. 

Assume that the generic point $\eta$ of $\Tilde X_K$ above arises from a generic point  $\eta$
of $R[[T]]\{T^{-1}\}\otimes _RK$. In particular, the generic point $\eta$ 
determines then an algebraic closure $\overline M$ of $M$. We have a natural exact sequence
of Galois groups
$$1\to \Gal (\overline M/M.\overline K)\to \Gal (\overline M/M)\to 
\Gal (\overline K/K)\to 1.$$

Let $I'\defeq I_{(\pi)}'\subset \Gal (\overline M/\overline K.M)$ be the normal 
closed subgroup which is generated by the inertia subgroups above the ideal 
$(\pi)$ of $R[[T]]\{T^{-1}\}$, which is generated by $\pi$. Write
$$\overline \Delta' \defeq \Gal (\overline M/\overline K.M)/I'.$$

Note that, by definition, the profinite group $\overline \Delta'$ classifies
finite Galois covers 
$\Cal Y_L\to \Cal X_L\defeq  \Cal X\times _RL,$ 
where $L$ is a finite extension of $K$ with valuation ring $R'$, $L$ is 
algebraically closed in $\Cal Y_L$, $\pi'$ is a uniformising parameter of $L$, 
and the natural morphism $\Cal Y\to \Cal X'\defeq \Cal X\times _RR'$ 
where $\Cal Y$ is the normalisation of $\Cal X'$ in $\Cal Y_L$ is 
\'etale. 

The natural morphism
$$\Spec R[[T]]\{T^{-1}\}\to \Spec R[[T]]$$
induces a natural homomorphism 
$$\overline \Delta ' \to \overline \Delta.$$ 

Let 
$$\Delta' \defeq {\overline \Delta' }^p$$ 
be the maximal pro-$p$ quotient of $\overline \Delta'$. We have a natural homomorphism
$$\Delta'\to \Delta.$$

Our next technical result in this section is the following.

\proclaim {Theorem 2.4.1} There exists a natural homomorphism $\Delta '\to \Delta$
which makes $\Delta'$ into a direct factor of the free pro-$p$ group
$\Delta$. In particular, $\Delta'$ is a free pro-$p$ group [this can also be deduced 
from the fact that the maximal pro-$p$ quotient of the absolute Galois group of the field
$k((t))$ is a free pro-$p$ group].
\endproclaim

\demo{Proof} One has to verify the cohomological criterion in Proposition 1.1.3
for being a direct factor. 

Let $f': \Cal Y\to \Cal X$
be an \'etale $\Bbb Z/p\Bbb Z$-torsor. One has to construct [eventually after a finite extension
of $K$] a finite generically Galois cover $f: \Tilde Y\to \Tilde X\defeq \Spec R[[T]]$ of degree $p$
which induces above $\Cal X$, by pull-back via the natural morphism
$\Cal X\defeq \Spec R[[T]]\{T^{-1}\}\to \Tilde X\defeq \Spec R[[T]]$,
the $\Bbb Z/p\Bbb Z$-torsor $f'$. 

The torsor $f'$ induces naturally a finite generically Galois cover 
$\bar f':\Cal Y_k\to \Cal X_k=\Spec k[[t]]$ of degree $p$.
There exists [as is easily verified (cf. also [Ka])] a finite Galois cover 
$\bar g:Y_k\to \Bbb P^1_k$ of degree $p$ which is ramified above a unique point 
$\infty \in \Bbb P^1_k$, and such that the Galois cover induced by 
$\bar g$ above the formal completion $\Spec \hat {\Cal O}_{\Bbb P^1_k,\infty}$ 
of $\Bbb P^1_k$ at $\infty$ is isomorphic to $\bar f'$. Let
$\bar g':Y_k'\to \Bbb A^1_k\defeq \Bbb P^1_k\setminus \{\infty\}$ be the restriction
of $\bar g$ which is an \'etale $\Bbb Z/p\Bbb Z$-torsor above $\Bbb A^1_k$. The \'etale torsor $\bar g'$
lifts [uniquely up to isomorphism] to an \'etale $\Bbb Z/p\Bbb Z$-torsor
$g':Y'\to D_{0}\defeq \Sp K<S>$ [where $D_0$ is the closed disc centered at $S=0$],
by the theorems of liftings of \'etale covers (cf. [Gr]), 
whose restriction $\tilde g':\tilde Y'\to C$ to the annulus $C$ is an \'etale 
$\Bbb Z/p\Bbb Z$-torsor, which corresponds to a continuous homomorphism 
$\psi:\pi_1(C,\eta)^{\geo,p} \twoheadrightarrow \Bbb Z/p\Bbb Z$.

The geometric fundamental group
 $\pi_1(C,\eta)^{\geo,p}$ being a direct factor of $\underset{S} 
\to {\varprojlim}\ \pi_1(\Bbb P^1_{\overline K}\setminus  S_{\overline K}, 
\eta)^p$, the above homomorphism $\psi$ arises [by restriction] from a
continuous homomorphism $\psi':  \underset{S} 
\to {\varprojlim}\ \pi_1(\Bbb P^1_{\overline K}\setminus  S_{\overline K}, 
\eta)^p\twoheadrightarrow \Bbb Z/p\Bbb Z$, and the later gives rise naturally to
a Galois cover $g:Y\to \Bbb P^1_K$ of degree $p$
[this cover only exists a priory over a finite extension of $K$
but we can, without loss of generality, assume that it is defined over $K$]
whose restriction to the annulus $C$ is isomorphic [by construction]
to the above Galois cover $\Tilde g'$. The Galois cover $g$ induces naturally
a Galois cover above $\Tilde X\defeq \Spec R[[T]]$ (i.e. above the formal completion at $T=\infty$),
which induces above the boundary $\Cal X$ the torsor $f'$ as required.
\qed  
\enddemo

\definition {2.5}
In [Ga], Garuti investigated the problem of lifting of Galois covers between 
smooth curves.
In this sub-section we will prove a refined version of the main result in [Ga], using Theorem 2.3.1 and Theorem 2.4.1.
\enddefinition

First, we recall the following main result of Garuti.

\proclaim {Theorem 2.5.1 (Garuti)} Let $X$ be a proper, smooth, and geometrically connected
$R$-curve. Let 
$$f_k:Y_k\to X_k\defeq X\times _Rk$$ 
be a finite [possibly ramified] Galois cover between smooth $k$-curves  
with Galois group $G$. Then there exists a finite extension $R'/R$ and a finite
morphism 
$$f':Y'\to X'\defeq X\times _RR',$$ 
which is generically \'etale and Galois with Galois group $G$, satisfying the 
following properties:

\noindent
(i)\ $Y'$ is a proper and normal $R'$-curve.

\noindent
(ii)\ The natural morphism $f'_k:Y'_k\defeq Y'\times _Rk\to X'_k=X_k$ is generically 
\'etale and Galois with Galois group $G$. Moreover, there exists a $G$-equivariant 
birational morphism $\nu:Y_k\to Y'_k\defeq Y'\times _Rk$ such that the following 
diagram is commutative:
$$
\CD
Y_k  @>{\nu}>>   Y_k'  \\
@V{f_k}VV       @V{f'_k}VV  \\
X_k   @>{\id_{X_k}}>>   X_k \\
\endCD
$$
and the morphism $\nu$ is an isomorphism outside the divisor of ramification 
in the morphism $f_k:Y_k\to X_k$.

\noindent
(iii)\ The special fiber $Y'_k$ is reduced, unibranche, and the morphism
$\nu:Y_k\to Y'_k$ is a morphism of normalisation. In particular, $Y_k$ and $Y'_k$
are homeomorphic.
\endproclaim

\demo {Proof} (cf. [Ga], Proof of Th\'eor\`eme 2).
\qed
\enddemo

In light of the above result, we define Garuti liftings as follows.

\definition {Definition 2.5.2 (Garuti Liftings of Galois Covers between Smooth Curves)} 
Let $X$ be a proper, smooth, and geometrically connected $R$-curve. Let 
$$f_k:Y_k\to X_k\defeq X\times _Rk$$ 
be a finite [possibly ramified] Galois cover with Galois group $G$. Let
$$f':Y'\to X'\defeq X\times _RR'$$ 
be as in Theorem 2.5.1, for some finite extension $R'/R$. We call $f'$ a Garuti 
lifting of the Galois cover $f_k$ [defined over $R'$]. 

We say that $f'$ is a smooth lifting of $f_k$, if 
$Y'$ is a smooth $R$-curve, which is equivalent to the above morphism
$\nu: Y_k\to Y'_k$ being an isomorphism. 

Note that, by definition, a Garuti lifting is 
defined [a priori] over a finite extension of $R$. Also, if $f_k$ is \'etale,
then a smooth lifting of $f_k$ always exists over $R$ as follows from the theorems
of liftings of \'etale covers (cf. [Gr]).
\enddefinition

The following Theorem is a refined version of the above result of Garuti.

\proclaim {Theorem 2.5.3} Let $X$ be a proper, smooth, and geometrically connected $R$-curve. Let
$$f_k:Y_k\to X_k\defeq X\times _Rk$$ 
be a finite [possibly ramified] Galois cover with Galois group $G$ between 
smooth $k$-curves. Assume that the finite group $G$ sits in an exact sequence
$$1\to H'\to G\to H\to 1.$$ 
Let 
$$Y_k @>{g_k}>> Z_{k} @>{h_k}>> X_k$$ 
be the corresponding factorisation of 
the Galois cover $f_k$. Thus, $h_k:Z_{k}\to X_k$ is a finite Galois cover with
Galois group $H$ between smooth $k$-curves. Let 
$$h':Z'\to X'\defeq X\times _RR'$$ 
be a  Garuti lifting of the Galois cover $h_k$ defined over the finite extension $R'/R$ (cf. Definition 2.5.2). 

Then there exists a finite extension $R''/R'$, and a Garuti lifting 
$$f'':Y''\to X''\defeq X\times _RR''$$ 
of the Galois cover $f_k$ over $R''$, which dominates $h'$, i.e. we have a factorisation 
$$f'':Y''@>{g''}>> Z''\defeq Z\times _{R'}R'' @>{h''\defeq h'\times _{R'}R''}>> X'',$$
where $g'':Y''\to Z''$ is a finite morphism between normal $R''$-curves.
\endproclaim

\demo{Proof}
The proof is [in some sense] similar to the proof of Th\'eor\`eme 2 in [Ga] using the 
above Theorem 2.4.1. More precisely, using the techniques of formal paching 
(cf. [Ga], and Proposition 1.2.2)
the proof of Theorem 2.5.3 follows directly from the following local result in Theorem 2.5.5.
\qed
\enddemo

Before stating our main local result, we first define the local analog of Garuti liftings.

\definition {Definition 2.5.4 (Local Garuti Liftings)} Let 
$\Tilde {X}\defeq \Spec R[[T]]$, and  $\Tilde {X}_k\defeq \Spec k[[t]]$.  
Let $G$ be a finite group and 
$$f_k:\Tilde {Y}_k\to \Tilde {X}_k$$ 
a finite morphism, which is generically Galois with Galois group $G$, with $\Tilde {Y}_k$ connected 
and normal. We call a Garuti lifting of the Galois cover $f_k$, over the finite extension $R'/R$, a finite Galois cover
$$f':\Tilde {Y}'\to \Tilde {X}'\defeq X\times _R R'$$
with Galois group $G$, where $R'/R$ is a finite extension, 
the morphism $f'_k:\Tilde {Y}'_k\to 
\Tilde {X}_k$ is generically Galois with Galois group $G$, there exists a birational 
$G$-equivariant morphism $\nu: \Tilde {Y}_k\to \Tilde {Y}'_k$ which is a morphism of normalisation,
and a factorisation 
$$f_k:\Tilde {Y}_k@>{\nu}>> \Tilde {Y}'_k@>{f'_k}>> \Tilde {X}_k.$$

Moreover, we say that $f$ is a smooth lifting of $f_k$ if $\Tilde Y'$ is smooth over $R'$, or equivalently 
if the above morphism $\nu$ is an isomorphism.
\enddefinition

The following is our main result which is a refined version of the local version of Garuti's main Theorem 2.5.1.

\proclaim{Theorem 2.5.5} Let $\Tilde X\defeq \Spec R[[T]]$, and $\Tilde X_k\defeq 
\Spec k[[t]]$.  Let $G$ be a finite group and 
$$f_k:\Tilde Y_k\to \Tilde X_k$$ 
a finite morphism 
which is generically Galois with Galois group $G$, with $\Tilde Y_k$ normal and connected. 
Let $H$ be a quotient of $G$ and 
$$h_k:\Tilde Z_k\to \Tilde  X_k$$ 
the corresponding Galois 
sub-cover with Galois group $H$. Let 
$$h':\Tilde Z'\to \Tilde  X'\defeq \Tilde X\times _R R'$$ 
be  a Garuti lifting of $h_k$ over a finite extension $R'/R$ (cf. Definition 2.5.4). 
Then there exists  a finite extension $R''/R'$,
and a Garuti lifting 
$$f'':\Tilde Y''\to \Tilde X''\defeq \Tilde X\times _R R''$$ 
of $f_k$ over $R''$ which dominates $h'$, i.e. we have a 
factorisation: 
$$f'':\Tilde Y''\to \Tilde Z''\defeq Z'\times _{R'}R''@>{h''\defeq h'\times _{R'}R''}>> \Tilde  X''.$$
\endproclaim

\demo{Proof}  The Galois group $G$ is a solvable group which is a semi-direct
product of a cyclic group of order prime to $p$ by a $p$-group.
 By similar arguments as the ones used by Garuti in [Ga],
 it suffices to treat the case where $G$ is a $p$-group (see the arguments used in [Ga],
 Th\'eor\`eme 2.13, and Corollaire 1.11). In this case the proof follows from Theorem 2.4.1.
 
 More precisely, assume that $G$ is a $p$-group [hence $H$ is also a $p$-group].
 The Galois cover $f_k:\Tilde Y_k\to \Tilde X_k$ is generically given by an \'etale Galois cover
 $\Spec k((s))\to \Spec k((t))$ with Galois group $G$. This \'etale cover 
 lifts uniquely to an \'etale Galois cover $\Cal Y\to \Cal X\defeq \Spec R[[T]]\{T^{-1}\}$
 above the boundary of the open disc $\Tilde X$, which is Galois with Galois group $G$,
 and which corresponds to a continuous homomorphism $\psi_2: \Delta'\to G$.
 
 Let $N$ be a complement of $\Delta '$ in $\Delta$ (cf. Theorem 2.4.1).
 The local Garuti lifting $h':\Tilde Z'\to \Tilde  X'\defeq \Tilde X\times _R R'$
 corresponds to a continuous homomorphism $\phi: \Delta \to H$, which restricts to continuous 
 homomorphisms $\psi_1:\Delta'\to H$, and $\phi_1:N\to H$. The above homomorphism
 $\psi_2$ dominates by construction the homomorphism $\psi _1$. The pro-$p$ group $N$ being free
 one can lift the homomorphism $\phi_1$ to a continuous homomorphism  $\phi_2:N\to G$ 
 which dominates $\phi_1$. The pro-$p$ group $\Delta$ being isomorphic to the direct free product
 $\Delta' \star N$, both $\psi_2$ and $\phi_2$ give rise to a continuous homomorphism
 $\phi':\Delta \to G$ which dominates the above morphism $\phi$. The homomorphism
 $\phi'$ in turn corresponds to a Galois cover $\Tilde Y''\to \Tilde X''\defeq \Tilde X\times _R R''$
 over some finite extension $R''/R$, which is a Garuti lifting of $f_k:\Tilde Y_k\to \Tilde X_k$, and
 which by construction dominates the Garuti lifting $h':\Tilde Z'\to \Tilde  X'\defeq \Tilde X\times _R R'$
 of the sub-cover $h_k:\Tilde Z_k\to \Tilde  X_k$ as required.
 \qed
\enddemo

\definition{Remark 2.5.6} We assumed in this section that $R$ is of unequal characteristics. In fact 
the main results of this section: Theorem 2.3.1, Theorem 2.4.1, Theorem 2.5.3, and Theorem 2.5.5, are also valid 
in the case of a complete discrete valuation ring $R$ of equal characteristics $p>0$. Indeed,
the result of Garuti  (cf. Proposition 2.2.3)  that we use in the proof of Theorem 2.3.1, and Theorem 2.4.1, is valid in this case (cf. [Ga]).
\enddefinition

\subhead
\S 3. Fake Liftings of Cyclic Covers between Smooth Curves
\endsubhead
In this section we use the same notations as in $\S 2$, 2.1. We will 
investigate the problem of lifting of cyclic [of p-power order] Galois 
covers between smooth curves.

\definition {3.1 The Oort Conjecture}
 First, we recall the following main conjecture 
which was formulated by F. Oort, and several of its variants.
In what follows $R$ is as in the Notations 2.1.
\enddefinition

\definition  {The Original Oort Conjecture [Conj-O]} (cf. [Oo], and [Oo1]) Let
$$f_k:Y_k\to X_k$$
be a finite [possibly 
ramified] Galois cover between smooth $k$-curves, with Galois group 
$G\isom \Bbb Z/m\Bbb Z$
a cyclic group. Then there exists a finite extension $R'/R$, and
a Galois cover
$$f:Y'\to X'$$
with Galois group $G$, where $X'$ and $Y'$ are smooth $R'$-curves,
which lifts the Galois cover $f_k$,
i.e. the morphism induced by $f$ at the level of special fibres is [Galois] isomorphic to $f_k$.
\enddefinition

In the original version of the conjecture, one doesn't fix $R$, but fixes $k$, $f_k$, and asks for the existence of a local domain $R$ dominating 
the ring of Witt vectors $W(k)$, over which a lifting of $f_k$ exists, 
as part of the conjecture (cf. [Oo]).

One can formulate several variants of the above conjecture, that we will list below.

\definition {[Conj-O1]} Let $X$ be a proper, smooth, geometrically connected
$R$-curve, and 
$f_k:Y_k\to X_k\defeq X\times _Rk$
a finite [possibly 
ramified] Galois cover between smooth $k$-curves, with Galois group 
$G\isom \Bbb Z/m\Bbb Z$
a cyclic group. Then there exists a smooth lifting of $f_k$ (cf. Definition 2.5.2), i.e. there exists
a finite extension $R'/R$, and a Galois cover 
$f':Y'\to X'\defeq X\times _RR'$ 
between smooth $R'$-curves,  
with Galois group $G$, such that the special fiber $X'_k\defeq X'\times _Rk$ 
(resp. $Y'_k\defeq Y'\times _Rk$) equals $X_k$ (resp. is isomorphic to $Y_k$), and the 
natural morphism $f'_k \defeq f'\times _Rk:Y'_k\to X'_k=X_k$ which is induced by $f'$
on the level of special fibres is 
isomorphic to $f_k$. 
\enddefinition

\definition {[Conj-O2]} Let 
$f_k:Y_k\to \Bbb P^1_k$ be
a finite 
ramified Galois cover, with $Y_k$ a smooth $k$-curve , and Galois group 
$G\isom \Bbb Z/m\Bbb Z$
a cyclic group. Then there exists a smooth lifting of $f_k$ (cf. Definition 2.5.2), i.e. there exists
a finite extension $R'/R$, a finite Galois cover 
$f':Y'\to \Bbb P^1_{R'}$, 
with $Y'$ a smooth $R'$-curve, 
with Galois group $G$, and such that the 
natural morphism $f'_k \defeq f'\times _Rk:Y'_k\to \Bbb P^1_k$ which is induced by $f'$
on the level of special fibres is 
isomorphic to $f_k$. 
\enddefinition

\definition {[Conj-O3]} Let 
$f_k:Y_k\to \Bbb P^1_k$
be a finite Galois cover, with $Y_k$ a smooth $k$-curve, and Galois group 
$G\isom \Bbb Z/m\Bbb Z$
a cyclic group, 
which is [totally] ramified above a unique point $\infty \in \Bbb P^1_k$.
Then there exists a smooth lifting of $f_k$ (cf. Definition 2.5.2), i.e. there exists
a finite extension $R'/R$, a finite Galois cover 
$f':Y'\to \Bbb P^1_{R'}$, 
with $Y'$ a smooth $R'$-curve, 
with Galois group $G$, and such that the 
natural morphism $f'_k \defeq f'\times _Rk:Y'_k\to \Bbb P^1_k$ which is induced by $f'$
on the level of special fibres is 
isomorphic to $f_k$. 
\enddefinition

\definition {[Conj-O4]} Let $\Tilde X\defeq \Spec R[[T]]$, and $\Tilde X_k\defeq 
\Spec k[[t]]$.  Let $f_k:\Tilde Y_k\to \Tilde X_k$ be 
a finite morphism 
which is generically Galois with Galois group $G\isom \Bbb Z/m\Bbb Z$ a cyclic group, 
with $\Tilde Y_k$ normal and connected. 
Then there exists  a finite extension $R'/R$,
and a smooth lifting 
$f':\Tilde Y'\to \Tilde X'\defeq \Tilde X\times _R R'$ 
of $f_k$, i.e. $\Tilde Y'\isom \Spec R'[[T']]$ is $R'$-smooth, and the natural
morphism $f'_k: \Tilde Y'_k\to \Tilde X'_k=\Tilde X_k$ which is induced by $f'$ at the level of special fibres is
isomorphic to $f_k$. 
\enddefinition

Moreover, in the above conjectures {\bf [Conj-O1]},  {\bf [Conj-O2]},   {\bf [Conj-O3]},   and {\bf [Conj-O4]}, 
one predicts that $R'$ can be chosen to be the minimal extension of $R$ which contains a primitive $m$-th root of $1$. 

In fact all the above variants of the Oort conjecture turn out to be equivalent. More precisely, we have the following.

\proclaim {Lemma 3.1.1} With the above notations, the various conjectures
{\bf [Conj-O]},  {\bf [Conj-O1]},  {\bf [Conj-O2]},   {\bf [Conj-O3]},   and {\bf [Conj-O4]}, are all equivalent.     
Moreover, in order to solve the above conjecture(s), it suffices to treat the case where $G\isom \Bbb Z/p^n\Bbb Z$ is a cyclic $p$-group.
\endproclaim

\demo {Proof} Follows easily from the local-global principle for the lifting of Galois covers between curves (cf. Proposition 1.2.4), the result
of approximation of local extensions by global extensions du to Katz, Gabber, and Harbater, (cf. [Ha], and [Ka]), and the formal patching result in Proposition
1.2.2. The last assertion can also be easily verified (see for example the arguments in [Gr-Ma], 6). 
\qed
\enddemo

Oort conjecture holds true in the case where the Galois cover $f_k$ is \'etale, as 
follows from the theorems of liftings of \'etale covers (cf. [Gr]). In this case
the statement of the conjecture is true for any finite group $G$ [not necessarily cyclic],
and a smooth lifting exists over $R$. In the case where $G\isom \Bbb Z/p^n\Bbb Z$
is a cyclic $p$-group, the conjecture has been verified in the cases 
where $n=1$, and $n=2$ (cf. [Se-Oo-Su] for the case $n=1$, and [Gr-Ma] for the case $n=2$).

In this paper, and in light of Theorem 2.5.3,
we propose the following refined version of the Oort conjecture. More precisely, we will formulate a
refined version of {\bf [Conj-O1]}, which is equivalent to {\bf [Conj-O]} by Lemma 3.1.1.

\definition {Oort Conjecture Revisited [Conj-O1-Rev]}
Let $X$ be a proper, smooth,  geometrically connected $R$-curve, and 
$f_k:Y_k\to X_k\defeq X\times _Rk$ 
a finite [possibly ramified] Galois cover between 
smooth $k$-curves, with Galois group 
$G\isom \Bbb Z/m\Bbb Z$ a cyclic group. 
Let $H$ be a quotient of $G$, and $g_k:Z_k\to X_k$
the corresponding Galois sub-cover of $f_k$ with Galois group $H$. 
Then there exists a smooth Galois lifting 
$$g:Z'\to X'\defeq X\times _RR'$$ 
of $g_k$, 
over some finite extension $R'/R$ [i.e. $g$ is a Galois cover with Galois group $H$
between smooth $R'$-curves which is a lifting of $g_k$]. 

Furthermore, for every smooth lifting $g$ of the Galois sub-cover $g_k$ of $f_k$ as above,
there exists a finite extension $R''/R'$, and a
finite Galois cover
$$f:Y''\to X''\defeq X\times _RR''$$ 
between smooth $R''$-curves, 
with Galois group $G$, which is a smooth lifting of $f_k$ (cf. Definition 2.5.2), and
such that $f$ dominates $g$, i.e. we have a factorisation 
$$f:Y''\to Z''\defeq Z'\times _{R'}R''@>{g \times _{R'}R''}>> X''.$$ 

Moreover, $R''$ can be chosen to be the minimal extension 
of $R'$ which contains a primitive $m$-th root of $1$. 
\enddefinition

\definition {Remark 3.1.2}  In a similar way, one can revisit the above [equivalent] variants of the original Oort 
conjecture, and formulate the revisited versions 
{\bf [Conj-O2-Rev]},   {\bf [Conj-O3-Rev]},   and {\bf [Conj-O4-Rev]},
which turn out to be all equivalent to {\bf [Conj-O1-Rev]} (use similar arguments as in the proof of Lemma 3.1.1). 
Moreover, and in order to solve these revisited versions, one can reduce to the case 
where $G\isom \Bbb Z/p^n\Bbb Z$ is a cyclic $p$-group.
In the case where $n=1$ [i.e. $G$ is a cyclic group of cardinality $p$] 
the revisited Oort conjecture is clearly true, since the [original] Oort conjecture 
is true in this case (see [Se-Oo-Su]). Both the original and the revisited conjectures are clearly equivalent in this case.
\enddefinition

\definition {3.2}  
Next, we give examples where the revisited Oort conjecture can be verified in the case where
$G\isom \Bbb Z/p^2\Bbb Z$.
\enddefinition

We assume that $G\isom \Bbb Z/p^2\Bbb Z$ is cyclic of order $p^2$. 
We will work within the framework of {\bf [Conj-O4-Rev]} (cf. Remark 3.1.2).
More precisely, let $\Tilde X\defeq \Spec R[[T]]$, and $\Tilde X_k\defeq \Spec k[[t]]$
its special fiber [where $t= T \mod \pi$]. Let 
$$f_k:Y_k\to \Tilde X_k$$
be a cyclic Galois cover of degree $p^2$, with $Y_k$ normal,  and 
$$h_k:Y'_k\to \Tilde X_k$$ 
its unique Galois sub-cover of degree $p$.
A smooth local lifting of $f_k$ [cf. Definition 2.5.4] exists by [Gr-Ma], Theorem 5.5,  over $R$ if $R$ contains the $p^2$-th roots of $1$.
From now on we will assume in this sub-section that $R$ contains a primitive $p^2$-th root of $1$.
Let 
$$h:Y'\to \Tilde X$$ 
be a smooth Galois lifting of  $h_k$, i.e.  
$h$ is a Galois cover of degree $p$,
$Y'\defeq \Spec A'$, $A'\isom 
R[[Z']]$ is an open disc, and $h$ induces the Galois cover $h_k$ on the level of special fibres. 
Then, in order to verify the {\bf [Conj-O4-Rev]} for the Galois cover $f_k$ and the smooth lifting $h$ of $h_k$, 
it suffices to show that there exists a smooth Galois lifting 
$$f:Y\to \Tilde X$$ 
of $f_k$, i.e. $f$ is a cyclic Galois cover of degree $p^2$, 
 $Y\defeq \Spec A$, $A\isom R[[Z]]$ is an open disc, and $f$ induces the Galois cover $f_k$
 on the level of special fibres,
 which dominates $h$: 
i.e. such that we have a factorisation 
$$f:Y\to Y' @>h>> \Tilde X.$$ 

The Galois cover $f_k$ is generically given, 
for an appropriate choice of the parameter $t$, by the equations:

$$x_1^p-x_1=t^{-m_1},\tag *$$
and
$$x_2^p-x_2=c(x_1^p,-x_1)+\sum _{0\le s < m_1(p-1)} a_st^{-s}+\sum
_{0\le j<m_1}t^{-jp}\sum _{0<i<p}(x_1^p-x_1)^ip_{j,p-i}(x_1^p-x_1)^p, \tag **$$
where $a_i\in k$, $p_{j,p-i}\in k[x]$ are polynomials of respective degrees
$d_{j,p-i}$, $\gcd (m_1,p)=1$, and 
$$c(x,y)\defeq \frac {(x+y)^p-x^p+(-y^p)}{p}.$$
(See [Gr-Ma], Lemma 5.1).
Moreover, the degree of the different in the Galois cover $f_k$ is
$$d_s\defeq (m_1+1)(p-1)p+(m_2+1)(p-1),$$
where
$$m_2\defeq \underset {0<i<p}\to {\max _{0\le j<m_1}}(p^2m_1,p(jp+(i+pd_{j,p-i})m_1))-(p-1)m_1,$$
(cf. loc. cit).

Let $\zeta _2\in R$ be a primitive $p^2$-th root of $1$. Let 
$\zeta_1\defeq \zeta _2^p$, and $\lambda \defeq \zeta_1-1$.
The smooth lifting $h:Y'\to \Tilde X$ 
of $h_k$  is generically given [by the Oort-Sekiguchi-Suwa theory (cf. [Se-Oo-Su])] by an equation
$$\frac {(\lambda X_1+1)^p-1}{\lambda\ ^p}=f(T),$$
where 
$$f(T)=\frac {h(T)}{g(T)},$$
$h(T)\in R[[T]]$, $g(T)\in R[T]$ is a distinguished  polynomial [i.e. its highest coefficient is a unit in $R$],
the degree of $g(T)$ is $m$, the Weierstrass degree of $h(T)$ is $m'$, $m\ge m'$, and  $m-m'=m_1$.
Furthermore, 
$$\frac {h(T)}{g(T)}= T^{-m_1}\mod \pi.$$

The smoothness of $Y'$ is equivalent, by the local criterion for smoothness (cf. [Gr-Ma], 3.4), 
to the fact that the Galois cover $h_K:Y'_K\to \Tilde X_K$ which is induced by $h$ between generic fibres, 
and which is given by the equation
$$(\lambda X_1+1)^p=\frac {\lambda ^ph(T)+g(T)}{g(T)},$$
is ramified above $m_1+1$ distinct geometric points of $\Tilde X_K$. 
Moreover,  $Z'\defeq X_1^{-\frac {1}{m_1}}$
is a parameter for the open disc $Y'$, as follows easily from arguments similar to the one 
given in the proof of Theorem 4.1 in [Gr-Ma] (cf. also [Gr-Ma], proof of 3.4).

We will consider two cases, depending on the lift $h$ of $h_k$,  
where we can prove the revisited Oort conjecture {\bf [Conj-O4-Rev]} for the smooth lifting
$h:Y'\to \Tilde X$ [i.e. we can dominate $h$ by a smooth lifting $f$ of $f_k$]. 
These two cases are considered separately in the following lemmas
3.2.1 and 3.2.2.

\proclaim {Lemma 3.2.1} With the same notations as above. Assume that
in the second equation $(**)$ above defining the Galois cover $f_k$
we have 
$$\sum _{0\le s < m_1(p-1)} a_st^{-s}+\sum
_{0\le j<m_1}t^{-jp}\sum _{0<i<p}(x_1^p-x_1)^ip_{j,p-i}(x_1^p-x_1)^p=0,$$ 
and also assume that the degree of
$g(T)$ above equals $m_1$. [In particular, $h(T)\in R[[T]]$ above is a unit in this case]. 
Then there exists a smooth lifting $f$ 
of $f_k$ which dominates the smooth lifting $h$ of $h_k$.
In particular, {\bf [Conj-O4-rev]} is true under these conditions for the Galois cover $f_k$, and the 
smooth lifting $h$ of the sub-cover $h_k$.
\endproclaim

\demo{Proof}
Consider the cover
$$f:Y\to \Tilde X$$ 
which is generically given by the equations
$$\frac {(\lambda X_1+1)^p-1}{\lambda\ ^p}=f(T),\tag i$$
where $f(T)={h(T)}/{g(T)}$ is as above, and
$$(\lambda X_2+\Exp_p (\mu X_1))^p=(\lambda X_1+1)\Exp _p(\mu^pY),\tag ii$$ 
where 
$$\Exp_p X\defeq 1+X+...+\frac {X^{p-1}}{(p-1)!}$$ 
is the truncated exponential, 
$$\mu\defeq \log_p(\zeta_2)=1-\zeta _2+......+(-1)^{p-1}\frac {\zeta_2^{p-1}}{p-1},$$ 
[$\Exp_p$ and $\log_p$ denote the truncation of the exponential and the logarithm, 
respectively, by terms of degree $>p-1$], 
and
$$Y\defeq \frac {(\lambda X_1+1)^p-1}{\lambda ^p}=\frac {h(T)}{g(T)}.$$
Then $f$ is a cyclic Galois cover of degree $p^2$ which lifts the Galois cover $f_k$
(cf. [Gr-Ma], the discussion in the beginning of 3, and Lemma 5.2). 

We claim that 
$Y$ is smooth over $R$. Indeed, the degree of the different in the morphism $f_k:Y_k\to \Tilde X_k$ in this case
is 
$$d_s=(m_1+1)(p-1)p+(p^2m_1-(p-1)m_1+1)(p-1).$$ 
Moreover, the above second equation (ii) defining the lifting $f$ is
$$(X_2')^p=(\lambda X_1+1)\Exp _p(\mu^pY)=(1+\lambda X_1)(1+\mu^p\frac {h(T)}{g(T)}+....+\frac {\mu^{p(p-1)}}{(p-1)!}\frac {h(T)^{(p-1)}}{g(T)^{(p-1)}})$$
and
$$1+\mu^p\frac {h(T)}{g(T)}+....+ \frac {\mu^{p(p-1)}}{(p-1)!}    \frac {h(T)^{(p-1)}}{g(T)^{(p-1)}} $$

equals
$$\frac {(p-1)!g(T)^{p-1}+\mu^p(p-1)! h(T)g(T)^{p-2}+....+\mu^{p(p-1)}h(T)^{(p-1)}}     {(p-1)!g(T)^{p-1}}.$$
Furthermore, 
$$(p-1)!g(T)^{p-1}+\mu^p(p-1)! h(T)g(T)^{p-2}+....+\mu^{p(p-1)}h(T)^{(p-1)}$$
can be written as a series in $X_1^{-\frac {1}{m_1}}$, whose Weierstrass degree is
$pm_1(p-1)$ [since we assumed the degree of $g(T)$ to be $m_1$]. 
From this we deduce that the degree of the generic different $d_{\eta}$ in the cover $f_K:Y_K\to \Tilde X_K$ satisfies
$$d_{\eta}\le (m_1+1)(p^2-1)+pm_1(p-1)^2,$$ 
which implies $d_{\eta}\le d_s$. One then concludes that $d_{\eta}=d_s$, hence that $Y$ is smooth over $R$,  
since in general we must have $d_s\le d_{\eta}$. Moreover, we have [by construction] a natural
factorisation $g:Y\to Y'@>h>> \Tilde X$.
\qed
\enddemo

\proclaim {Lemma 3.2.2} With the same notations as above. Assume that $g(T)=T^{m_1}$. [Thus, in particular, $h(T)\in R[[T]]$ is a unit]. 
[This case is rather special, since the corresponding smooth lifting $h$ of the Galois sub-cover $h_k$ has the property that
all branched points are equidistant in the $p$-adic topology of $K$]. Then there exists a smooth lifting $f$ 
of $f_k$ which dominates the smooth lifting $h$ of $h_k$. In particular, {\bf [Conj-O4-rev]} is true under these 
conditions for the Galois cover $f_k$, and the 
smooth lifting $h$ of the sub-cover $h_k$.
\endproclaim

\demo{Proof}
Consider the lifting 
$$f:Y\to \Tilde X$$ 
of the Galois 
cover $f_k:Y_k\to \Tilde X_k$,  which is generically given by the equations
$$\frac {(\lambda X_1+1)^p-1}{\lambda\ ^p}=f(T),\tag i'$$
where 
$$f(T)=\frac {h(T)}{T^{m_1}}$$ 
satisfies the above condition in this Lemma, and
$$[\lambda X_2+\Exp_p (\mu X_1)(1+\underset {0<i<p} \to {\sum _{0\le j<m_1}}T^{-j}\mu^i(p-i)!P_{j,p-i}
(g(T)))]^p\tag ii'$$
$$= (G(T^{-1})+p\mu^p\sum _{0<s<r}A_sT^{-s})(\lambda X_1+1),$$ 
where 
$$\Exp_p X\defeq 1+X+...+\frac {X^{p-1}}{(p-1)!}$$ 
is the truncated exponential, and
$$\mu\defeq \log_p(\zeta_2)=1-\zeta _2+......+(-1)^{p-1}\frac {\zeta_2^{p-1}}{p-1},$$ 
are as in the proof of Lemma 3.2.1 above, the polynomial 
$$G\defeq G(\frac {(\lambda X_1+1)^p-1}{\lambda\ ^p})$$ 
is defined in a similar way as in [Gr-Ma], Lemma 5.4,
$P_{j,p-i}\in R[X]$ are primitive polynomials which lift the $p_{j,p-i}\in k[x]$, 
and $A_s\in R$ lift the $a_s$ (cf. loc. cit). Then $f:\Tilde Y\to X$ is a Galois cover
with a cyclic Galois group [isomorphic to $\Bbb Z/p^2\Bbb Z$] and 
$Y$ is smooth over $R$, as follows from the local criterion
for good reduction (cf. [Gr-Ma], 3.4), by using Lemma 5.4 in [Gr-Ma] 
[where among others the degree of $G$ in $T^{-1}$ is computed], and the same argument
as in the proof of Theorem 5.5 in loc. cit. [The key points here are that $X_1^{-\frac {1}{m_1}}$
is a parameter for the disc $Y'$, and the key Lemma 5.4 in [Gr-Ma] is valid by replacing
$G\defeq G(T^{-m_1})$ there by $G\defeq G(f(T))$ in our case (formally speaking only the degree in $T^{-1}$ of 
$f(T)$, which is $m_1$, plays a role in loc. cit)].  Moreover, we have [by construction] a natural
factorisation $g:Y\to Y'@>h>> \Tilde X$.
\qed
\enddemo

\definition {3.3} Next, we will introduce the notion of fake liftings of cyclic Galois covers
between curves. We will work within the framework of {\bf [Conj-O2-Rev]}. 
\enddefinition
Let $n\ge 1$ be a positive integer.  Let
$$f_k:Y_k\to \Bbb P^1_k$$ 
be a finite ramified Galois cover, where $Y_k$ is a smooth $k$-curve, 
with Galois group $G\isom \Bbb Z/p^n\Bbb Z$
a cyclic $p$-group with cardinality $p^n$. 
We denote by 
$$g_k:X_k\to \Bbb P^1_k$$  
the unique sub-cover of $f_k$ which is Galois with Galois group 
$$H\isom \Bbb Z/p^{n-1}\Bbb Z.$$
We have a canonical factorisation 
$$f_k: Y_k @>h_k>> X_k @>g_k>> \Bbb P^1_k,$$ 
where $h_k: Y_k\to X_k$ is a cyclic Galois cover between smooth $k$-curves
of degree $p$.  

We assume that 
the Galois cover $g_k:X_k\to \Bbb P^1_k$ can be lifted to a Galois cover
between smooth $R$-curves [in other words admits a smooth lifting 
over $R$ (cf. Definition 2.5.2)], i.e. there exists a finite Galois cover
$$g:\Cal X\to \Bbb P^1_R$$
with Galois group $H$, where $\Cal X$ is smooth over $R$, $\Cal X_k \defeq \Cal X\times _Rk$ 
is isomorphic to $X_k$, and such that the morphism induced by $g$ at the level of 
special fibers 
$$g_k:\Cal X_k\to \Bbb P^1_k,$$ 
is isomorphic to the Galois cover $g_k:X_k\to \Bbb P^1_k$. 

By Theorem 2.5.3 
there exists a Garuti lifting (cf. Definition 2.5.2) of the Galois cover $f_k$ which dominates $g$.  
We assume [for simplicity] that such a Garuti lifting is defined over $R$, 
i.e. there exists a finite Galois cover
$$\tilde f:\Cal Y\to \Bbb P^1_R$$ 
with Galois group $G$,  and $\Cal Y$ normal, which dominates $g$, i.e. we have a factorisation 
$$\tilde f:\Cal Y @>{\tilde h}>>  \Cal X @>{g}>> \Bbb P^1_R,$$ 
and such  that the morphism 
$$\tilde f_k:\Cal Y_k\defeq \Cal Y\times _Rk \to \Bbb P^1_k$$ 
between special fibers is generically \'etale, Galois with Galois group
$G$, dominates $g_k$ [i.e. we have a factorisation $\tilde f_k:\Cal Y_k\to X_k
@>{g_k}>> \Bbb P^1_k$], the normalisation $\Cal Y_k^{\nor}$ of $\Cal Y_k$ is 
isomorphic to $Y_k$ [in particular, $\Cal Y_k$ is irreducible], 
and the natural morphism between the normalisations 
$$\Cal Y_k^{\nor} \to \Bbb P^1_k$$ 
[which is Galois] is isomorphic to $f_k$. 

Let $\delta _{\eta}\defeq \delta_{\tilde f_K}$ (resp. $\delta _s\defeq \delta_{f_k}$) be the degree 
of the different in the morphism $\tilde f_K:\Cal Y_K\defeq \Cal Y\times _RK\to \Bbb P^1_K$ 
between generic fibres (resp. in the morphism $f_k:Y_k\to \Bbb P^1_k$). It is well-known 
[and easy to verify] that we have the inequality 
$$\delta_{\eta}\ge \delta_{s}.$$
Furthermore, the equality 
$$\delta _{\eta}=\delta_{s}$$ 
holds if and only if $\Cal Y$ is smooth over $R$ [which is equivalent to $\Cal Y_k$ being 
isomorphic to  $Y_k$], as follows from the local criterion for good reduction 
(cf. [Gr-Ma], 3.4). 

We will consider the following assumption.

\definition {3.3.1 Assumption (A)}  Let $n\ge 1$ be a positive integer, and 
$f_k:Y_k\to \Bbb P^1_k$
a cyclic Galois cover with Galois group $G\isom \Bbb Z/p^n\Bbb Z$, with $Y_k$ a smooth $k$-curve.
Let
$g_k:X_k\to \Bbb P^1_k$
be the unique Galois sub-cover of $f_k$ of degree $p^{n-1}$. Assume that $g_k$ 
has a smooth Galois lifting 
$g:\Cal X\to \Bbb P^1_{R'}$ 
[over some finite extension $R'/R$] (cf. Definition 2.5.2).

We say that the Galois cover 
$f_k:Y_k\to \Bbb P^1_k$
satisfies the assumption {\bf (A)}, with respect to the smooth lifting $g$ of the sub-cover $g_k$,
if for all possible Garuti liftings 
$\tilde f:\Cal Y\to \Bbb P^1_{R''}$ of the Galois cover 
$f_k:Y_k\to \Bbb P^1_k$ which dominate $g$ [see preceding discussion], and are defined over a finite 
extension $R''/R'$ [the existence of such an $\tilde f$ is guaranteed by Theorem 2.5.3],
the strict inequality 
$$\delta _{\eta}\defeq \delta_{\tilde f_{K''}} >\delta _s\defeq \delta_{f_k}$$ 
[where $K''\defeq \Fr (R'')$] holds. 

In other words the assumption {\bf (A)} is satisfied 
if there doesn't exist a smooth lifting of $f_k$ which dominates the given smooth lifting $g$ of the 
sub-cover $g_k$ of $f_k$.

Note that if the above revisited version of Oort's conjecture {\bf [Conj-O2-Rev]} (cf Remark 3.1.2) is true
then no Galois cover $f_k:Y_k\to \Bbb P^1_k$ as above satisfies the assumption (A).
\enddefinition

Next, we introduce the notion of fake liftings of cyclic Galois covers between curves, which naturally arise
if cyclic Galois covers satisfy the above assumption (A).

\definition {Definition 3.3.2 (Fake liftings of Cyclic Covers between Curves)} 
Assume that the Galois cover $f_k:Y_k\to \Bbb P^1_k$ satisfies the assumption {\bf (A)},
with respect to the smooth lifting $g$ of the sub-cover $g_k$ (cf. 3.3.1). 
Let
$$\delta \defeq  \min \{\delta_{\tilde f_{K''}}\},$$
where the minimum is taken among all possible Garuti liftings $\tilde f:\Cal Y\to 
\Bbb P^1_{R''}$ of $f_k$ as above, which dominate the smooth lifting $g:\Cal X
\to \Bbb P^1_{R'}$ of the sub-cover $g_k:X_k\to \Bbb P^1_k$. 
[Note that $\delta > \delta _s$ by assumption].

We call a lifting $\tilde f:\Cal Y\to \Bbb P^1_{R''}$ as above satisfying the equality
$$\delta_{\tilde f_{K''}}=\delta$$ 
a fake lifting of the Galois cover $f_k:Y_k\to \Bbb P^1_k$, relative to
the smooth lifting $g$ of the sub-cover $g_k$. Note that if $\tilde f:\Cal Y\to \Bbb P^1_{R''}$ 
is a fake lifting of the Galois cover $f_k$ then $\Cal Y$ is [by definition] not smooth over $R''$.
\enddefinition

\definition {Remark 3.3.3} Fake liftings as in Definition 3.3.2 won't exist if the 
revisited Oort conjecture {\bf [Conj-O2-Rev]} is true, hence the reason we call them fake.
Moreover,  in order to prove the [revisited] Oort conjecture it suffices
to prove that fake liftings do not exist, as follows from the various definitions above.
\enddefinition

\definition {3.4} In this sub-section we introduce some terminology related to the semi-stable geometry of curves,
which will be used in the next sub-section 3.5, where
we investigate the geometry of the [minimal] semi-stable models
of fake liftings of cyclic Galois [of $p$-power order] covers between smooth curves. 
\enddefinition
Let 
$f_k:Y_k\to \Bbb P^1_k$ 
be a finite ramified Galois cover with Galois group $G\isom \Bbb Z/p^n\Bbb Z$ 
a cyclic group of order $p^n$, $n\ge 1$. Let
$G\twoheadrightarrow H\isom \Bbb Z/p^{n-1}\Bbb Z$ 
be the [unique] quotient of $G$ with cardinality $p^{n-1}$. Let 
$g_k:X_k\to \Bbb P^1_k$ 
be the cyclic sub-cover of $f_k$ with Galois group $H$. Assume that there exists
$g:\Cal X\to \Bbb P^1_R$ 
a smooth Galois lifting of $g_k$ over $R$ (cf. Definition 2.5.2). Let 
$\Tilde f:\Cal Y\to \Bbb P^1_R$ 
be a fake lifting of the Galois cover $f_k:Y_k\to \Bbb P^1_k$  [with respect to the smooth lifting $g$ of $g_k$],
which dominates the smooth lifting $g$ of $g_k$ (cf. Definition 3.3.2). [We assume that both $\Tilde f$ and $g$ are defined
over $R$ for simplicity]. We have a natural factorization 
$\Tilde f:\Cal Y@>h>> \Cal X @>g>> \Bbb P^1_R$
where $h:\Cal Y\to \Cal X$ is a finite Galois cover of degree $p$, with 
$\Cal Y$ normal and non smooth over $R$.

Next, we assume that $\Cal Y$ admits a semi-stable model over $R$. 
[It follows from the semi-stable reduction theorem for curves (cf. [De-Mu], and [Ab1])
that $\Cal Y$ admits a semi-stable model after eventually a finite extension
of $R$]. More precisely, we assume that there exists a birational morphism
$$\sigma: \Cal Y'\to \Cal Y$$ 
with $\Cal Y'$ semi-stable, i.e. the special fiber $\Cal Y'_k\defeq \Cal Y'
\times _Rk$ of $\Cal Y'$ is reduced, and its only singularities are ordinary 
double points. We also assume that the ramified points in the morphism $\Tilde f_K:\Cal Y_K\to 
\Bbb P^1_K$ specialise in smooth distinct points of $\Cal Y'_k$. Moreover, we will 
assume that the birational morphism $\sigma$ is minimal with respect to 
the above properties. In particular, the action of the Galois group $G\isom \Bbb Z/p^n\Bbb Z$ 
on $\Cal Y$ extends to an action of $G$ on  $\Cal Y'$. Let 
$$\Cal P\defeq \Cal Y'/G$$ 
be the quotient of $\Cal Y'$ by $G$, and 
$$\tilde f':\Cal Y'\to \Cal P$$
the natural 
morphism [which is Galois with Galois group $G$]. 
Let 
$$\tilde g: \Cal X'\to \Cal P$$
be the unique sub-cover of 
$\tilde f'$ which is Galois with Galois group $H$ [$\Cal X'$
is the quotient of $\Cal Y'$ by the unique subgroup of $G$ with cardinality $p$]. 
Then $\Cal P$ and $\Cal X'$ are semi-stable $R$-curves (cf. [Ra], appendice), 
and we have the following commutative diagram:

$$
\CD
\Cal Y @>h>>  \Cal X @>g>> \Bbb P^1_R \\
@A{\sigma}AA     @AAA  @AAA  \\
\Cal Y'   @>{\tilde h}>>  \Cal X' @>{\tilde g}>> \Cal P \\
\endCD 
$$ 
where the vertical maps are birational morphisms, and the horizontal maps are
finite morphisms. 

To the special fiber $\Cal Y'_k\defeq \Cal Y'\times_Rk$ of $\Cal Y'$
[which is a semi-stable $k$-curve] one associates a graph $\Gamma$ 
whose vertices 
$$\Ver (\Gamma)\defeq \{Y_i\}_{i=0}^{n'}$$ 
are the irreducible components of $\Cal Y'_k$, and edges are the double points 
$$\Edg (\Gamma)\defeq \{y_j\}_{j\in J}$$ 
of $\Cal Y'_k$. A double 
point $y_j\in Y_t\cap Y_s$ defines and edge linking the vertices $Y_t$ and 
$Y_s$. We assume that $Y_0$ is the strict transform of $\Cal Y_k$ [which is irreducible] 
in $\Cal Y'$. 

One also associates to the special fibre 
$\Cal X'_k\defeq \Cal X'\times _R k$ of $\Cal X'$ [which is a semi-stable $k$-curve] 
a graph $\Gamma'$  whose vertices 
$$\Ver (\Gamma')\defeq \{X_i\}_{i=0}^{m}$$ are the irreducible components of 
$\Cal X'_k$, and edges are the double points 
$$\Edg (\Gamma')\defeq \{x_j\}_{j\in J'}$$ 
of $\Cal X_k'$. We assume that $X_0$ is the strict transform of $\Cal X_k\isom X_k$ 
in $\Cal X'$. Then it follows easily [from the fact that $\Cal X$ is 
smooth] that the graph $\Gamma'$ is a tree, and all the irreducible components of 
$\Cal X_k'$ which are distinct from $X_0$ are isomorphic to $\Bbb P^1_k$. 
We choose an orientation of $\Gamma'$ starting from $X_0$ towards the end 
vertices of the tree $\Gamma'$. We have a natural morphism of 
graphs 
$$\Gamma \to \Gamma'.$$ 

Similarly one associates to the special fibre 
$\Cal P_k\defeq \Cal P\times _R k$ of $\Cal P$ [which is a semi-stable $k$-curve] 
a graph $\Gamma''$  whose vertices 
$$\Ver (\Gamma'')\defeq \{P_i\}_{i=0}^{n}$$ 
are the irreducible components of $\Cal P_k$, and edges are the double points 
$$\Edg (\Gamma'')\defeq \{\tilde x_j\}_{j\in J''}$$ 
of $\Cal P_k$. We assume that $P_0$ is the strict transform of 
$\Bbb P^1_k$ [the special fibre of $\Bbb P^1_R$]
in $\Cal P$. The graph $\Gamma''$ is a tree and all the irreducible 
components of $\Cal P_k$ are isomorphic to $\Bbb P^1_k$. We choose 
an orientation of $\Gamma''$ starting from $P_0$ towards the end vertices 
of the tree $\Gamma''$. We have natural morphisms of graphs 
$$\Gamma \to \Gamma'\to \Gamma''.$$ 

The morphism $\Gamma \to \Gamma''$ (resp. $\Gamma'\to \Gamma''$) is
$G$-equivariant (resp. $H$-equivariant). [The graph $\Gamma$ (resp. $\Gamma'$)
is naturally endowed with an action of the group $G$ (resp. $H$)]. 

Let $Y_i$ be a vertex of the graph $\Gamma$. To $Y_i$ one associates two 
subgroups of the Galois group $G\isom \Bbb Z/p^n\Bbb Z$ of the Galois cover 
$\Tilde f:\Cal Y\to \Bbb P^1_R$: the decomposition subgroup $D_i\subseteq G$, and the 
inertia subgroup $I_i\subseteq D_i$, at the generic point of $Y_i$ in the Galois 
cover $\Tilde f$.  We call the [irreducible
component] vertex $Y_i$ of $\Gamma$ an end vertex [or end component] 
of $\Gamma$ if the graph $\Gamma$ is a tree, and if $Y_i$ is an end vertex
 of this tree. We call 
$Y_i$ a separable vertex of $\Gamma$ if the inertia  subgroup $I_i$ 
which is associated to $Y_i$ is trivial. Finally, we call the irreducible 
component $Y_i$ a ramified vertex if there exists a ramified point in the 
morphism $f_K:\Cal Y_K\to \Bbb P^1_K$ which specialises in the component $Y_i$.

Similarly let $X_i$ be a vertex of the graph $\Gamma'$. To $X_i$ one associates two 
subgroups of the Galois group $H\isom \Bbb Z/p^{n-1}\Bbb Z$ of the Galois cover 
$g:\Cal X\to \Bbb P^1_R$: the decomposition subgroup $\Tilde D_i\subseteq H$, and the 
inertia subgroup $\Tilde I_i\subseteq \Tilde D_i$, at the generic point of $X_i$ in the Galois 
cover $g$.  We call the vertex $X_i$ of $\Gamma'$ an end vertex 
of $\Gamma'$ if $X_i$ is an end vertex of the tree $\Gamma'$. We call $X_i$ an internal 
vertex of $\Gamma'$ if $X_i$ is distinct from $X_0$, and the end vertices of $\Gamma'$.
We call $X_i$ a separable vertex of $\Gamma'$ if 
the inertia  subgroup $\Tilde I_i$ which is associated to $X_i$ is trivial. 
Finally, we call the irreducible component $X_i$ a ramified vertex if there exists a 
ramified point in the morphism $g_K:\Cal X_K\to \Bbb P^1_K$ which specialises in the 
component $X_i$.

Finally, By a geodesic in a finite tree linking two vertices we mean the path, or sub-tree, with smallest 
length which links the two vertices.

\definition {3.5} 
In this sub-section we first establish in the next Proposition 
some properties of the [not necessarily minimal] semi-stable model $\Cal X'\to \Cal X$ of 
the smooth lifting $g:\Cal X\to \Bbb P^1_R$ of the Galois sub-cover 
$g_k:X_k\to \Bbb P^1_k$ of $f_k:Y_k\to \Bbb P^1_k$.
\enddefinition

\proclaim {Proposition 3.5.1} let $g_k:X_k\to \Bbb P^1_k$ be a finite ramified Galois cover 
with Galois group $H\isom \Bbb Z/p^{n-1}\Bbb Z$ [$n> 1$], and $X_k$ a smooth $k$-curve. Let 
$g:\Cal X\to \Bbb P^1_R$ be a smooth Galois lifting of $g_k$ over $R$ (cf. Definition 2.5.2). 
Assume that there exists a birational morphism $\Cal X'\to \Cal X$ such that
$\Cal X'$ is semi-stable, the action of $H$ on $\Cal X$ extends to an action on 
$\Cal X'$, and the ramified points in the 
Galois cover $g_K:\Cal X_K\to \Bbb P^1_K$ specialise in smooth distinct points of 
$\Cal X'_k$. [We do not assume that $\Cal X'$ is minimal with respect to the above properties]. 
Let $\Cal P \defeq \Cal X'/H$ be the quotient of $\Cal X'$ by $H$. We have a commutative digram: 
$$
\CD
\Cal X @>g>> \Bbb P^1_R \\
@AAA  @AAA  \\
\Cal X' @>{\tilde g}>> \Cal P \\
\endCD 
$$ 
where $\Cal P$ is a semi-stable $R$-curve, and the vertical maps are birational morphisms.

Let $\Gamma'$ (resp. $\Gamma ''$) be the graph associated to the semi-stable
$k$-curve  $\Cal X'_k$ (resp. $\Cal P_k$).  Let $\Ver (\Gamma')\defeq \{X_i\}_{i=0}^{m}$
(resp. $\Ver (\Gamma'')\defeq \{P_i\}_{i=0}^{n'}$) be the set of vertices of $\Gamma'$
(resp. of $\Gamma''$). Then we have a natural morphism $\Gamma'
\to \Gamma''$ of graphs and  the followings hold.

(i)\ The graphs $\Gamma'$ and $\Gamma''$ are trees. Furthermore, each vertex
$X_i$ (resp. $P_i$) of $\Gamma'$ (resp. of $\Gamma''$) which is distinct from the strict transform of 
$\Cal X_k$ (resp. distinct from the strict transform of the special fibre of $\Bbb P^1_R$) 
is isomorphic to $\Bbb P^1_k$.

Let $X_0$ be the strict transform of $\Cal X_k\isom X_k$ in $\Cal X'$. We choose an 
orientation of the tree $\Gamma'$ starting form $X_0$ towards the end vertices
 of $\Gamma'$. For a vertex
$X_i$ of $\Gamma'$ we will denote by $\Tilde D_i$ (resp. $\Tilde I_i\subseteq \Tilde D_i$) 
the decomposition (resp. inertia) subgroup of $H$ at the generic point of $X_i$. Then:

(ii)\ $\Tilde D_0=H$ and $\Tilde I_0=\{1\}$.

(iii)\ Let $X_i$ be an internal vertex of $\Gamma'$ 
[i.e. $X_i$ is distinct from $X_0$ and from the end vertices of $\Gamma'$], and 
$X_j$ an adjacent vertex to $X_i$ in the direction moving towards the end vertices of $\Gamma'$.
Then the following two cases occur:

(1)\  Either $\Tilde D_i=\Tilde I_i$. In this case $\Tilde D_j=\Tilde D_i$.

(2)\ Or $\Tilde I_i\subsetneq \Tilde D_i$. In this case $\Tilde D_j=\Tilde I_i$
and we have an exact sequence
$$1\to \Tilde D_j\to \Tilde D_i\to \Bbb Z/p\Bbb Z \to 0.$$

Furthermore, in the case (2) if $\Tilde X_i$ denotes the image of $X_i$ in the quotient
$\Cal X'/ \Tilde I_i$ of $\Cal X'$ by $\Tilde I_i$ then the natural morphism $\Tilde X_i\to P_i$,
where $P_i\isom \Bbb P^1_k$ is the image of $X_i$ in $\Gamma''$, is a Galois cover of degree $p$
ramified above a unique point $\infty \in P_i$ [which is the edge of the geodesic linking 
$P_i$ to $P_0$, which is linked to $P_i$]  with Hasse conductor $m=1$ at $\infty$. 

In particular, when we move in the graph  $\Gamma'$ starting from $X_0$ towards the end 
vertices of $\Gamma'$ then the cardinality of the decomposition group $\Tilde D_i$ (resp. 
the cardinality of the inertia subgroup $\tilde I_i$) 
of a vertex $X_i$ decreases. 
More precisely, if when moving from a vertex $X_i$ towards the end vertices  of $\Gamma'$ we encounter a vertex 
$X_j$ then $\Tilde D_j \subseteq \Tilde D_i$ and $\Tilde I_j\subseteq \Tilde I_i$.

(iv)\ Let $X_i$ be a separable vertex of $\Gamma'$ [i.e. $\Tilde I_i=\{1\}$]
which is distinct form $X_0$.
Then either $X_i$ is an internal vertex [of $\Gamma'$] which is adjacent to an end vertex of the graph $\Gamma'$. 
Furthermore, $\Tilde D_i=\Bbb Z/p\Bbb Z$ in this case
and $X_i$ is a Galois cover of $\Bbb P^1_k$ ramified above a unique point $\infty \in \Bbb P^1_k$
with Hasse conductor $m=1$ at $\infty$. [In this case if $X_j$ is the end vertex of $\Gamma'$ which is adjacent
to $X_i$ then $\Tilde D_j=\{1\}$ (cf. (ii), (2))].
Or, $X_i$ is an end vertex of $\Gamma'$, and two cases can occur: either
$\Tilde D_i=\Bbb Z/p\Bbb Z$ and $X_i$ is a Galois cover of $\Bbb P^1_k$ ramified above a unique point $\infty \in \Bbb P^1_k$
[which is the point linking $X_i$ to the rest of the tree $\Gamma'$]
with Hasse conductor $m=1$ at $\infty$, or $\Tilde D_i=\{1\}$ and $X_i$ is adjacent to a [unique] internal separable vertex $X_j$ 
with $\Tilde D_j\isom \Bbb Z/p\Bbb Z$, $\Tilde I_j=\{1\}$,  and $X_j$ is a Galois cover of $\Bbb P^1_k$ ramified above a unique point $\infty \in \Bbb P^1_k$
 [which is the edge of the geodesic linking $P_j$ to $P_0$, which is linked to $P_j$] with Hasse conductor $m=1$ at $\infty$.  

Let $0<j \le n-1$ be an integer. Let $x\in \Cal X_K$ be a ramified point in the morphism
$g_K:\Cal X_K\to \Bbb P^1_K$. We say that the ramified point $x$ is of type $j$ if the
inertia subgroup $\Tilde I_x\subseteq H$ at $x$ is isomorphic to $\Bbb Z/p^j\Bbb Z$.
A vertex $X_i$ of $\Gamma'$ is called a ramified vertex 
of type $j$ if there exists a ramified point $x$ of type $j$ in the morphism $g_K:\Cal X_K\to \Bbb P^1_K$ 
which specialises in the component $X_i$.

(v) Let $X_i$ be a ramified component of $\Gamma'$. Then $X_i$ is of type $j$ for a unique integer
$0<j \le n-1$. In other words if $0<j<j' \le n-1$ are integers then ramified points $x\in \Cal X_K$  
(resp. $x'\in \Cal X_K$) of type $j$ (resp. type $j'$) in the morphism $g_K:\Cal X_K\to \Bbb P^1_K$ 
specialise in distinct irreducible components of $\Cal X_k$. More precisely, if $X_i$ is a ramified 
vertex of type $j$ then the inertia subgroup $\Tilde I_i$ which is associated to $X_i$
has cardinality $p^j$, i.e. $\Tilde I_i\isom \Bbb Z/p^j\Bbb Z$. [In other words the type $j$ of a ramified component $X_i$
is uniquely determined by $X_i$].  

Furthermore, let $P_i$ be the image of $X_i$ in $\Cal P$. Then the natural
morphism $X_i\to P_i$ has the structure of a $\mu_{p^j}$-torsor outside the double points supported by $P_i$,
and the specialisation of the branched points in $P_i$ [in this case $\Tilde D_i=\Tilde I_i$].

(vi) Let $X_i$ be a ramified vertex of $\Cal X_k$ of type $j$.  Then when moving in the graph $\Gamma'$ from $X_i$ 
towards the end vertices of $\Gamma'$ we encounter at most a unique ramified vertex $X_{i'}\neq X_i$. Moreover, in such a 
component $X_{i'}$ specialises a unique ramified point in the morphism $f_K:\Cal X_K\to \Bbb P^1_K$,
and the component $X_{i'}$ is necessarily of the same type $j$ as $X_i$. [In other words the graph $\Gamma'$ separates the directions 
of the ramified vertices of $\Gamma'$ which are of distinct types].

(vii) Assume that $\Cal X$ is minimal [with respect to its defining properties above]. Then the ramified vertices in
the graph $\Gamma'$ are the end vertices of the tree $\Gamma'$.
\endproclaim

\demo{Proof} Assertion (i) is clear and follows immediately from the fact that 
$\Cal X$ is smooth. 

Assertion (ii)
is also clear since $\Cal X_k$ is irreducible and the natural morphism
$\Cal X_k\to \Bbb P^1_k$ [which is isomorphic to $g_k :X_k \to \Bbb P^1_k$] 
is generically Galois with Galois group $H$. 

Next, we prove (iii). Let $X_i$ be an internal vertex of $\Gamma'$, and $X_j$ an adjacent 
vertex to $X_i$ in the direction moving towards the end vertices  of $\Gamma'$. Let $P_i$ (resp. $P_j$)
 be the image of $X_i$ (resp. $X_j$) in $\Cal P$. 
 
 Assume first that $\Tilde D_i=\Tilde I_i$, we will show
 that $\Tilde D_j=\Tilde D_i$ in this case. Let $\Cal X_1\defeq \Cal X'/\Tilde D_i$ be the quotient of $\Cal X'$ by $\Tilde D_i$.
 Then $\Cal X_1$ is a semi-stable $R$-curve, and the configuration of the special fibre $(\Cal X_1)_k$ of $\Cal X_1$
 is a tree-like (cf. (i)). 
 The natural morphism $\Cal X_1\to \Cal P$ is by assumption completely split above the irreducible component
 $P_i$ of $\Cal P_k$, hence [a fortiori] is also completely split above $P_j$. This shows that $\Tilde D_j
 \subseteq \Tilde D_i$.  Assume that $\Tilde D_j\subsetneq \Tilde D_i$. Let $x\defeq X_i\cap X_j$ which is a
 double point of $\Cal X'$  and $x'\defeq P_i\cap P_j$ its image in $\Cal P$.   Let $\Cal X''\defeq \Cal X'/\Tilde D_j$
 be the quotient of $\Cal X'$ by $\Tilde D_j$ [$\Cal X''$ is a semi-stable $R$-curve and the configuration
 of the special fibre $\Cal X''_k$ of $\Cal X''$ is a tree-like], and $X''_i$ the image of $X_i$ in $\Cal X''$. The natural morphism
 $\Cal X''\to \Cal P$ is by assumption completely split above $P_j$, hence also completely split above the double point
 $x'$. In particular, the natural morphism $X''_i\to P_i$ is \'etale above $x'$ and is generically Galois with Galois group $\Tilde D_i/\Tilde D_j$.
 This contradicts the fact that $\Tilde D_i=\Tilde I_i$. Hence $\Tilde D_j=\Tilde D_i$ necessarily.

Assume now that $\Tilde I_i\subsetneq  \Tilde D_i$ and write $D_i'\defeq \Tilde D_i/
\Tilde I_i \neq \{1\}$. Let $\Tilde X_i$ be the image of $X_i$ in the quotient
$\Cal X'/\Tilde I_i$ of $\Cal X'$ by $\Tilde I_i$.
 We have a natural morphism $\Tilde X_i\to P_i$ which is 
generically Galois with Galois group $D_i'$. The vertex $P_i\in \Ver \Gamma''$ 
is an internal vertex of the tree $\Gamma ''$ [as is easily seen since $X_i$ is an internal vertex
of $\Gamma'$], hence is linked to more than one double 
point of $\Gamma''$.  More precisely, $P_i$ is linked to a unique double point $x'$
which links $P_i$ to the geodesic joining $P_i$ and the vertex $P_0$ [$P_0$ is the image
of $X_0$ in $\Cal P$], and [at least another] other double points linking $P_i$ to the geodesics
joining $P_i$ and some of the end vertices of the graph $\Gamma''$.

If the natural morphism $\Tilde X_i\to P_i$ is unramified above the double point $x'$ 
then it is easy to see that this would introduce loops in the configuration of
$\Gamma'$ hence the later won't be a tree. Thus, the morphism $\Tilde X_i\to P_i$ 
must [totally] ramify above the double point $x'$. In particular, this morphism is necessarily
unramified above the remaining double points linking $P_i$ to the end vertices of $\Gamma''$.
Indeed, for otherwise the genus of $\Tilde X_i$ [hence that of $X_i$] would be $>0$, since the degree of 
this morphism is a power of $p$, as follows easily from the Riemann-Hurwitz genus formula, 
and this would contradict the second assertion in (i)].  

Also the degree of the morphism $\Tilde X_i\to P_i$ is necessarily equal to $p$,
and this morphism is only ramified above the double point $x'$ with Hasse conductor $m=1$ at $x'$
[for otherwise the genus of $\Tilde X_i$ [hence that of $X_i$] would be $>0$ for similar reasons as above].
This also shows that $\Tilde D_j\subset \Tilde I_i$ 
[indeed, the natural morphism $\Cal X'/\Tilde I_i\to \Cal P$ is easily seen to be completely
split above the component $P_j$ which is the image of $X_j$ in $\Cal P$],
and that we have a natural exact sequence
$$1\to \Tilde I_i\to \Tilde D_i\to \Bbb Z/p\Bbb Z \to 0.$$

Now we show that $\Tilde D_j=\Tilde I_i$. Assume that $\Tilde D_j\subsetneq \Tilde I_i$.
Let $\Tilde {\Cal X}' \defeq \Cal X'/ \Tilde D_j$ (resp. $\Tilde {\Cal X}'' \defeq \Cal X'/ \Tilde I_i$)
be the quotient of $\Cal X'$ by $\Tilde D_j$ (resp. the quotient of $\Cal X'$ by $\Tilde I_i$),
and $\Tilde X_i'$ (resp. $\Tilde X_i''$) the image of $X_i$ in $\Tilde {\Cal X'}$ (resp. 
$\Tilde {\Cal X}''$). By assumption the natural morphism $\Tilde X_i'\to \Tilde X_i''$
[which is of degree $\ge p$] must be  on the one hand  a homeomorphism, and on the 
other hand completely split above the image of the double point $x\defeq X_i\cap X_j$. This is a contradiction. 
Hence we necessarily have the equality $\Tilde I_i=\Tilde D_j$.
This proves the assertions 1 and 2 in (iii). The remaining assertion in (iii) follows easily from this.

The assertion (iv) follows easily from (iii), and the fact that if in a generically
Galois cover $f:C\to \Bbb P^1_k$ with Galois group a cyclic $p$-group  we have $C\isom \Bbb P^1_k$,
then $f$ has necessarily degree $p$ and is ramified above a unique point $\infty \in \Bbb P^1_k$
with Hasse conductor $m=1$ [as follows easily from the Riemann-Hurwitz genus formula, and Artin-Schreier-Witt theory]. 

Next, we prove (v). Let $0<j \le n-1$ be an integer. Let $x\in \Cal X_K$ be a ramified point in the morphism
$g_K:\Cal X_K\to \Bbb P^1_K$ of type $j$ which specialises in the irreducible component $X_i$ of 
$\Cal X'_k$. We will show that $\Tilde I_i=\Tilde I_x$, where $\Tilde I_x\isom \Bbb Z/p^j\Bbb Z$ is the inertia subgroup at $x$.

Let $\Cal X_2\defeq \Cal X'/\Tilde I_x$ be the quotient of $\Cal X'$ by $\Tilde I_x$, and $\Tilde X_i$ the image of $X_i$ in
$\Cal X_2$. The natural morphism $X_i\to \Tilde X_i$ is a radicial morphism, as follows from [Sa], Corollary 4.1.2,
hence $\Tilde I_x\subset \Tilde I_i$.
Assume that $\Tilde I_x\subsetneq \Tilde I_i$. Let $\Cal X'_2\defeq \Cal X'/\Tilde I_i$, and $\Tilde X'_i$ the image of 
$X_i$ in $\Cal X'_2$. The natural morphism $\Tilde X_i\to \Tilde X'_i$ [which has degree bigger than $1$] is by assumption
on the one hand radicial, and on the other hand unramified above the image of the specialisation of the ramified 
point $x$ in $\Tilde X'_i$, which is a contradiction. Hence we necessarily have $\Tilde I_x=\Tilde I_i$.
The last assertion in (v) follows from Lemma 3.5.5 (see end of $\S3$), and the corresponding assertion in the case where 
$G\isom \Bbb Z/p\Bbb Z$ in [Sa], Corollary 4.1.2.

Assertion (vi) follows directly from the next Lemma 3.5.2, by passing to the quotient of $\Cal X'$ by
the unique subgroup $H'$ of $H$ with cardinality $p$.

Next, we prove (vii). Assume that  $\Cal X$ is minimal with respect to its defining properties.
Let $X_i$ be a ramified vertex of the tree $\Gamma'$. We will show that $X_i$ is necessarily an
end vertex of $\Gamma'$. Assume that $X_i$ [which is distinct from $X_0$] is an internal vertex of
$\Gamma'$.  Let $X_{\Tilde i}$ be an end vertex of $\Gamma'$ which we encounter when moving
in $\Gamma'$ from $X_i$ towards the end vertices of $\Gamma'$, and $\gamma$ the geodesic linking
$X_i$ and $X_{\Tilde i}$. All vertices of $\gamma$ are projective lines (cf. (i)). 

In $\gamma$ there exists 
at most a unique vertex $X_j\neq X_i$ which is a ramified vertex (cf. (vi)). All vertices
of $\gamma$ which are not ramified vertices can be contracted in $\Cal X$ 
without destroying the defining properties of $\Gamma'$.
Thus, we deduce that $\gamma$ contains a unique vertex which is distinct from $X_i$, namely $X_j$, and 
the later $X_j=X_{\Tilde i}$ is an end vertex of $\Gamma$. By (vi) the vertex $X_j$ is of the same type as the vertex $X_i$,
and there exists a unique ramified point in the morphism $\Cal X_K\to \Bbb P^1_K$ which specialises 
in [a smooth point of] $X_j$. The vertex $X_j$ can also be contracted in a [smooth] point of $\Cal X'$
which is supported by $X_i$ and in this point will specialise [after contracting $X_j$] a unique ramified point, which doesn't
destroy the defining properties of $\Cal X'$. But this would contradict the minimality of $\Cal X'$.
Thus, $X_i$ is necessarily a terminal vertex to start with. 
\qed
\enddemo

The following lemma is used in the proof of assertion (vi) in Proposition 3.5.1.

\proclaim{Lemma 3.5.2}  Let $f:\Cal X\to \Cal Y$ be a finite Galois cover between smooth $R$-curves
with Galois group $H'\isom \Bbb Z/p\Bbb Z$, such that the morphism $f_K:\Cal X_K\to \Cal Y_K$
between generic fibres is ramified. Assume that there exists a birational morphism $\Cal X'\to
\Cal X$ such that $\Cal X'$ is a semi-stable $R$-curve, the action of the Galois group $H'$ on 
$\Cal X$ extends to an action of $H'$ on $\Cal X'$, and the ramified points in the morphism
$f_K:\Cal X_K\to \Cal Y_K$ specialise in smooth distinct points of $\Cal X'_k$. Then the graph
$\Gamma'$ associated to the special fibre $\Cal X'_k$ of $\Cal X'$ is a tree. Let $X_0$ be the
strict transform of $\Cal X_k$ in $\Cal X'$. Choose an orientation of $\Gamma'$ starting from 
$X_0$ towards  the end vertices of $\Gamma'$. 

Let $X_i$ be a vertex of $\Gamma'$. 
Assume that $X_i$ is a ramified vertex of $\Gamma'$ [i.e. there exists a 
ramified point in the morphism $f_K:\Cal X_K\to \Cal Y_K$ which specialises in $X_i$].
Then when moving in the graph $\Gamma'$ from $X_i$ towards the end vertices we encounter 
at most a unique ramified vertex $X_j\neq X_i$. Moreover, in such a component $X_j$ specialises 
a unique ramified point in the morphism $f_K:\Cal X_K\to \Cal Y_K$.
\endproclaim

\demo{Proof} The fact that the graph $\Gamma'$ is a tree follows immediately from the fact that
$\Cal X$ is smooth over $R$. Let $X_0$ be the strict transform of $\Cal X_k$ in $\Gamma'$. Let $X_i$ be a 
ramified  component of  $\Cal X_k$. Then $X_i\neq X_0$ as follows from [Sa], Corollary 4.1.2.
Thus, $X_i$ is either an internal or an end component of $\Gamma'$. Assume that $X_i$
is an internal component. Let $X_j$ be an irreducible component of $\Gamma'$ 
which is a ramified vertex and that we encounter when moving from $X_i$ towards the end vertices of $\Gamma'$. 
We will show that only a unique ramified point in the morphism $f_K:\Cal X_K\to \Cal Y_K$
specialises in such a component $X_j$, and that such a component is unique.

After eventually contracting all the irreducible components which form the vertices of the 
geodesics of $\Gamma'$ which link $X_i$ to the end vertices of $\Gamma'$
we can assume that $X_i$ is an end vertex of $\Gamma'$. The component
$X_j$ then contracts to a smooth point $x$ of $X_i$ [which is the specialisation
of some ramified points in the morphism $f_K:\Cal X_K\to \Cal Y_K$].
Let $P_i$ be the image of $X_i$ in the quotient $\Cal Y'\defeq \Cal X'/H'$ of $\Cal X'$
by $H'$, and $y$ the image of $x$ in $\Cal Y'$ which is a smooth point. 
The natural morphism $X_i\to P_i$ is a $\mu_p$-torsor (cf. loc. cit). 
Furthermore, the natural morphism
 $\hat \Cal O_{\Cal X,x}\to \hat \Cal O_{\Cal Y,y}$ between the formal completions at the smooth points $x$
 and $y$ has a degeneration on the boundary of the formal completion 
$\hat \Cal O_{\Cal Y,y}$ of type $(\mu_p,0,h)$  (cf. [Sa], Corollary 4.1.2), and there is a unique ramified point
which specialises in $x$ (cf. loc. cit).
\qed
\enddemo

Proposition 3.5.1 has the following local analog, which describes the geometry of a [minimal] semi-stable 
model of an order $p^n$ automorphism of a $p$-adic open disc [over $K$] without inertia at $\pi$
(cf. [Gr-Ma], 1), and which was proven in [Gr-Ma1] in the case of an order $p$-automorphism. [Though we state
our result in terms of Galois covers between formal germs of smooth curves].

\proclaim{Proposition 3.5.3} let $f:\Tilde {\Cal X}\defeq \Spf A\to \Tilde {\Cal Y}\defeq \Spf B$ be a Galois cover between
connected formal germs of smooth $R$-curves (i.e. $A\isom B \isom R[[T]]$)
which is Galois with Galois group $G\isom \Bbb Z/p^n\Bbb Z$, $n\ge 1$,  and such that the natural morphism
$f_k: \Tilde {\Cal X}_k\defeq \Spec {A}/{\pi A} \to \Tilde {\Cal Y}_k\defeq \Spec {B}/{\pi B}$ 
between special fibres is generically separable. Assume that there exists a birational morphism
$\Tilde {\Cal X}'\to \Tilde {\Cal X}$ such that the ramified points 
in the morphism $f_K:\Tilde {\Cal X}_K\defeq \Spec (A\otimes_RK)\to \Tilde {\Cal Y}_K\defeq \Spec (B\otimes _RK)$
specialise in smooth distinct points of $\Tilde {\Cal X}'_k$, and the action of $G$ on $\Tilde {\Cal X}$
extends to an action of $G$ on $\Tilde {\Cal X}'$. [We do not assume that $\Tilde {\Cal X}'$ is minimal with respect to the above properties].
Let $\Tilde {\Cal Y}'\defeq \Tilde {\Cal X}'/G$ be the quotient
of $\Tilde {\Cal X}'$ by $G$. Then $\Tilde {\Cal Y}'$ is semi-stable (cf. [Ra], Appendice). 
We have a commutative digram: 
$$
\CD
\Tilde {\Cal X} @>f>> \Tilde {\Cal Y}\\
@AAA  @AAA  \\
\Tilde {\Cal X}' @>{\tilde f}>> \Tilde {\Cal Y}' \\
\endCD 
$$ 
where the vertical maps are birational morphisms.

Let $\Gamma'$
(resp. $\Gamma''$) be the graph associated to the special fibre $\Tilde {\Cal X}'$ (resp. $\Tilde {\Cal Y}'$).
Let $\Ver (\Gamma')\defeq \{X_i\}_{i=0}^{m}$
(resp. $\Ver (\Gamma'')\defeq \{Y_i\}_{i=0}^{n'}$) be the set of vertices of $\Gamma'$
(resp. of $\Gamma''$). Then we have a natural morphism $\Gamma'
\to \Gamma''$ of graphs and  the followings hold.

(i)\ The graphs $\Gamma'$ and $\Gamma''$ are trees. Furthermore, each vertex
$X_i$ (resp. $Y_i$) of $\Gamma'$ (resp. of $\Gamma''$) which is distinct from the strict transform of 
[the generic point of] $\Tilde {\Cal X}_k$ in $\Tilde {\Cal X}'$ (resp. distinct from the strict transform of 
[the generic point of] $\Tilde {\Cal Y}_k$ in $\Tilde {\Cal Y}'$) is isomorphic to $\Bbb P^1_k$.

Let $X_0$ be the strict transform of [the generic point of] $\Tilde {\Cal X}_k$ in $\Tilde {\Cal X}'$. 
We choose an orientation of the tree $\Gamma'$ starting form $X_0$ towards the end vertices
 of $\Gamma'$. For a vertex $X_i$ of $\Gamma'$ we will denote by $\Tilde D_i$ (resp. $\Tilde I_i\subseteq \Tilde D_i$) 
the decomposition (resp. inertia) subgroup of $H$ at the generic point of $X_i$. Then:

(ii)\ $\Tilde D_0=H$ and $\Tilde I_0=\{1\}$.

(iii)\ Let $X_i$ be an internal vertex of $\Gamma'$ 
[i.e. $X_i$ is distinct from $X_0$ and from the end vertices of $\Gamma'$], and 
$X_j$ an adjacent vertex to $X_i$ in the direction moving towards the end vertices of $\Gamma'$.
Then the following two cases occur:

(1)\  Either $\Tilde D_i=\Tilde I_i$. In this case $\Tilde D_j=\Tilde D_i$.

(2)\ Or $\Tilde I_i\subsetneq \Tilde D_i$. In this case $\Tilde D_j=\Tilde I_i$
and we have a natural exact sequence
$$1\to \Tilde D_j\to \Tilde D_i\to \Bbb Z/p\Bbb Z \to 0.$$
Furthermore, in the case (2) if $\Tilde X_i$ denotes the image of $X_i$ in the quotient
$\Tilde {\Cal X}'/ \Tilde I_i$ of $\Tilde {\Cal X}'$ by $\Tilde I_i$ then the natural morphism $\Tilde X_i\to P_i$,
where $P_i\isom \Bbb P^1_k$ is the image of $X_i$ in $\Gamma''$, is a Galois cover of degree $p$
ramified above a unique point $\infty \in P_i$ [which is the edge of the geodesic linking 
$P_i$ to $P_0$, which is linked to $P_i$] with Hasse conductor $m=1$ at $\infty$. 

In particular, when we move in the graph  $\Gamma'$ starting from $X_0$ towards the end 
vertices of $\Gamma'$ then the cardinality of the decomposition group $\Tilde D_i$ (resp. 
the cardinality of the inertia subgroup $\tilde I_i$) 
of a vertex $X_i$ decreases. 
More precisely, if when moving from a vertex $X_i$ towards the end vertices  of $\Gamma'$ we encounter a vertex 
$X_j$ then $\Tilde D_j \subseteq \Tilde D_i$ and $\Tilde I_j\subseteq \Tilde I_i$.

(iv)\ Let $X_i$ be a separable vertex of $\Gamma'$ [i.e. $\Tilde I_i=\{1\}$]
which is distinct form $X_0$.
Then either $X_i$ is an internal vertex [of $\Gamma'$] which is adjacent to an end vertex of the graph $\Gamma'$. 
Furthermore, $\Tilde D_i=\Bbb Z/p\Bbb Z$ in this case
and $X_i$ is a Galois cover of $\Bbb P^1_k$ ramified above a unique point $\infty \in \Bbb P^1_k$
with Hasse conductor $m=1$ at $\infty$. [In this case if $X_j$ is the end vertex of $\Gamma'$ which is adjacent
to $X_i$ then $\Tilde D_j=\{1\}$ (cf. (ii), (2))].
Or, $X_i$ is an end vertex of $\Gamma'$, and two cases can occur: either
$\Tilde D_i=\Bbb Z/p\Bbb Z$ and $X_i$ is a Galois cover of $\Bbb P^1_k$ ramified above a unique point $\infty \in \Bbb P^1_k$
[which is the point linking $X_i$ to the rest of the tree $\Gamma'$]
with Hasse conductor $m=1$ at $\infty$, or $\Tilde D_i=\{1\}$ and $X_i$ is adjacent to a [unique] internal separable vertex $X_j$ 
with $\Tilde D_j\isom \Bbb Z/p\Bbb Z$, $\Tilde I_j=\{1\}$,  and $X_j$ is a Galois cover of $\Bbb P^1_k$ ramified above a unique point $\infty \in \Bbb P^1_k$
 [which is the edge of the geodesic linking $P_i$ to $P_0$, which is linked to $P_i$] 
 with Hasse conductor $m=1$ at $\infty$.  

Let $0<j \le n$ be an integer. Let $x\in \Tilde {\Cal X}_K$ be a ramified point in the morphism
$f_K:\Tilde {\Cal X}_K\to \Tilde {\Cal Y}_K$. We say that the ramified point $x$ is of type $j$ if the
inertia subgroup $\Tilde I_x\subseteq G$ at $x$ is isomorphic to $\Bbb Z/p^j\Bbb Z$.
A vertex [irreducible component] $X_i$ of $\Gamma'$ is called a ramified vertex 
of type $j$ if there exists a ramified point $x$ of type $j$ in the morphism $f_K:\Tilde {\Cal X}_K\to \Tilde {\Cal Y}_K$
which specialises in the component $X_i$.

(v) Let $X_i$ be a ramified component of $\Gamma'$. Then $X_i$ is of type $j$ for a unique integer
$0<j \le n$. In other words if $0<j<j' \le n$ are integers then ramified points $x\in \Tilde {\Cal X}_K$  
(resp. $x'\in \Tilde {\Cal X}_K$) of type $j$ (resp. type $j'$) in the morphism $g_K:\Cal X_K\to \Bbb P^1_K$ 
specialise in distinct irreducible components of $\Cal X_k$. More precisely, if $X_i$ is a ramified 
vertex of type $j$ then the inertia subgroup $\Tilde I_i$ which is associated to $X_i$
has cardinality $p^j$, i.e. $I_i\isom \Bbb Z/p^j\Bbb Z$. [In other words the type $j$ of a ramified component $X_i$
is uniquely determined by $X_i$].  Furthermore, let $Y_i$ be the image of $X_i$ in $\Gamma''$. Then the natural
morphism $X_i\to Y_i$ has the structure of a $\mu_{p^j}$-torsor outside the specialisation of the branched points
in $Y_i$, and the double points of $\Tilde {\Cal Y}'_k$ which are supported by $Y_i$.

(vi) Let $X_i$ be a ramified vertex of $\Tilde {\Cal X}'_k$ of type $j$.  Then when moving in the graph $\Gamma'$ from $X_i$ 
towards the end vertices of $\Gamma'$ we encounter at most a unique ramified vertex $X_{i'}\neq X_i$. Moreover, in such a 
component $X_{i'}$ specialises a unique ramified point in the morphism $f_K:\Tilde {\Cal X}'_K\to \Tilde {\Cal Y}'_K$,
and the component $X_{i'}$ is necessarily of the same type $j$ as $X_i$. [In other words the graph $\Gamma'$ separates the directions 
of the ramified components of $\Gamma'$ which are of distinct types].

(vii) Assume that $\Tilde {\Cal X}'$ is minimal [with respect to its defining properties above]. Then the ramified vertices in
the graph $\Gamma'$ are the end vertices of the tree $\Gamma'$.
\endproclaim

\demo{Proof}
Similar to the proof of Proposition 3.5.1.
\qed
\enddemo

Our main result in this section is the following, which describes the semi-stable 
reduction of fake liftings of cyclic Galois covers between smooth curves [assuming they exist], and shows that
fake liftings [if they exist] have semi-stable models with some very specific properties which in some 
sense are reminiscent to the properties of semi-stable models of smooth liftings
of cyclic Galois covers between curves (cf. Proposition 3.5.1).

\proclaim{Theorem 3.5.4} Let $f_k:Y_k\to \Bbb P^1_k$ be a finite ramified Galois cover 
with Galois group $G\isom \Bbb Z/p^n\Bbb Z$ a cyclic group of order $p^n$, 
$n\ge 1$, with $Y_k$ a smooth $k$-curve. Let $g_k:X_k\to 
\Bbb P^1_k$ be the [unique] cyclic sub-cover of $f_k$ with Galois group $H \isom  \Bbb Z/p^{n-1}\Bbb Z$
of cardinality $p^{n-1}$. Assume that there exists a smooth Galois lifting $g:\Cal X\to \Bbb P^1_R$ of $g_k$ defined over $R$ (cf. Definition 2.5.2),
and that $f_k$ satisfies the assumption {\bf (A)} [with respect to the smooth lifting $g$ of $g_k$] 
(cf. 3.3.1). 
Let $\Tilde f:\Cal Y\to \Bbb P^1_R$ be a fake lifting [relative to the smooth lifting $g$ of $g_k$] of the Galois cover 
$f_k:Y_k\to \Bbb P^1_k$ which 
dominates the smooth lifting $g$ of $g_k$, and which we suppose defined
over $R$ (cf. Definition 3.3.2). 

Assume that there exists a minimal birational morphism 
$\Cal Y'\to \Cal Y$ with $\Cal Y'_k\defeq \Cal Y'\times _Rk$ semi-stable,
and such that the ramified points in the morphism $f_K:\Cal Y_K\to 
\Bbb P^1_K$ specialise in smooth distinct points of $\Cal Y'_k$.
Let $\Gamma$ be the graph associated to the 
semi-stable $k$-curve $\Cal Y'_k$.  Write $Y_0$ for the [irreducible component]
vertex of $\Gamma$ which is the strict transform of $\Cal Y_k$ [$\Cal Y_k$ is irreducible] 
in $\Cal Y'_k$. For a vertex $Y_i$ of $\Gamma$ we denote by $D_i$ 
(resp. $I_i\subseteq D_i$) the decomposition (resp. inertia) subgroup of $G$ at the generic point of $Y_i$. 
Then the followings hold.

(i)\ The graph $\Gamma$ is a tree.

(ii)\ The vertex $Y_0\in \Ver (\Gamma)$ is a separable vertex [i.e. $I_0=\{1\}$], and 
$D_0=G$.

Let $H'\isom \Bbb Z/p\Bbb Z$ be the unique subgroup of $G$ with cardinality
$p$. Let $\Cal X'\defeq \Cal Y'/H'$, and $\Cal P\defeq \Cal Y'/G$. Then $\Cal X'$ and $\Cal P$
are semi-stable $R$-curves, and we have a commutative diagram where the vertical maps are 
birational morphisms:

$$
\CD
\Cal Y @>h>>  \Cal X @>g>> \Bbb P^1_R \\
@AAA     @AAA  @AAA  \\
\Cal Y'   @>{\tilde h}>>  \Cal X' @>{\tilde g}>> \Cal P \\
\endCD 
$$

Let $\Gamma'$ (resp. $\Gamma''$) be the graph associated to
the semi-stable $k$-curve $\Cal X'_k$ (resp. $\Cal P_k$). Then the graphs $\Gamma'$ and
$\Gamma''$ are trees (cf. Proposition 3.5.1, (i)), and we have natural morphisms of graphs [actually these are morphisms of trees
by Proposition 3.5.1 (i), and (i) above] 
$$\Gamma\to \Gamma'\to \Gamma''.$$

Let $Y_i$ be a vertex of $\Gamma$ which is distinct from $Y_0$. Let $X_i$ (resp. $P_i$)
be the image of $Y_i$ in $\Cal X'$ (resp. $\Cal P$). Let $\Tilde D_i$ (resp. $\tilde I_i\subseteq \Tilde D_i$)
be the decomposition subgroup (resp. inertia subgroup) of the Galois group $H\defeq G/H'$
which is associated to the generic point of the irreducible component $X_i$.

(iii)\  We have a natural exact sequence
$$0\to H'\to D_i\to \Tilde D_i\to 0.$$
Furthermore, either we have an exact sequence
$$0\to H'\to I_i\to \Tilde I_i\to 0.$$
In particular, $H'\subseteq I_i$ in this case. Or 
$$I_i=\Tilde I_i=\{1\},$$ 
and the inertia subgroups $I_i$ and $\Tilde I_i$ are trivial. The later case can occur only if 
$X_i$ is adjacent, or equal,  to an end vertex of $\Gamma'$ (cf. Proposition  3.5.1, (iv)). 
[See (v) below for a more precise statement related to this case].

Let $0<j \le n$ be an integer. Let $y\in \Cal Y_K$ be a ramified point in the morphism
$f_K:\Cal Y_K\to \Bbb P^1_K$. We say that the ramified point $y$ is of type $j$ if the
inertia subgroup $I_y\subseteq G$ at $y$ is isomorphic to $\Bbb Z/p^j\Bbb Z$. A vertex [irreducible 
component] $Y_i$ of $\Gamma$ is called a ramified vertex of type $j$ if there exists a 
ramified point $y$ of type $j$ in the morphism $\Tilde f_K:\Cal Y_K\to \Bbb P^1_K$ 
which specialises in the component $Y_i$.  The followings hold.

(iv) Let $Y_i$ be a ramified vertex of $\Gamma$. Then $Y_i$ is of type $j$ for a unique integer
$0<j \le n$. In other words if $0<j<j' \le n$ are integers then ramified points $y\in \Cal Y_K$  
(resp. $y'\in \Cal Y_K$) of type $j$ (resp. type $j'$) in the morphism $f_K:\Cal Y_K\to \Bbb P^1_K$ 
specialise in distinct irreducible components of $\Cal Y_k$.  Furthermore, $D_i=I_i\isom \Bbb Z/p^j\Bbb Z$ 
in this case, and the natural morphism $Y_i\to P_i$ has the structure of a $\mu_{p^j}$-torsor outside the specialisation
of the branched points in $P_i$ and the double points of $\Cal P_k$ which are supported by $P_i$.

(v)\  The set of separable vertices of $\Gamma$ which are distinct from $Y_0$ is non empty.
Furthermore, let $Y_i$ be a separable vertex of $\Gamma$ [i.e. $I_i=\{1\}$ is trivial]
which is distinct from $Y_0$. Then $Y_i$ is an 
end vertex of $\Gamma$, and either $D_i\isom \Bbb Z/p\Bbb Z$  or 
$D_i\isom \Bbb Z/p^2\Bbb Z$. [In other words the cardinality of $D_i$ is $\le p^2$].
In the second case the natural morphism $Y_i\to P_i$ is Galois with group $D_i\isom \Bbb Z/p^2\Bbb Z$,
$X_i\to P_i$ is its unique Galois sub-cover of degree $p$, and $X_i$ is ramified
above a unique point $\infty$ of $P_i$ with Hasse conductor $1$ at $\infty$.
[In particular, $X_i\isom \Bbb P^1_k$ in this case]. Moreover, the genus of $Y_i$ is $>0$.
Moreover, no separable vertex of $\Gamma$ is a ramified vertex.

(vi) When we move in the tree $\Gamma$ from a given vertex towards the end vertices
of $\Gamma$ we encounter either ramified vertices or separable vertices of $>0$ genus
[the later are necessarily end components by (v) above].
In particular, an end vertex of the graph $\Gamma$ [which is a tree by (i)] is either a ramified vertex or 
a separable vertex of $\Gamma$.

\endproclaim

\demo {Proof}  The assertion in (ii) is clear since the natural 
morphism $Y_0\to \Bbb P^1_k$ is generically Galois with Galois group $G$. 

Next, we will prove the assertion (iii). 
Let $Y_i$ be a vertex of $\Gamma$ which is distinct from $Y_0$. Let $X_i$ (resp. $P_i$)
be the image of $Y_i$ in $\Cal X'$ (resp. $\Cal P$). Let $\Tilde D_i$ (resp. $\tilde I_i$)
be the decomposition subgroup (resp. inertia subgroup) of the Galois group $H\defeq G/H'$
which is associated to the generic point of the irreducible component $X_i$. 

The image of the decomposition group $D_i$ in $G/H$ via the natural
morphism $G\twoheadrightarrow G/H$ coincides with $\Tilde D_i$. Hence we necessarily
either have an exact sequence $0\to H'\to D_i\to \Tilde D_i\to 0$, since the group $G$ is cyclic,
or we have $D_i=\Tilde D_i=\{1\}$  [if $D_i\cap H'=\{1\}$ then $D_i=\{1\}$ is trivial]
in which case the vertex $X_i$ (resp. $Y_i$) is an end vertex of $\Gamma'$ (resp. of $\Gamma$) (cf. Proposition 3.5.1, (iv)).
The later case can not occur for otherwise the irreducible component
$Y_i$ would be a projective line which is an end vertex of $\Gamma$, and
is not a ramified vertex of $\Gamma$ [as is easily seen since $I_i=\Tilde I_i=\{1\}$ (cf. [Sa], Proposition 4.1.1)],
hence can be contracted in the semi-stable model $\Cal Y'$ without destroying the defining properties of $\Cal Y'$,
and this would contradict the minimal character of $\Cal Y'$. 
Also the image of the subgroup $I_i$ in $G/H$ via the natural
morphism $G\twoheadrightarrow G/H$ coincides with $\Tilde I_i$. Hence we either 
have an exact sequence $0\to H'\to I_i\to \Tilde I_i\to 0$, or the inertia groups $I_i=\Tilde I_i=\{1\}$ 
are trivial, since the group $G$ is cyclic. The later case can occur only if 
$X_i$ is adjacent, or equal,  to an end vertex of $\Gamma'$ (cf. Proposition 3.5.1, (iv)). 

Next, we prove the first assertion in (v). Assume that the set of separable vertices of $\Gamma$ which are 
distinct from $Y_0$ is empty. Let $Y_i$ be a vertex of $\Gamma$ which is distinct from $Y_0$, and
$X_i$ its image in $\Gamma'$. The inertia subgroup $I_i\neq \{1\}$ is non trivial by assumption and
we have a natural exact sequence $0\to H'\to I_i\to \Tilde I_i\to 0$ (cf. (iii)). In particular, the natural
morphism $Y_i\to X_i$ is radicial hence a homeomorphism. Thus, $Y_i$ is a projective line.
Moreover, the natural morphism of graphs $\Gamma\to \Gamma'$ is a homeomorphism in this case,
and the graph $\Gamma$ is a tree. In particular, the arithmetic genus of the special fibre $\Cal Y'_k$ 
is equal to the genus of $Y_k$.  Hence the genera of $Y_K$ and $Y_k$ are equal.
This implies that $Y_K$ has good reduction, which contradicts the fact that $\Cal Y$ is a fake lifting of
$f_k$ [more precisely this contradicts the fact that $\Cal Y$ is not smooth over $R$ (cf. Definition 3.3.2)].

Next, we prove the assertion (i).  In the course of proving (i) we will also prove the second assertion
in (v).
Let's move in the graph $\Gamma'$ starting from the origin vertex
$X_0$ towards a given end vertex $X_{\Tilde i}$ [of $\Gamma'$] along the geodesic $\gamma$ of
$\Gamma'$ which links  $X_0$ and $X_{\Tilde i}$.  Let $X_i$ be a vertex of $\gamma$ which is distinct
from both $X_0$ and $X_{\Tilde i}$. Then $X_i$ is an internal vertex of $\Gamma'$, and the pre-image
of $X_i$ in $\Gamma$ via the natural morphism $\Gamma \to \Gamma'$ consists of  a unique vertex $Y_i$
(cf. (iii) above, more precisely the exact sequence $0\to H'\to D_i\to \Tilde D_i\to 0$). Moreover, the natural
morphism $Y_i\to X_i$ is either radicial [this occurs only if $H'\subseteq I_i$], or is a separable morphism 
in which case $I_i=\Tilde I_i=\{1\}$,
and $X_i$ is adjacent to an end vertex of $\Gamma'$  as follows from (iii).  In fact we will show below that
the later case can not occur. Let now $Y_{\Tilde i}$ be the unique  vertex of $\Gamma$ which is in the pre-image 
of the end vertex $X_{\Tilde i}$ of $\Gamma'$.
The following two cases occur. Either the inertia subgroup $I_{\Tilde i}\neq \{1\}$ [of the group $G$] which is associated
to the vertex $Y_{\Tilde i}$ is non trivial, in which case we have an exact sequence $0\to H'\to I_{\Tilde i} \to \Tilde I_{\Tilde i}\to 0$
, or the inertia subgroups  $I_{\Tilde i}=\Tilde I_{\Tilde i}= \{1\}$ are trivial. In the first case the natural morphism 
$Y_{\Tilde i}\to X_{\Tilde i}$ is radicial, hence a homeomorphism.  

In summary two cases occur: either for every vertex $X_i$ of the geodesic $\gamma$
which is distinct from $X_0$ [in particular $X_i$ may be equal to $X_{\Tilde i}$] and its unique pre-image $Y_i$
in $\Gamma$ we have $I_i\neq \{1\}$ [in particular, $H\subseteq I_i$ in this case], or there exists a vertex 
$X_i$ of $\gamma$ which is distinct from $X_0$
and its unique pre-image $Y_i$ in $\Gamma$ such that $I_i=\Tilde I_i=\{1\}$. 

In the first case the natural
morphism $Y_i\to X_i$ is radicial and the natural morphism $\Tilde h^{-1}(\gamma)\to \gamma$, where
$\Tilde h^{-1}(\gamma)$ is the pre-image of $\gamma$ in $\Gamma$, is a homeomorphism. In particular,
$\Tilde h^{-1}(\gamma)$ is a tree in this case. More precisely, in this case $\Tilde h^{-1}(\gamma)$ is a geodesic
which links $Y_0$ to the unique vertex $Y_{\Tilde i}$ in the pre-image of $X_{\Tilde i}$ which is an end vertex of 
$\Gamma$. Moreover, all vertices of $\Tilde h^{-1}(\gamma)$ which are distinct from $X_0$ are projective lines in this case
and the vertex $Y_{\Tilde i}$ is necessarily a ramified vertex. For otherwise the component $Y_{\Tilde i}$ would be a [non ramified]
projective line hence can be contracted in the semi-stable model $\Cal Y'$ without destroying the defining properties of $\Cal Y'$, and this
would contradict the minimal character of $\Cal Y'$.  Now we shall investigate the second case.  

Assume that the second case above occurs.
In order to show that the graph $\Gamma$ is a tree it suffices to show that the pre-image $\Tilde h^{-1}(\gamma)$ of the geodesic $\gamma$ 
is also a tree in this case [for every possible choice of $\gamma$]. More precisely, we will show that the natural map $\Tilde h^{-1}(\gamma)\to \gamma$ 
is a homeomorphism of trees.
Let $X_i$ be the first vertex of $\gamma$ that we encounter when moving from $X_0$ towards $X_{\Tilde i}$, and $Y_i$ the unique pre-image of 
$X_i$ in $\Gamma$, such that the inertia groups $I_i=\Tilde I_i=\{1\}$ are trivial. We will show that $X_i=X_{\Tilde i}$ is necessarily the end vertex
of $\gamma$ and that the natural morphism $Y_i\to X_i$, which is generically Galois [with Galois group $H'$], is only [totally] ramified above
the unique double point $x_i$ of $\Cal X'_k$ which is supported by $X_i$.  This will complete the proof of the assertion that $\Gamma$ is a tree, and
will also prove the second assertion in (v).

Assume the contrary that $X_i\neq X_{\Tilde i}$ is not the end vertex of $\gamma$. Then $X_i$ is an internal vertex of $\Gamma$, which is
linked to  a unique double point $x_i$
which is an edge of the geodesic which links $X_i$ to $X_0$, and is linked to [at least] another double point $x_{i'}$ which is an edge of the geodesic which links $X_i$
to $X_{\Tilde i}$ [there may be more double points linked to $X_i$ which are edges of the possible geodesics linking $X_i$ to other end vertices of $\Gamma'$].
Moreover, $\Tilde D_{i}\isom \Bbb Z/p\Bbb Z$ in this case (cf. Proposition 3.5.1, (iv))
which necessarily implies that $D_i\isom \Bbb Z/p^2\Bbb Z$, and the natural morphism $X_i\to P_i$
[where $P_i$ is the image of $X_i$ in $\Cal P$] is a Galois cover of degree $p$ ramified above a unique point $\infty \in P_i$ 
[which is the image of the double point $x_i$ in $\Cal P$]
with Hasse conductor $m=1$ at $\infty$
(cf.  Proposition 3.5.1 (iii)). [In particular, $X_i\isom \Bbb P^1_k$ is a projective line].

The natural morphism $Y_i\to X_i$ is a generically Galois morphism with Galois group $\Bbb Z/p\Bbb Z$, and is ramified above the double point $x_i$ with
Hasse conductor $m_i$ at this point [if $X_j$ is the vertex of $\gamma$ such that $x_i=X_i\cap X_j$ and $Y_j$ its unique pre-image in $\Gamma$ then $I_j\neq \{1\}$
by assumption]. Above the double point $x_{i'}$ this morphism is either ramified with Hasse conductor $m_{i'}$ or is unramified. In both cases the double point
$x_{i'}$ produces a non trivial contribution to the arithmetic genus of $\Cal Y'_k$. More precisely, in the first case the contribution of $x_{i'}$ to the arithmetic genus is 
$p-1$, and in the second case it is $\frac {(m_{i'}+1)(p-1)}{2}$.  

We will construct, in order to contradict the above assumption, a new Garuti lifting $f_1:\Cal Y_1\to \Bbb P^1_R$ of the Galois cover 
$f_k:Y_k\to \Bbb P^1_k$ which dominates the smooth lifting $g:\Cal X\to \Bbb P^1_R$ of the Galois subcover $g_k:X_k\to \Bbb P^1_R$ of degree $p^{n-1}$,
and such that the degree of ramification $\delta _1\defeq \delta _{f_{1,K}}$ in the morphism $f_{1,K}:\Cal Y_{1,K}\to \Bbb P^1_K$ between generic fibres
satisfies the inequality $\delta_1 < \delta\defeq \delta _{\Tilde f_K}$. This would contradict the
minimality of $\delta$, i.e. contradicts the fact that $\Tilde f$ is a fake lifting of $f_k$. To simplify the arguments below we will assume that $G=D_i=\Bbb Z/p^2\Bbb Z$.
[The construction of $f_1$ in the general case is done  in a similar fashion by using induced covers from $D_i$ to $G$ (cf.  the construction of Garuti in [Ga], 3,
for similar arguments)].  

Let $X_{1,k}$ be the semi-stable $k$-curve which is obtained from $\Cal X'_k$ by removing the geodesic of the graph $\Gamma'$
which links $X_i$ to the terminal vertex $X_{\tilde i}$, with the vertex $X_i$ removed. Thus, $X_{1,k}$ is a semi-stable $k$-curve with the same arithmetic genus 
as $\Cal X'_k$ [which is the same as that of $X_k$].  Moreover, the graph associated to the semi-stable $k$-curve $X_{1,k}$ is a tree with origin vertex $X_0$,
and the irreducible component $X_i$ is an end vertex of this tree. 
Let $P_{1,k}$ be the image of $X_{1,k}$ in $\Cal P$ [here we view $X_{1,k}$ as a closed sub-scheme of $\Cal X_k$],
 and $Y_{1,k}$ the pre-image of $X_{1,k}$ in $\Cal Y'_k$.  We have natural finite morphisms $Y_{1,k}\to X_{1,k}\to P_{1,k}$ between semi-stable $k$-curves. 

 One can construct a
new finite morphism $Y'_{1,k}\to X_{1,k}\to P_{1,k}$ which above $P_{1,k}\setminus P_i$ coincides with the finite cover which is induced by the above cover
$Y_{1,k}\to X_{1,k}\to P_{1,k}$, above $P_i$ is a generically separable Galois cover with Galois group $D_i=G$ which is ramified only above the unique double 
point $\infty $ of $P_{1,k}$ linking $P_i$ to the geodesic of $\Gamma''$ which links $P_i$ and $P_0$ [the point $\infty$ is the image of $x_i$ in $\Cal P$] , and which above the formal completion of $P_{1,k}$ at the double point $\infty$ coincides with the cover that is induced by the morphisms $Y_{1,k}\to X_{1,k}\to P_{1,k}$. 
In other words in this new cover we eliminate all the irreducible components of the geodesic $\gamma$ that we encounter when moving from
 $X_i$ in the direction of $X_{\Tilde i}$, and we also eliminate the ramification in the morphism $Y_i\to X_i$ which may arise above points of $X_i$ which are distinct from 
 the double point $x_i$ (cf. discussion above).

The finite morphisms $Y'_{1,k}\to X_{1,k}\to P_{1,k}$ can be lifted [uniquely] to finite morphisms $\Tilde {\Cal Y}_1\to \Cal X_1\to \Cal P_1$,
where $\Tilde {\Cal Y}_1\to \Cal P_1$ is a Galois cover with Galois group $G$ which lifts the finite morphism $Y'_{1,k}\to P_{1,k}$, 
and $\Cal X_1\to \Cal P_1$ is the unique sub-cover with Galois group $H$ which lifts the finite morphism $X_{1,k}\to P_{1,k}$, as follows.

First,  we have a natural Galois lifting of the finite morphism $Y'_{1,k}\setminus Y_1'\to P_{1,k}\setminus P_1$ which is the restriction of 
the finite Galois morphism $\Cal Y'\to \Cal P$ to the formal fibre of  $P_{1,k}\setminus P_1$ in $\Cal P$. The restriction of the finite morphism
$\Cal Y'\to \Cal P$ to the formal fibre at the double point $\infty$ [above] provides a natural lifting of the cover above the formal completion of $P_{1,k}$
at the double point $\infty$ which is induced by $Y_{1,k}\to X_{1,k}\to P_{1,k}$.  Second, the restriction of the finite morphism $Y_{1,k}'\to X_{1,k}\to P_{1,k}$
to the irreducible component $P_1\setminus \{\infty \}$  [which is an \'etale torsor]
 can be lifted to an \'etale torsor of the formal fibre of $P_1\setminus \{\infty \}$ in $\Cal P_1$ with Galois group $G$ by the theorems of liftings of
 \'etale covers (cf. [Gr]).
Theses liftings can be patched using formal patching techniques to construct the required Galois cover $\Tilde {\Cal Y}_1\to \Cal X_1\to \Cal P_1$ (cf. [Ga], 
and Proposition 1.2.2). 

Let's now contract [in a Galois  equivariant  fashion] in $\Tilde {\Cal Y}_1$ (resp. in $\Cal X_1$) all the irreducible components of the special fibre
$\Tilde {\Cal Y}_{1,k}$ (resp. $\Cal X_{1,k}$) which are distinct from  $Y_0$ (resp. distinct from $X_0$). We then obtain a normal $R$-curve $\Cal Y_1$
(resp. obtain the smooth $R$-curve $\Cal X$). We have natural finite Galois morphisms $f_1:\Cal Y_1\to \Cal X@>g>> \Bbb P^1_R$, and the Galois cover
$f_1:\Cal Y_1\to \Bbb P^1_R$ is by construction a Garuti lifting of   the Galois cover $f_k:Y_k\to \Bbb P^1_k$ [the fact that $f_1$ dominates the smooth lifting
$g:\Cal X\to \Bbb P^1_R$ of $g_k$ is easily verified, and follows from the above construction]. Let $\delta_1\defeq \delta_{f_{1,K}}$ be the degree of 
the different in the cover $f_{1,K}:\Cal Y_{1,K}\to \Bbb P^1_K$ between generic fibres. 
Then clearly [by construction] we have $\delta_1<\delta$, since the only point of the irreducible 
component $X_i$ of $\Cal X_{1,k}$ which contributes to the arithmetic genus of $\Cal Y_{1,k}$ is the double point $x_i$ [and this contribution is 
the same contribution as in the original cover $\Cal Y'_k\to \Cal P_k$ by construction] (cf. the discussion above). 
But this contradicts the minimality of $\delta$,
and the fact that $\Tilde f:\Cal Y\to \Bbb P^1_R$ is a fake lifting of the Galois cover $f_k$. 

This shows that the irreducible component $X_i=X_{\Tilde i}$ is necessarily an 
end vertex of the geodesic $\gamma$ [hence also an end vertex of the graph $\Gamma'$].  A similar argument shows that the natural morphism
$Y_i\to X_i$ [which is generically separable] is only ramified above the unique double point $x_i$ of $X_i$. This, in particular, shows that 
$\Tilde h^{-1} (\gamma)$ is a tree, and the natural morphism $h^{-1}(\gamma)\to \gamma$ is a homeomorphism of trees. 
Thus, the graph $\Gamma$ is a tree as claimed. Furthermore, $Y_i$ can not be a ramified component
by [Sa], Corollary 4.1.2, which proves the last assertion in (v). 

The proof of (iv) is similar to the proof of Proposition 3.5.1 (v). 

The proof of the second assertion in (v) follows from Proposition 3.5.1, (iv), and [Sa], Corollary 4.1.2.

Finally, we prove (vi). Let $Y_i$ be an internal vertex of $\Gamma$, and $Y_{\Tilde i}$ and end vertex of $\Gamma$ which we encounter when moving from $Y_i$
towards the end vertices.
Let $\gamma$ be the geodesic of $\Gamma$ which links $Y_i$ and $Y_{\Tilde i}$. Assume that $Y_{\Tilde i}$ is neither a ramified component nor a separable component.
Then all vertices of $\gamma$ are projective lines as is easily seen, and can be contracted in $\Cal Y'$ without destroying the defining properties of $\Cal Y'$, which would
contradict the minimal character of $\Cal Y'$. Thus, $Y_i$ is a terminal vertex as claimed.
\qed
\enddemo

The following Lemma 3.5.5 is used in the proof of Proposition 3.5.1 (v), and Theorem 3.5.4 (vi).

\proclaim {Lemma 3.5.5} Let $\Cal X\defeq \Spf A$ be a connected smooth $R$-formal affine scheme. Let $f:\Cal Y\to \Cal X$ be a finite Galois cover
between smooth $R$-formal schemes with $\Cal Y$ connected, with Galois group $G\isom \Bbb Z/p^n\Bbb Z$, $n\ge 1$, and such that the natural 
morphism $f_K:\Cal Y_K\to \Cal X_K$ between generic fibres is \'etale.
[Here the generic fibres $\Cal Y_K$ and $\Cal X_K$ denote the rigid analytic spaces associated to $\Cal Y$ and $\Cal X$ respectively (cf. [Ab])].
Let $\eta$ be the generic point of the special fibre  of $\Cal X$  and $\delta$ the degree of the different in the morphism $f$ above $\eta$.
Assume that $\delta = v_K(p)(1+p+p^2+...+p^{n-1})$.  Then the natural morphism $f_k:\Cal Y_k\to \Cal X_k\defeq \Spec {A}/{\pi A}$
between special fibres has the structure of a $\mu_{p^n}$-torsor.
\endproclaim

\demo{Proof} The Galois cover $f$ has a natural factorization 
$$f:\Cal Y=\Cal Y_n@>f_{n-1}>> \Cal Y_{n-1}\to....\to \Cal Y_2@>{f_1}>> \Cal Y_1\defeq \Cal X,$$
where $f_i:\Cal Y_{i+1}\to \Cal Y_i$ is a Galois cover of degree $p$.  Let $\delta_i$ be the degree of the different in 
the morphism $f_i$ above the generic point $\eta _i$ of $\Cal Y_i$. Then $\delta_i\le v_K(p)$ (cf. [Sa], Proposition 2.3).
The assumption on $\delta$ implies that $\delta_i=v_K(p),\ \ \forall i\in \{1,...,n-1\}$. Hence $f_i:\Cal Y_{i+1}\to \Cal Y_i$ is a torsor under the group scheme
$\mu_{p,R}$ (cf. loc. cit.).  In fact this later property is equivalent to $\delta_i=v_K(p)$.
This implies in particular that the Galois cover $f$ is given by an equation $Z^{p^n}=u$ where $u\in A$
is a unit whose image $\bar u$ in ${A}/{\pi A}$ is not a $p$-th power and hence has the structure of a $\mu_{p^n,R}$-torsor. 
\qed
\enddemo

\subhead
\S 4. The Smoothening Process
\endsubhead 
In this section we introduce the process of smoothening of fake liftings of 
cyclic Galois covers between smooth curves. The idea of smoothening
of fake liftings already germs in the proof of Theorem 3.5.4. 
The smoothening process ultimately aims to show that fake liftings 
as introduced in $\S3$ do not exist. This in turn would imply the validity of the [revisited]
Oort conjecture (cf. Remark 3.3.3).

We use the same notations as in $\S2$, and $\S3$. Especially the Notations 2.1.

\subhead {4.1}
\endsubhead
Let $n\ge 1$ be a positive integer. Let 
$f_k:Y_k\to \Bbb P^1_k$ 
be a finite ramified Galois cover with Galois group 
$G\isom \Bbb Z/p^n\Bbb Z,$ 
a cyclic group of order $p^n$, with $Y_k$ smooth over $k$.  Let 
$g_k:X_k\to \Bbb P^1_k$ 
be the [unique] cyclic sub-cover of $f_k$ with Galois group 
$H \isom  \Bbb Z/p^{n-1}\Bbb Z$
of cardinality $p^{n-1}$. Assume that there exists a smooth lifting 
$g:\Cal X\to \Bbb P^1_R$ 
of $g_k$ defined over $R$ (cf. Definition 2.5.2). 

Assume that $f_k$ satisfies the 
assumption {\bf (A)} in 3.3.1 [with respect to the smooth lifting $g$ of $g_k$].
Let 
$\Tilde f:\Cal Y\to \Bbb P^1_R$ 
be a fake lifting of the Galois cover 
$f_k:Y_k\to \Bbb P^1_k$ 
[with respect to the smooth lifting $g$ of $g_k$], which 
dominates the smooth lifting $g$ of $g_k$, and which we suppose defined
over $R$ (cf. Definition 3.3.2). We assume that there exists a minimal birational morphism 
$\Cal Y'\to \Cal Y$ 
with $\Cal Y'_k\defeq \Cal Y'\times _Rk$ semi-stable,
and such that the ramified points in the morphism $\Tilde f_K:\Cal Y_K\to 
\Bbb P^1_K$ specialise in smooth distinct points of $\Cal Y'_k$.
Let $\Gamma$ be the graph associated to the 
semi-stable curve $\Cal Y'_k$ which is a tree by Theorem 3.5.4 (i). 

Let $H'\isom \Bbb Z/p\Bbb Z$ be the unique subgroup of $G$ with cardinality
$p$. Let 
$\Cal X'\defeq \Cal Y'/H'$,
and 
$\Cal P\defeq \Cal Y'/G$, 
be the quotient of $\Cal Y'$ by $H'$, and the quotient of $\Cal Y'$ by $G$, respectively.
Then $\Cal X'$ and $\Cal P$
are semi-stable $R$-curves, and we have a natural Galois morphism
$f':\Cal Y'\to \Cal P$ with Galois group $G$.
We have a commutative diagram where the vertical maps are 
birational morphisms:

$$
\CD
\Cal Y @>h>>  \Cal X @>g>> \Bbb P^1_R \\
@AAA     @AAA  @AAA  \\
\Cal Y'   @>{\tilde h}>>  \Cal X' @>{\tilde g}>> \Cal P \\
\endCD 
$$

Let $\Gamma'$ (resp. $\Gamma''$) be the graph associated to
the semi-stable $k$-curve $\Cal X'_k$ (resp. $\Cal P_k$). Then the graphs $\Gamma'$ and
$\Gamma''$ are trees (cf. Proposition 3.5.1, (i)), and we have natural morphisms of trees
$$\Gamma\to \Gamma'\to \Gamma''.$$

Let $Y_0$ be the origin vertex of $\Gamma$ [which is the strict transform of $\Cal Y_k$ in $\Cal Y'$],
and let $P_0$ be its image in $\Gamma''$ which is the origin vertex of $\Gamma''$.

\subhead {4.1.1 The semi-stable curve $\Cal P_i$ associated to an internal vertex $P_i$}
\endsubhead
Let $P_i$ be an internal vertex of $\Gamma''$. Let $P_{i,k}$ be the semi-stable
$k$-curve of arithmetic genus $0$, which is obtained from the semi-stable $k$-curve $\Cal P_k\defeq \Cal P\times _Rk$ 
by removing all the geodesics of $\Gamma''$ which link the vertex $P_i$ to the 
end vertices of $\Gamma''$, excluding the vertex $P_i$.  
The graph associated to the semi-stable
curve $P_{i,k}$ is a tree $\Gamma''_i$ in which the vertex $P_i$ is a terminal vertex. Denote by $\infty$ the unique double
point of $P_{i,k}$ which is supported by $P_i$, and which links $P_i$ to the geodesic of $\Gamma''_i$
joining $P_i$ and $P_0$.

Let $\Cal P_i$ be the semi-stable $R$-model
of $\Bbb P^1_R$ which is obtained from the semi-stable $R$-model $\Cal P$ by contracting all the irreducible components of
$\Cal P_k\setminus P_{i,k}$ [here we view $P_{i,k}$ as a closed sub-scheme of $\Cal P_k$].  Then the  special fibre
$\Cal P_{i,k}\defeq \Cal P_i\times _R k$ of $\Cal P_i$ equals $P_{i,k}$, and we have
natural birational morphisms 
$$\Cal P\to \Cal P_i\to \Bbb P^1_R.$$ 
Let $\Cal P_i'$ be the formal fibre of 
$P_i\setminus \{\infty\}$ in $\Cal P_i$. Then 
$$\Cal P_i'\isom \Spf R<S>$$ 
is a formal closed disc. Let 
$\Cal P''_i$ be the formal fibre of $P_{i,k}\setminus \{P_i\}$ in
$\Cal P_i$, and $\Cal P_{i,\infty}$ the formal fibre of $\Cal P_i$ at $\infty$ which is a
 formal open annulus, i.e. 
 $$\Cal P_{i,\infty}\isom \Spf \frac {R[[S,T]]}{(ST-\pi^e)},$$
 for some integer $e\ge1$ [actually $e$ is necessarily divisible by a suitable power of $p$].

Note that the semi-stable $R$-curve $\Cal P_i$ is obtained by patching $\Cal P_i'$ 
and $\Cal P''_i$ along the open annulus $\Cal P_{i,\infty}$.

Next, we define the important concept of a removable vertex in Definition 4.1.2, 
and the smoothening process in Definition 4.1.3.

\definition {Definition 4.1.2 (Removable Vertex of $\Gamma''$)}  [We use the same notations and assumptions as above]. 
We say that $P_i$ is a removable vertex of the tree $\Gamma''$ if there exists a finite Galois cover
$$f'_1:\Cal Y_1'\to \Cal P_i$$ 
[where $\Cal P_i$ is as in 4.1.1] with Galois group $G$, satisfying the following three conditions.
\par
(i)\ The restriction of the Galois cover 
$f'_1$ to $\Cal P_i''$ (resp. to $\Cal P_{i,\infty}$)
is isomorphic to the restriction of the Galois cover 
$$f':\Cal Y'\to \Cal P$$ 
[which is the semi-stable minimal model of the fake lifting $\Tilde f:\Cal Y\to \Bbb P^1_R$ of $f_k$]
above $\Cal P_i''$ (resp. above $\Cal P_{i,\infty}$).
\par
(ii) Let $g'_1:\Cal X_1'\to \Cal P_i$ be the unique Galois sub-cover of $f_1'$ of degree $p^{n-1}$.
Then $g'_1$ is generically [Galois] isomorphic to the Galois cover $g:\Cal X\to \Bbb P^1_R$ which is the given smooth lifting of $g_k$.
\par
(iii)\ The arithmetic genera $g$ (resp. $g_1$)
of the special fibres $\Cal Y'_k$ (resp. $\Cal Y_{1,k}'\defeq \Cal Y_1'\times _R k$) satisfy the inequality 
$$g_1<g.$$
\enddefinition

\definition{Definition 4.1.3 (Smoothening of a Fake Lifting)} [We use the same notations and assumptions as above].
Assume that $P_i$ is a removable vertex in the sense of Definition 4.1.2. Let 
$f'_1:\Cal Y_1'\to \Cal P_i$ 
be the corresponding Galois cover with Galois group $G$ (which is given by Definition 4.1.2).
Let $\Cal Y_1$ be the normal $R$-curve which is obtained from $\Cal Y_1'$ by contracting all the irreducible components
of $\Cal Y_{1,k}'$ which are distinct from $Y_0$. The Galois cover $f'_1$ induces naturally a Galois cover 
$$f_1:\Cal Y_1\to \Bbb P^1_R$$
with Galois group $G$ [since the above contraction procedure is Galois  equivariant].

The inequality $g_1<g$ implies [in fact is equivalent to the fact] that the degree of the [generic] 
different $\delta_1\defeq \delta_{f_{1,K}}$ 
in the natural morphism 
$$f_{1,K}:\Cal Y_{1,K}\defeq \Cal Y_1\otimes _R K\to \Bbb P^1_K$$ 
between generic fibres satisfies the inequality
$$\delta_1<\delta\defeq \delta _{\Tilde f_K}.$$

We call the Galois cover $f_1:\Cal Y_1\to \Bbb P^1_R$ a smoothening of the fake lifting 
$\Tilde f:\Cal Y\to \Bbb P^1_R$.

Note that [by property (ii) in Definition 4.1.2] the Galois cover $f_1:\Cal Y_1\to \Bbb P^1_R$
is a Garuti lifting of the Galois cover $f_k:Y_k\to \Bbb P^1_k$, which dominates the
smooth lifting $g:\Cal X\to \Bbb P^1_R$ of the Galois sub-cover $g_k:X_k\to \Bbb P^1_k$.
[This last property may be used to define the notion of a smoothening of a fake lifting
independently from Definition 4.1.2]
\enddefinition

\subhead {4.2} 
\endsubhead
The existence of a removable vertex in the tree $\Gamma''$, which implies [by definition] the existence of a 
smoothening  $f_1:\Cal Y_1\to \Bbb P^1_R$ of the fake lifting $\Tilde f:\Cal Y\to \Bbb P^1_R$ [more precisely, the 
above inequality $\delta_1<\delta$] (cf. Definition 4.1.3), 
contradicts the fact that $\Tilde f$ is a fake lifting [i.e. contradicts the minimality
of the generic different $\delta$ of $\Tilde f$], hence will prove the [revisited] Oort 
conjecture for the Galois cover $f_k:Y_k\to \Bbb P^1_k$ [and the smooth lifting $g$ of $g_k$] (cf. Remark 3.3.3). 
More precisely, we have the following.

\proclaim {Proposition  4.2.1} Let  $f_k:Y_k\to \Bbb P^1_k$ 
be a finite ramified Galois cover with Galois group 
$G\isom \Bbb Z/p^n\Bbb Z$, and $Y_k$ is a smooth $k$-curve. Let 
$g_k:X_k\to \Bbb P^1_k$ be the Galois sub-cover of $f_k$ with Galois group 
$H \isom  \Bbb Z/p^{n-1}\Bbb Z$. Assume that there exists
a smooth Galois lifting $g:\Cal X\to \Bbb P^1_R$ of $g_k$ defined over $R$ (cf. Definition 2.5.2). 
Assume that $f_k$ satisfies the 
assumption {\bf (A)} in 3.3.1 [with respect to the smooth lifting $g$ of $g_k$].
Let $\Tilde f:\Cal Y\to \Bbb P^1_R$ be a fake lifting of the Galois cover 
$f_k:Y_k\to \Bbb P^1_k$ [with respect to the smooth lifting $g$ of $g_k$],
which dominates the smooth lifting $g$ of $g_k$, which we suppose defined
over $R$ (cf. Definition 3.3.2). We assume that there exists a minimal birational morphism 
$\Cal Y'\to \Cal Y$  with $\Cal Y'_k\defeq \Cal Y'\times _Rk$ semi-stable,
and such that the ramified points in the morphism $f_K:\Cal Y_K\to 
\Bbb P^1_K$ specialise in smooth distinct points of $\Cal Y'_k$.
Let $\Cal P\defeq \Cal Y'/G$ be the quotient of $\Cal Y'$ by $G$, and $\Gamma''$ the tree which is associated to the 
special fibre $\Cal P_k$ of $\Cal P$. 

Under these assumptions suppose that there exists an internal vertex $P_i$ of the tree $\Gamma''$ 
which is a removable vertex of $\Gamma''$ in the sense of Definition
4.1.2, or equivalently that there exists a smoothening $f_1:\Cal Y_1\to \Bbb P^1_R$ of
the fake lifting $\Tilde f:\Cal Y\to \Bbb P^1_R$ in the sense of Definition 4.1.3. 
Then the [revisited] Oort conjecture {\bf [Conj-O2-Rev]} is true for the Galois cover $f_k:Y_k\to \Bbb P^1_k$, and the 
smooth lifting $g$ of the Galois sub-cover $g_k$.
\endproclaim

One can show that fake liftings of cyclic Galois covers between smooth curves, assuming they exist, always
admit a smoothening in the case of cyclic Galois covers of degree $p$. This provides
an alternative proof of the Oort conjecture in the case of a cyclic Galois group $G\isom \Bbb Z/p
\Bbb Z$ of order $p$. This proof doesn't use the equation describing the degeneration of the Kummer
equation of degree $p$ to the Artin-Schreier equation (as in [Se-Oo-Su], and [Gr-Ma]), 
but rather uses the degeneration of the 
Kummer equation to a radicial equation (see proof of Proposition 4.2.2). 
More precisely, we have the following.
 
\proclaim {Proposition  4.2.2} Assume that $R$ contains a primitive $p$-th root of unity. 
Let  $f_k:Y_k\to \Bbb P^1_k$ 
be a finite ramified Galois cover with Galois group 
$G\isom \Bbb Z/p\Bbb Z$, and $Y_k$ is a smooth $k$-curve. Assume that $f_k$ satisfies the 
assumption {\bf (A)} in 3.3.1. 
The assumption {\bf (A)} in this case means that $f_k$ admits 
no smooth lifting, and a fake lifting means a Garuti lifting with minimal generic different.
Let $\Tilde f:\Cal Y\to \Bbb P^1_R$ be a fake lifting of the Galois cover 
$f_k:Y_k\to \Bbb P^1_k$ [which we suppose defined
over $R$] (cf. Definition 3.3.2). We assume that there exists a minimal birational morphism 
$\Cal Y'\to \Cal Y$  with $\Cal Y'_k\defeq \Cal Y'\times _Rk$ semi-stable,
and such that the ramified points in the morphism $f_K:\Cal Y_K\to 
\Bbb P^1_K$ specialise in smooth distinct points of $\Cal Y'_k$.
Let $\Cal P\defeq \Cal Y'/G$ be the quotient of $\Cal Y'$ by $G$ [$\Cal P$ is a semi-stable $R$-model of $\Bbb P^1_R$] ,
and $\Gamma'$ the tree which is associated to the 
special fibre $\Cal P_k$ of $\Cal P$.  

Then there exists an internal vertex $P_i$ of the 
tree $\Gamma'$ which is a removable vertex of $\Gamma'$ in the sense of Definition
4.1.2. In particular, the [revisited] Oort conjecture is true for the Galois cover $f_k:Y_k\to \Bbb P^1_k$
(cf. Proposition 4.2.1).
\endproclaim

\demo{Proof} We can assume, without loss of generality, that the morphism $f_k$ is ramified above a unique point $\infty$ of 
$\Bbb P^1_k$, i.e. work within the framework of {\bf [Conj-O3]}.
Let $P_0$ be the origin vertex of the tree $\Gamma'$, and $P_1$ the [unique] vertex of $\Gamma'$ which is adjacent
to $P_0$. We will show that $P_1$ is a removable vertex of $\Gamma'$. 

The semi-stable $R$-curve
$\Cal P_1$ (cf. 4.1.1) in this case has a special fibre $\Cal P_{1,k}$ which consists of the two
irreducible [smooth] components $P_0$ and $P_1$, which meet at the unique double point $\infty$.  

Let $\Cal P_1'\isom \Spf R<\frac {1}{T}>$ be the formal fibre of $P_1\setminus \{\infty\}$ in $\Cal P_1$, $\Cal P_{1,\infty}$
the formal completion of $\Cal P_1$ at $\infty$, and $\Cal P''_1$ the formal fibre of $\Cal P_{1,k}\setminus P_1$ in
$\Cal P_1$. The natural Galois morphism $\Cal Y'\to \Cal P$ restricts to Galois morphisms $\Cal Y''_1\to \Cal P''_1$,
and $\Cal Y'_{y}\to \Cal P_{1,\infty}$, where $\Cal Y'_{y}$ is the formal completion of $\Cal Y'$ at the unique
double point $y$ above $\infty$. 

The degeneration type of the Galois cover $\Cal Y'_{y}\to \Cal P_{1,\infty}$ 
on the boundary which is linked to $\Cal P'_1$ is necessarily radicial of type
$(\alpha_p,-m,0)$ where $m>0$ is an integer prime to $p$
[since $P_1$ is an internal vertex of $\Gamma'$], or of type $(\mu_p,-m,0)$ where $m$ is as above.
We only treat the first case, the second case is treated in a similar way (cf. [Sa], Proposition 3.3.1, (a2)).

In the first case the Galois cover $\Cal Y'_{y}\to \Cal P_{1,\infty}$ induces a Galois cover
on the boundary which is linked to $\Cal P'_1$ given by an equation
$X^p=1+\pi ^{tp}T^m$, for a suitable choice of $T$ as above, and $t<v_K(\lambda)$ (cf. Proposition 1.3.2). 
Here $\lambda=\zeta_1-1$, and $\zeta_1$ is a primitive $p$-th root of $1$.

Consider the Galois cover $\Cal Y'_1\to \Cal P_1$ which is generically given by the equation
$X^p=T^{-\alpha}(T^{-m}+\pi ^{pt})$ where $\alpha$ is an integer such that $\alpha+m\equiv0\mod p$. 
Then $\Cal Y'_1$ is smooth over $R$, and the natural morphism $\Cal Y'_{1,k}\to \Cal P_{1,k}$ between special fibres is radicial
(cf. [Sa], Proposition 3.3.1, (b)).
The above coverings can be patched using formal patching techniques to construct a Galois cover $\Tilde {\Cal Y}_1\to \Cal P_1$ 
with Galois group $G$ between semi-stable $R$-curves (cf. [Ga], and Proposition 1.2.2), and by construction
the arithmetic genus $g_1$ of the special fibre $\Tilde {\Cal Y}_{1,k}$ [which is in fact equal to that of $Y_k$]
satisfies the inequality $g_1<g$ as required. 
\qed
\enddemo

\subhead {4.3}
\endsubhead
Next, we will give some sufficient conditions for the existence of removable vertices in the 
case where the Galois group $G\isom \Bbb Z/p^2\Bbb Z$ has order $p^2$.

\proclaim{Proposition 4.3.1} Assume that $R$ contains a primitive $p^2$-th root of unity.
Let  $f_k:Y_k\to \Bbb P^1_k$ 
be a finite ramified Galois cover with Galois group 
$G\isom \Bbb Z/p^2\Bbb Z$, and $Y_k$ a smooth $k$-curve. 
Let 
$g_k:X_k\to \Bbb P^1_k$ be the Galois sub-cover of $f_k$ with Galois group 
$H \isom  \Bbb Z/p\Bbb Z$, of cardinality $p$. Assume that there exists a smooth Galois lifting
$g:\Cal X\to \Bbb P^1_R$ of $g_k$ defined over $R$ (cf. Definition 2.5.2). 
Assume that $f_k$ satisfies the 
assumption {\bf (A)} in 3.3.1 [with respect to the smooth lifting $g$ of the Galois sub-cover $g_k$].
Let $\Tilde f:\Cal Y\to \Bbb P^1_R$ be a fake lifting of the Galois cover 
$f_k:Y_k\to \Bbb P^1_k$ which 
dominates the smooth lifting $g$ of $g_k$ [which we suppose defined
over $R$] (cf. Definition 3.3.2). We assume that there exists a minimal birational morphism 
$\Cal Y'\to \Cal Y$  with $\Cal Y'_k\defeq \Cal Y'\times _Rk$ semi-stable,
and such that the ramified points in the morphism $f_K:\Cal Y_K\to 
\Bbb P^1_K$ specialise in smooth distinct points of $\Cal Y'_k$.

Let $\Cal P\defeq \Cal Y'/G$ be the quotient of $\Cal Y'$ by $G$
[which is a semi-stable $R$-model of $\Bbb P^1_R$],
and $\Gamma''$ the tree which is associated to the 
special fibre $\Cal P_k$ of $\Cal P$. Assume that there exists an internal vertex $P_i$
of $\Gamma''$ which satisfies the following properties.

(i)\ The pre-image of $P_i$ in $\Gamma$ contains no ramified vertex.
 
(ii)\  When moving in the tree $\Gamma''$ from the vertex $P_i$ towards the end vertices of 
$\Gamma''$ we encounter a vertex [necessarily terminal by Theorem 3.5.4 (v)] whose pre-image in $\Gamma$ contains 
a separable vertex. 

(iii)\ When moving in the tree $\Gamma''$ from the vertex $P_i$ towards the end vertices of 
$\Gamma''$ we encounter a unique vertex  whose pre-image in $\Gamma$ contains [in fact consists of]
a ramified vertex of type $2$. 

(iv)\ When moving in the tree $\Gamma''$ from the vertex $P_i$ towards the end 
vertices of $\Gamma''$ we encounter no  vertex whose pre-image in $\Gamma$ contains a ramified vertex of type $1$.

Then $P_i$ is a removable vertex of $\Gamma''$ in the sense of Definition 4.1.2, and the revisited Oort conjecture
{\bf [Conj-O2-Rev]} is true for the Galois cover
$f_k:Y_k\to \Bbb P^1_k$, and the smooth lifting $g$ of the Galois sub-cover $g_k$. 
\endproclaim

\demo{Proof}  Let $Y_i$ be a vertex of the graph $\Gamma$ which is in the pre-image of the vertex 
$P_i$, and $D_i$ (resp. $I_i$) the decomposition (resp. inertia) subgroup of the Galois group $G$
at the generic point of $Y_i$.  Then $I_i\neq \{1\}$, since the vertex $Y_i$ is not terminal (cf. Theorem 3.5.4 (v)).
Moreover, $I_i=D_i=G$, for otherwise we will contradict the assumption (iii) satisfied by $P_i$ above.

Let $\Cal P_i$, $\Cal P'_{i}$, $\Cal P''_i$
and $\Cal P_{i,\infty}$ be as in 4.1.1. Let $H'\subset G$ be the unique
subgroup of $G$ with cardinality $p$, and $\Cal X'\defeq \Cal Y'/H'$ the quotient of $\Cal Y'$ by $H'$.
We have natural morphisms $f':\Cal Y'\to \Cal X'\to \Cal P$.

The Galois cover $\Cal X'\to \Cal P$ induces above the irreducible component $P_i$ of $\Cal P_k$, 
outside the specialisation of the branched points, and the double points of $\Cal P$ supported by $P_i$,
an $\Cal H_{pt,R}$-torsor (cf. 1.3.1), 
where $pt<v_K(\zeta_1-1)$, and $\zeta _1$ is a primitive p-th root of 1.
This torsor is generically given by an equation
$$Z^{p}=1+\pi^{tp^2}g(T),$$ 
where $1+\pi ^{tp^2}g(T)\in \Fr (R<\frac {1}{T}>)$ has $m+1$ distinct geometric zeros in $\Cal P'_1$,
which we may assume without loss of generality specialise in the point $\frac {1}{t}=0$ at infinity 
[the later follows from the uniqueness of the ramified vertex of type 2
in the assumption (iii)].  
We will assume for simplicity that $g(T)=T^m$. [The general case is treated in a similar fashion].

The above Galois cover $f':\Cal Y'\to \Cal P$ induces a cyclic Galois cover  
$\Cal  Y'_{\infty} \to  \Cal P_{i,\infty}$ of degree $p^2$
above the formal open annulus
$\Cal P_{i,\infty}$, with  $\Cal  Y'_{\infty}$ connected [since $D_i=I_i=G$], 
which induces a cyclic Galois cover 
$$f'_{\infty,1}:\Cal Y'_{\infty,1}\to \Cal X'_{\infty,1}\to  \Cal P_{i,\infty,1}\isom \Spf R[[T]]\{T^{-1}\}$$ 
of degree $p^2$ above the formal boundary $\Cal P_{i,\infty,1}\isom \Spf R[[T]]\{T^{-1}\}$ of $\Cal P_{i,\infty}$
which is linked to $\Cal P_i'\isom \Spf R<\frac {1}{T}>$.
We will give an explicit description of the Galois cover $f'_{\infty,1}$, using the assumptions satisfied by the vertex $P_i$.

The Galois cover 
$$\Cal X'_{\infty,1}\to  \Cal P_{i,\infty,1}$$ 
is a torsor under the group scheme $\Cal H_{pt,R}$ [where $t$ is as above]
which
has a degeneration type $(\alpha_p,-m,0)$, where $m>1$ is as above
[this results form the assumption (iii) satisfied by $P_i$], and is given by an equation
$$\frac {(\pi ^{pt}X_1+1)^p-1}{\pi ^{p^2 t}}=T^m,\tag *'$$
where $pt<v_K(\zeta_1-1)$, and $\zeta _1$ is a primitive p-th root of 1 as above [in general replace $T^m$ by $g(T)$ above].  
The $\alpha_p$-torsor
$$\Cal X'_{\infty,1,k}\to  \Cal P_{i,\infty,1,k}$$ 
at the level of special fibres is given by the equation 
$$x_1^p=t^m,$$ 
where $x_1= X_1\mod \pi$, and $t= T\mod \pi$.

From the above equation (*') we deduce that in $\Cal X'_{\infty,1}$, we have 
$$T=(X_1^{\frac {1}{m}})^p [1+\sum _{k=1}^{p-1}\binom {p}{k} \pi ^{pt(k-p)}X_1^{k-p}]^{\frac {1}{m}}.$$
In particular, $\Cal X'_{\infty,1,k}\isom \Spf R[[T_i]]\{T_i^{-1}\}$, and $X_1^{\frac {1}{m}}$ is a parameter of $\Cal X'_{\infty,1,k}$.

Moreover, the Galois cover $\Cal Y'_{\infty,1}\to \Cal X'_{\infty,1}$ is given by an equation
$$X_2^p=(1+\pi ^{pt}X_1)(1+\pi ^{ps}f(T))$$
where $f(T)\in \Fr (R<\frac {1}{T}>)$ is such that $(1+\pi ^{ps}f(T))$ is a unit in $\Cal P_i$, for otherwise we will contradict the assumption (iv) satisfied by $\Cal P_i$.
We can assume without loss of generality that $1+\pi^{pt}f(T)\in R<\frac {1}{T}>$.
We will give an explicit description [by equations] of the degeneration of the Galois cover $\Cal Y'_{\infty,1}\to \Cal X'_{\infty,1}$.

Assume for simplicity that $f(T)=T^{-m_1}$,  with $m_1>0$. 
The general case is treated in a similar fashion. 
Thus, our equation is
 $$X_2^p=(1+\pi ^{pt}X_1)(1+\pi ^{ps}T^{-m_1}).$$
Assume first that $t\le s$. Then on the level of special fibres the $\alpha_p$-torsor
$$\Cal Y'_{\infty,1,k}\to \Cal X'_{\infty,1,k}$$ 
is given [in the case where $s=t$ one has to eliminate $p$-powers]  by the equation
$$x_2^p=x_1=(x_1^{\frac {1}{m}})^m$$
where $x_1= X_1\mod \pi$ ($t^{-1}$ becomes a $p$-power in $\Cal X'_{\infty,1,k}$).
In this case the above cover $\Cal Y'_{\infty,1}\to \Cal X'_{\infty,1}$ 
is a torsor under the group scheme $\Cal H_{t,R}$, and has a degeneration of type
$(\alpha_p,-m,0)$. [Note that $x_1^{\frac {1}{m}}$ is a parameter of $\Cal X'_{\infty,1,k}$].

Assume now that $s<t$. Then
$$X_2^p=1+\pi^{ps}T^{-m_1}+\pi^{pt}X_1+\pi^{p(t+s)}X_1T^{-m_1},$$
which is not an integral equation for $\Cal Y_{\infty,1}'$, since $T^{-m_1}$ is a $p$-power mod $\pi$ in $\Cal X'_{\infty,1,k}$.  

To obtain an integral equation we need first to replace $T^{-m_1}$ by its expression, which is deduced from the above description of $T$,
$$T^{-m_1}=(X_1^{\frac {1}{m}})^{-m_1p} [1+\sum _{k=1}^{p-1}\binom {p}{k} \pi ^{pt(k-p)}X_1^{k-p}]^{\frac {-m_1}{m}}.$$
Thus,
$$X_2^p=1+\pi ^{ps}(X_1^{\frac {1}{m}})^{-m_1p}+.....+\pi^{pt}X_1+...,$$
where the remaining terms have coefficients with a valuation which is greater than $ps$.  After replacing $1+\pi ^{ps}(X_1^{\frac {1}{m}})^{-m_1p}$
by 
$$(1+\pi ^s (X_1^{\frac {1}{m}})^{-m_1})^p-.......,$$
and multiplying the above equation by $(1+\pi ^s (X_1^{\frac {1}{m}})^{-m_1})^{-p}$, we reduce to an equation
$$(X'_2)^p=1+\pi^{pt}(X_1^{\frac {1}{m}})^m+.....,$$
where the remaining terms have coefficients with a valuation which is greater than $pt$. 

In particular, the Galois cover
$\Cal Y'_{\infty,1}\to \Cal X'_{\infty,1}$ is a torsor under the group scheme $\Cal H_{t,R}$ and has a degeneration 
of type $(\alpha_p,-m,0)$. More precisely, the $\alpha_p$-torsor  
$\Cal Y'_{\infty,1,k}\to \Cal X'_{\infty,1,x}$ on the level of special fibres is given by an equation
$$\Tilde x^p_2=x_1=(x_1^{\frac {1}{m}})^m.$$

The Galois cover $f':\Cal Y'\to \Cal P$ restricts to Galois covers $\Cal Y_1''\to \Cal P''_i$, and
$\Cal Y_{1,\infty}'\to \Cal P'_{i,\infty}$, above $\Cal P''_i$, and $\Cal P_{i,\infty}$, respectively.
Consider the cyclic Galois cover $\Cal Y_1\to \Cal X_1\to \Cal P_i'$ of degree $p^2$ which is generically given 
by the equations
$$\frac {(\pi ^{pt}X_1+1)^p-1}{\pi ^{p^2 t}}=T^m,$$
[in general replace $T^m$ by $g(T)$ above], and
$$X_2^p=(1+\pi ^{pt}X_1)(1+\pi ^{ps}f(T)),$$
where $t$, s, and $f(T)$ are as above. This Galois
cover on the generic fibre is ramified only at ramified points of type $2$ [$(1+\pi^{ps} f(T))$ is a unit in $\Cal P_i'$].   
Furthermore, both $\Cal X_1$ and $\Cal Y_1$ are smooth, 
and the arithmetic genus of the special fibre $\Cal Y_{1,k}$ is $0$. Indeed, $\Cal X_1$ is smooth,
and the $\alpha_p$-torsor $\Cal Y_{1,k}\to \Cal X_{1,k}$ is given by an equation $\tilde x_2^p=x_1$
by arguments similar to the one above. [One also uses the fact that $x_1^{\frac {1}{m}}$ is a parameter on $\Cal X_{1,k}$].  

The above coverings can be patched using formal patching techniques to construct a Galois cover $\Tilde {\Cal Y}_1\to \Cal P_i$ 
with Galois group $G$ between semi-stable $R$-curves (cf. [Ga], and Proposition 1.2.2), and by construction
the arithmetic genera $g_1$ and $g$ of the special fibres $\Tilde {\Cal Y}_{1,k}$ and $\Cal Y'_k$  
satisfy the inequality $g_1<g$. Indeed, we have eliminated 
the contribution to the arithmetic genus of $\Cal Y'_k$ which arise from the separable end components of
$\Gamma$, that lie above the end components of $\Gamma''$ that we encounter when moving form the vertex $P_i$ towards the ends of $\Gamma''$, and
which exist by the assumption (ii) satisfied by $P_i$. 
This proves that $P_i$ is a removable vertex as claimed.
\qed
\enddemo

$$\text{References.}$$

\noindent
[Ab] Abbes, A. \'El\'ements de g\'eom\'etrie rigide I. Book in preparation.

\noindent
[Ab1] Abbes, A. R\'eduction semi-stable des courbes d'apr\`es Artin, Deligne,
Grothendieck, Mumford, Saito, Winters. Courbes semi-stables et groupe fundamental
en g\'eom\'etrie alg\'ebrique (Luminy, 1998), 59-110, Progr. Math., 187, Birkh\"auser, Basel,
2000.

\noindent
[De-Mu] Deligne, P., and Mumford, D. The irreducibility of the space of curves of given genus.
Inst. Hautes \'Etudes Sci. Publ. Math. No, 36 (1969) 75-109.

\noindent
[Ga] Garuti, M. Prolongements de rev\^etements galoisiens en g\'eom\'etrie 
rigide, Compositio Mathematica, tome 104, n 3 (1996), 305-331.

\noindent
[Gr-Ma] Green, B., and Matignon, M. Liftings of Galois covers of smooth curves. Compositio Mathematica 113, 
237-272, 1998.

\noindent
[Gr-Ma1] Green, B., and Matignon, M. Order $p$-automorphisms of the open disc of a $p$-adic field, J. Amer. Math.
Soc. 12 (1) (1999) 269-303.

\noindent
[Gr] Grothendieck, A. Rev\^etements \'etales et groupe fondamental, Lecture 
Notes in Math. 224, Springer, Heidelberg, 1971.

\noindent
[Ha] Harbater, D. Moduli of $p$-covers of curves. Comm. Algebra 8
(1980), no. 12, 1095-1122.

\noindent
[Ka] Katz, N. Local-to-global extensions of representations of fundamental groups. 
Ann. Inst. Fourier (Grenoble) 36 (1986), no. 4, 69-106.

\noindent
[Oo] Oort, F. Lifting algebraic curves, abelian varieties, and their endomorphisms to characteristic zero, Proceedings of Symposia in 
Pure Mathematics, Vol. 46 (1987).

\noindent
[Oo1] Oort, F. Some Questions in Algebraic Geometry, Utrecht Univ. Math. Dept. Preprint
Series, June 1995. Available in the Homepage of F. Oort, 
\newline
http://www.staff.science.uu.nl/~oort0109/.

\noindent
[Ra] Raynaud, M. $p$-groupes et r\'eduction semi-stable des courbes. The Grotendieck
Festchrift Vol. III, 179-197, Progr. Math., 88, Birk\"auser Boston, Boston, MA, 1990.

\noindent
[Ri-Za] Ribes, L., and Zalesskii, P. Profinite groups, Ergenisse der Mathematik 
und ihrer Grenzgebiete, Folge 3,  Volume 40.

\noindent
[Ro] Roquette, P. Absch\"atzung der Automorphismenzhal von Functionenk\"orper
bei Primzhalcharacteristic, Math. Z. 117 (1970), 157-163.

\noindent
[Sa] Sa\"\i di, M. Wild ramification and a vanishing cycles formula, J. Algebra 273 (2004), no. 1, 
108-128.

\noindent
[Se-Oo-Su] Sekiguchi, T., Oort, F., and Suwa, N. On the deformation of Artin-Schreier to
Kummer. Ann. Sci. \'Ecole Norm. Sup. (4), 22, (1989), no. 3, 345-375.

\noindent
[Se] Serre, J-P. Cohomologie Galoisienne, Lecture Notes in Math., 5, Springer Verlag, 
Berlin, 1994.

\bigskip

\noindent
Mohamed Sa\"\i di

\noindent
College of Engineering, Mathematics, and Physical Sciences

\noindent
University of Exeter

\noindent
Harrison Building

\noindent
North Park Road

\noindent
EXETER EX4 4QF 

\noindent
United Kingdom

\noindent
M.Saidi\@exeter.ac.uk

\end
\enddocument